\documentclass[amssymb,12pt]{amsart}
\usepackage{latexsym}
\newtheorem{propo}{Proposition}[section]

\newtheorem{defi}[propo]{Definition}
\newtheorem{conje}[propo]{Conjecture}
\newtheorem{lemma}[propo]{Lemma}
\newtheorem{corol}[propo]{Corollary}

\newtheorem{theor}[propo]{Theorem}

\newtheorem{examp}[propo]{Example}
\newtheorem{remar}[propo]{Remark}

\newcommand{\Ker}{\operatorname{Ker}}
\newcommand{\Aut}{{\mathrm {Aut}}}
\newcommand{\Out}{{\mathrm {Out}}}
\newcommand{\Mult}{{\mathrm {Mult}}}
\newcommand{\Inn}{{\mathrm {Inn}}}
\newcommand{\Irr}{{\mathrm {Irr}}}

\newcommand{\diag}{{\mathrm {diag}}}
\newcommand{\soc}{{\mathrm {soc}}}
\newcommand{\End}{{\mathrm {End}}}

\newcommand{\Hom}{{\mathrm {Hom}}}

\def\rank{\mathop{\mathrm{ rank}}\nolimits}

\newcommand{\Spec}{{\mathrm {Spec}}}

\newcommand{\Sym}{{\mathrm {Sym}}}

\newcommand{\codim}{{\mathrm {codim}}}

\newcommand{\Res}{{\mathrm {Res}}}
\newcommand{\Ann}{{\mathrm {Ann}}}
\newcommand{\CC}{{\mathbb C}}

\newcommand{\QQ}{{\mathbb Q}}
\newcommand{\ZZ}{{\mathbb Z}}

\newcommand{\SSS}{{\mathbb S}}
\newcommand{\AAA}{{\mathbb A}}
\newcommand{\HA}{\hat{{\mathbb A}}}
\newcommand{\FF}{{\mathbb F}}

\newcommand{\FQ}{\mathbb{F}_{q}}

\newcommand{\EC}{\mathcal{E}}
\newcommand{\XC}{\mathcal{G}}
\newcommand{\HC}{\mathcal{H}}
\newcommand{\AC}{\mathcal{A}}
\newcommand{\OC}{\mathcal{O}}
\newcommand{\OL}{O^{\ell'}}
\newcommand{\IL}{\mathcal{I}}
\newcommand{\CL}{\mathcal{C}}
\newcommand{\RL}{\mathcal{R}}
\newcommand{\SL}{\mathcal{S}}
\newcommand{\DL}{\mathcal{D}}
\newcommand{\LF}{{\mathfrak{M}}}
\newcommand{\BS}{\bar{S}}
\newcommand{\HS}{\hat{\SSS}}
\newcommand{\VC}{V_{\CC}}
\newcommand{\VD}{V^{\otimes 4}}
\newcommand{\VE}{V^{\otimes e}}
\newcommand{\Vbb}{V^{(2)}}
\newcommand{\Vcc}{V^{(3)}}
\newcommand{\Vdd}{V^{(4)}}

\newcommand{\eps}{\epsilon}
\newcommand{\lam}{\lambda}
\newcommand{\al}{\alpha}
\newcommand{\gam}{\gamma}
\newcommand{\om}{\varpi}
\newcommand{\Om}{\Omega}
\newcommand{\Omp}{\Om^{+}}

\newcommand{\hra}{\hookrightarrow}
\newcommand{\ua}{\uparrow}
\newcommand{\da}{\downarrow}
\newcommand{\bb}{\mathfrak{b}}
\newcommand{\dd}{\mathfrak{d}}
\newcommand{\an}{\alpha_{n}}
\newcommand{\bn}{\beta_{n}}
\newcommand{\zn}{\zeta_{n}}
\newcommand{\ran}{\rho^{1}_{n}}
\newcommand{\rbn}{\rho^{2}_{n}}
\newcommand{\GC}{\mathcal{G}}
\newcommand{\GCC}{\mathcal{G}_{\CC}}
\newcommand{\GCR}{G^{\circ}}
\newcommand{\LC}{L_{\CC}}
\newcommand{\TW}{\tilde{\wedge}}

\newcommand{\TWB}{\tilde{\wedge}^{2}}
\newcommand{\TSB}{\tilde{S}^{2}}
\newcommand{\SB}{S^{2}}
\newcommand{\SC}{S^{3}}
\newcommand{\TWC}{\tilde{\wedge}^{3}}
\newcommand{\TSC}{\tilde{S}^{3}}
\newcommand{\SD}{S^{4}}
\newcommand{\WB}{\wedge^{2}}
\newcommand{\WC}{\wedge^{3}}
\newcommand{\WD}{\wedge^{4}}
\newcommand{\MK}{M_{2k}(\XC)}

\newcommand{\MB}{M_{4}(\XC)}
\newcommand{\MBA}{M'_{4}(\XC)}
\newcommand{\MC}{M_{6}(\XC)}
\newcommand{\MCA}{M'_{6}(\XC)}
\newcommand{\MD}{M_{8}(\XC)}
\newcommand{\MDA}{M'_{8}(\XC)}
\newcommand{\MDAA}{M''_{8}(\XC)}
\newcommand{\ME}{M_{10}(\XC)}
\newcommand{\MF}{M_{12}(\XC)}
\newcommand{\MG}{M_{14}(\XC)}
\newcommand{\qn}{q_{0}}
\newcommand{\dr}{\boldsymbol{d}_{V}}
\newcommand{\ags}{\boldsymbol{\alpha}}
\newcommand{\sign}{{\bf{sgn}}}
\def\skipa{\vspace{-1.5mm} & \vspace{-1.5mm} & \vspace{-1.5mm}\\}
\newcommand{\tn}{\hspace{0.5mm}^{t}\hspace*{-0.2mm}}
\newcommand{\ta}{\hspace{0.5mm}^{2}\hspace*{-0.2mm}}
\newcommand{\tb}{\hspace{0.5mm}^{3}\hspace*{-0.2mm}}
\renewcommand{\mod}{\bmod \,}

\begin{document}
\title[Small Tensor Powers and Larsen's Conjecture]
{Decompositions of Small Tensor Powers and Larsen's Conjecture}
\author{Robert M. Guralnick}
\address{Department of Mathematics, University of Southern California,
Los Angeles, CA 90089-1113, USA}
\email{guralnic@math.usc.edu}
\author{Pham Huu Tiep}
\address{Department of Mathematics, University of Florida, Gainesville,
FL 32611, USA}
\email{tiep@math.ufl.edu}
\date{submitted Feb. 27, 2003; revised Dec. 15, 2004}



\thanks{The authors would like to thank Nick Katz for suggesting this problem to them and 
for insightful comments on the paper, Nolan Wallach and Peter Sin for discussions about the 
decompositions for various tensor powers of the natural module for classical groups. They
are thankful to Arjeh Cohen, Gerhard Hiss, and Frank L\"ubeck for helping them with 
some computer calculation. They also thank the referee for helpful comments that greatly 
improved the exposition of the paper.}

\thanks{The authors gratefully acknowledge the support of the NSF (grants
DMS-0236185 and DMS-0070647), and of the NSA (grant H98230-04-0066).}

\begin{abstract}
{We classify all pairs $(G,V)$ with $G$ a closed subgroup
in a classical group $\GC$ with natural module $V$ over $\CC$, such that
$\GC$ and $G$ have the same composition factors on
$V^{\otimes k}$ for a fixed $k \in \{2,3,4\}$. In particular,
we prove Larsen's conjecture stating that for $\dim(V) > 6$ and $k = 4$
there are no such $G$ aside from those containing the derived subgroup of
$\GC$. We also find all the examples where this fails for $\dim(V) \leq 6$.
As a consequence of our results, we obtain a short proof of a related
conjecture of Katz. These conjectures are used in Katz's recent works on
monodromy groups attached to Lefschetz pencils and to character sums over
finite fields. Modular versions of these conjectures are also studied,
with a particular application to random generation in finite groups
of Lie type.}
\end{abstract}
\maketitle

\section{Introduction: Larsen's Conjecture}

Let $V = \CC^{d}$ with $d > 4$. Fix a nondegenerate quadratic form and if $d$ is even fix a 
nondegenerate symplectic form on $V$, and let $\XC$ be one of $GL(V)$, $O(V)$ or $Sp(V)$ (the 
latter only when $n$ is even).
If $G$ is any subgroup of $GL(V)$, define $M_{2k}(G,V)$ to be the dimension of
$\End_{G}(V^{\otimes k})$, and let $\GCR$ denote the connected component of $G$. Abusing the language
we will also say that the trivial group is reductive. 

M. Larsen \cite{Lars} proved the following alternative:

\begin{theor}\label{closed}
{\sl Let $G$ be a Zariski closed subgroup of $\XC$. Assume that $\GCR$ is a reductive 
subgroup of positive dimension. Then $M_{4}(G,V) > M_{4}(\XC,V)$ unless $G$ contains $[\XC,\XC]$.
\hfill $\Box$}
\end{theor}

An easy proof of Theorem \ref{closed} can be obtained as follows (see \cite{Ka1}). Let $G$ be as
in the theorem and $M_{4}(\XC,V) = M_{4}(G,V)$. By Weyl's theorem, every finite dimensional 
complex representation $W$ of $G$ is completely reducible. Hence the equality 
$M_{4}(\XC,V) = M_{4}(G,V)$ implies that the adjoint module of $\XC$ must be irreducible over $G$. 
However the Lie algebra $L(G^{\circ})$ is a nonzero invariant submodule of $L(\XC^{\circ})$. Thus, 
$L(G^{\circ}) = L(\GC^{\circ})$ and so $G \geq [\GC,\GC]$.

\begin{propo}\label{reduction}
{\sl Suppose that $G$ is a finite subgroup of $\XC$ such that $M_{4}(G,V) = M_{4}(\XC,V)$.
Then $G$ preserves no tensor structure and acts primitively on $V$.  In particular,
either $G/Z(G)$ is almost simple or $G$ is contained in the normalizer of an
irreducible subgroup of symplectic type of $\XC$ (in the latter case, $\dim(V) = p^{a}$
for some prime $p$).
\hfill $\Box$}
\end{propo}

Indeed, using Aschbacher's theorem \cite{A}, it is quite easy to see that any finite subgroup not of 
the form described either preserves a tensor structure or is imprimitive. But any imprimitive 
non-monomial subgroup or any subgroup preserving a tensor structure is contained in a nontrivial 
positive dimensional subgroup. Furthermore, any monomial subgroup in $SL(V)$ normalizes a 
$(d-1)$-dimensional torus of $SL(V)$ and so has an invariant $(d-1)$-dimensional subspace in the  
adjoint module (if $d := \dim(V)$). Also see \S4 for more refined reductions.

We are interested in Larsen's conjecture:

\begin{conje}\label{larsen}
{\sl If $G$ is a finite subgroup of $\XC$, then $M_{2k}(G,V) > M_{2k}(\XC,V)$ for some
$k \leq 4$.}
\end{conje}

An interesting application of this conjecture comes from algebraic geometry \cite{Ka1}. Given
a projective smooth variety $X$ of dimension $n+1 \geq 1$ over a finite field $k$, and 
consider a Lefschetz pencil of smooth hypersurface sections of a given degree $d$. Then the monodromy 
group $G_{d}$ is defined to be the Zariski closure of the monodromy group of the local system 
${\mathcal F}_{d}$ on the space of all smooth degree $d$ hypersurface sections (see \cite{Ka1}). Fix 
a degree $d$ hypersurface $H$ which is tranverse to $X$, and let $V$ be the subspace spanned by the 
vanishing cycles in $H^{n}((X \otimes_{k}\bar{k}) \cap H,\overline{\QQ}_{\ell})$ \cite[(4.2.4)]{D}. 
Then the cup product induces a ($G_{d}$-invariant) nondegenerate bilinear form on $V$. Deligne
\cite[4.4]{D} showed that $G_{d} = Sp(V)$ if $n$ is odd; if $n \geq 2$ is even then 
$G_{d} = O(V)$ or $G_{d}$ is finite. One would like to be able to rule out the 
finite group possibility (under additional hypotheses).
Katz has shown in \cite{Ka1} that $G_{d}$ is a subgroup of $\GC := Sp(V)$ or $O(V)$ with
the same {\it fourth moment} as of $\GC$, that is, $M_{4}(G_{d},V) = M_{4}(\GC,V)$.
Larsen's alternative \ref{closed} says that in this case either $G_{d} \geq [\GC,\GC]$ (if 
$G_{d}$ is reductive) or $G_{d}$ is finite. Furthermore, Larsen's conjecture \ref{larsen} states 
that if $\GC$ is a classical Lie group on $V$ and $G < \GC$ is a closed reductive subgroup with the 
same eighth moment as of $\GC$, i.e. $M_{8}(G,V) = M_{8}(\GC,V)$, then $G \geq [\GC,\GC]$.

A more recent application is described in \cite{Ka2}, where Larsen's alternative and 
Larsen's conjecture, as well as {\it drop ratio conjectures} (cf. \cite[Chapter 2]{Ka2} and 
Theorem \ref{main4} below) play an important role in the determination of the geometric monodromy 
group attached to a family of character sums over finite fiels.
 
In fact, we will prove:

\begin{theor}\label{main1}
{\sl Let $V = \CC^{d}$ with $d \geq 5$ and $\GC$ be $GL(V)$, $Sp(V)$, or $O(V)$. Assume that $G$ is a 
closed subgroup of $\GC$ such that $\GCR$ is reductive. Then one of the following holds:

{\rm (i)} $M_{8}(G,V) > M_{8}(\GC,V)$;

{\rm (ii)} $G \geq [\GC,\GC]$;

{\rm (iii)} $d = 6$, $\GC = Sp(V)$, and $G = 2J_{2}$.}
\end{theor}

Notice that Theorem \ref{main1} fails if $G$ is not reductive: if $G$ contains a Borel subgroup $B$ 
of $\GC$ then $M_{2k}(G,V) = M_{2k}(\GC,V)$ for any $k$. (Since $\GC$ is completely reducible on its 
finite dimensional representations, it suffices to show that $B$ has no fixed points on any nontrivial 
finite dimensional $\GC$-irreducible complex module. Now if $0 \neq v \in W$ is fixed by $B$ 
then $v$ is a maximal vector, $W = \CC\GC v$, and the highest weight of $W$ is $0$, cf. 
\cite[p. 189, 190]{Hu}. It follows that the module $W$ is trivial.) 

In the case (iii) of Theorem \ref{main1}, $M_{2k}(G,V) = M_{2k}(\GC,V)$ for $k \leq 5$ but 
$M_{12}(G,V) > M_{12}(\GC,V)$ (the values for $M_{2k}(\GC,V)$ are given in Lemma \ref{dim}). On
the other hand, if $d > 6$ then the condition $M_{8}(G,V) = M_{8}(\XC,V)$ for (reductive) subgroups 
$G \leq \GC$ implies that $G \geq [\GC,\GC]$. Moreover, Theorem \ref{main2} (below) states that the 
only closed subgroups of $\GC$ with $d > 6$ that can be irreducible on every $\GC$-composition 
factor of $V^{\otimes 4}$ are the ones containing $[\GC,\GC]$. A related question is asked by 
B. H. Gross: which finite subgroups of a simple complex Lie group $\GC$ is irreducible on every 
fundamental representation of $\GC$; this question has been answered in \cite{MMT}.

Our next Theorem \ref{complex}, resp. Theorem \ref{complex-m6}, classifies all (reductive) closed 
subgroups $G$ of $\GC$ such that $G$ has the same fourth, resp, sixth, moment as of $\GC$. 

\begin{theor}\label{complex}
{\sl Let $V = \CC^{d}$ with $d \geq 5$, $\GC = GL(V)$, $Sp(V)$, or $O(V)$. Assume $G$ is a
closed subgroup of $\GC$. Set $\BS = S/Z(S)$ for $S := F^{*}(G)$ if $G$ is finite. Then $G$ is 
irreducible on every $\GC$-composition factor of $V \otimes V^{*}$ (this condition is equivalent to  
$M_{4}(G,V) = M_{4}(\GC,V)$ if $G^{\circ}$ is reductive) if and only if one of the following holds.

{\rm (A)} $G \geq [\GC,\GC]$.

{\rm (B)} (Lie-type case) One of the following holds.

\hspace{0.5cm}{\rm (i)} $\BS = PSp_{2n}(q)$, $n \geq 2$, $q = 3,5$, $G = Z(G)S$, and
$V \da_{S}$ is a Weil module of dimension $(q^{n} \pm 1)/2$.

\hspace{0.5cm}{\rm (ii)} $\BS = U_{n}(2)$, $n \geq 4$, and $V \da_{S}$ is a
Weil module of dimension $(2^{n} + 2(-1)^{n})/3$ or $(2^{n}-(-1)^{n})/3$.

{\rm (C)} (Extraspecial case) $d = p^{a}$ for some prime $p$, $p > 2$ if $\GC = GL(V)$ and 
$p = 2$ otherwise, $F^{*}(G) = Z(G)E$ for some extraspecial subgroup $E$ of order $p^{1+2a}$ of 
$\GC$, and one of the conclusions {\rm (i) -- (iii)} of Lemma \ref{ex-irr} holds.

{\rm (D)} (Exceptional cases) $(\dim(V), \BS, G, \GC)$ is as listed in Table {\rm I}.}
\end{theor}

\begin{theor}\label{complex-m6}
{\sl Let $V = \CC^{d}$ with $d \geq 5$, $\GC = GL(V)$, $Sp(V)$, or $O(V)$. Assume $G$ is a
closed subgroup of $\GC$. Then $G$ is irreducible on every $\GC$-composition factor of 
$V^{\otimes 3}$ (this condition is equivalent to $M_{6}(G,V) = M_{6}(\GC,V)$ if $G^{\circ}$ is 
reductive) if and only if one of the following holds.

{\rm (A)} $G \geq [\GC,\GC]$; moreover, $G \neq SO(V)$ if $d = 6$.

{\rm (B)} (Extraspecial case) $d = 2^{a}$ for some $a > 2$. If $\GC = GL(V)$ then 
$G = Z(G)E \cdot Sp_{2a}(2)$ with $E = 2^{1+2a}_{+}$. If $\GC = Sp(V)$, resp. $O(V)$, then 
$E \cdot \Omega^{\eps}_{2a}(2) \leq G \leq E \cdot O^{\eps}_{2a}(2)$, with 
$E = 2^{1+2a}_{\eps}$ and $\eps = -$, resp. $\eps = +$.

{\rm (C)} (Exceptional cases) $G$ is finite, with the unique nonabelian composition factor 
$$\BS \in \{L_{3}(4),U_{3}(3),U_{4}(3),J_{2},\AAA_{9},\Omega^{+}_{8}(2),U_{5}(2),G_{2}(4),Suz,
J_{3},Co_{2},Co_{1},F_{4}(2)\},$$ 
and $(\dim(V), \BS, G, \GC)$ is as listed in the lines marked by $^{(\star)}$ in Table {\rm I}.}
\end{theor}

\newpage
\vskip10pt
\centerline
{{\sc Table} I. Exceptional examples in the complex case in dimension $d \geq 5$}
\begin{center}
\begin{tabular}{|r||c|c|c|c|c|} \hline \skipa
     $d$ & $\BS$ & $G$ & $\GC$ & 
     {\small \begin{tabular}{c}The largest $2k$ with\\ $M_{2k}(G,V) = M_{2k}(\GC,V)$ \end{tabular}}& 
         \begin{tabular}{c}$M_{2k+2}(G,V)$ vs.\\$M_{2k+2}(\GC,V)$ \end{tabular}\\ 
          \skipa \hline \hline \skipa
     $6$ & $\AAA_{7}$ & $6\AAA_{7}$ & $GL_{6}$ & $4$ & $21$ vs. $6$\\
     $6$ & $L_{3}(4)$$^{\;(\star)}$ & $6L_{3}(4) \cdot 2_{1}$ & $GL_{6}$ & $6$ & $56$ vs. $24$\\
     $6$ & $U_{3}(3)$$^{\;(\star)}$ & $(2 \times U_{3}(3)) \cdot 2$ & $Sp_{6}$ & $6$ & $195$ vs. $104$\\
     $6$ & $U_{4}(3)$$^{\;(\star)}$ & $6_{1} \cdot U_{4}(3)$ & $GL_{6}$ & $6$ & $25$ vs. $24$\\
     $6$ & $J_{2}$$^{\;(\star)}$ & $2J_{2}$ & $Sp_{6}$ & $10$ & $10660$ vs. $9449$\\ \hline
     $7$ & $Sp_{6}(2)$ & $Sp_{6}(2)$ & $O_{7}$ & $4$ & $16$ vs. $15$\\ \hline
     $8$ & $L_{3}(4)$ & $4_{1} \cdot L_{3}(4)$ & $GL_{8}$ & $4$ & $17$ vs. $6$\\
     $8$ & $\AAA_{9}$$^{\;(\star)}$ & $\HA_{9}$ & $SO_{8}$ & $6$ & $191$ vs. $106$\\ 
     $8$ & $\Om^{+}_{8}(2)$$^{\;(\star)}$ & $2\Om^{+}_{8}(2)$ & $SO_{8}$ & $6$ & $107$ vs. $106$\\\hline
     $10$ & $U_{5}(2)$$^{\;(\star)}$ & $(2 \times U_{5}(2)) \cdot 2$ & $Sp_{10}$ & $6$ & $120$ vs. 
        $105$\\
     $10$ & $M_{12}$ & $2M_{12}$ & $GL_{10}$ & $4$ & $15$ vs. $6$\\ 
     $10$ & $M_{22}$ & $2M_{22}$ & $GL_{10}$ & $4$ & $7$ vs. $6$\\ \hline
     $12$ & $G_{2}(4)$$^{\;(\star)}$ & $2G_{2}(4) \cdot 2$ & $Sp_{12}$ & $6$ & $119$ vs. $105$\\
     $12$ & $Suz$$^{\;(\star)}$ & $6Suz$ & $GL_{12}$ & $6$ & $25$ vs. $24$\\ \hline
     $14$ & $\ta B_{2}(8)$ & $\ta B_{2}(8) \cdot 3$ & $GL_{14}$ & $4$ & $90$ vs. $6$\\
     $14$ & $G_{2}(3)$ & $G_{2}(3)$ & $O_{14}$ & $4$ & $21$ vs. $15$\\ \hline
     $18$ & $Sp_{4}(4)$ & $(2 \times Sp_{4}(4)) \cdot 4$ & $O_{18}$ & $4$ & $25$ vs. $15$\\ 
     $18$ & $J_{3}$$^{\;(\star)}$ & $3J_{3}$ & $GL_{18}$ & $6$ & $238$ vs. $24$\\ \hline
     $22$ & $McL$ & $McL$ & $O_{22}$ & $4$ & $17$ vs. $15$\\ \hline
     $23$ & $Co_{3}$ & $Co_{3}$ & $O_{23}$ & $4$ & $16$ vs. $15$\\ 
     $23$ & $Co_{2}$$^{\;(\star)}$ & $Co_{2}$ & $O_{23}$ & $6$ & $107$ vs. $105$\\ \hline 
     $24$ & $Co_{1}$$^{\;(\star)}$ & $2Co_{1}$ & $O_{24}$ & $6$ & $106$ vs. $105$\\ \hline 
     $26$ & $\ta F_{4}(2)'$ & $\ta F_{4}(2)'$ & $GL_{26}$ & $4$ & $26$ vs. $6$\\ \hline
     $28$ & $Ru$ & $2Ru$ & $GL_{28}$ & $4$ & $7$ vs. $6$\\ \hline  
     $45$ & $M_{23}$ & $M_{23}$ & $GL_{45}$ & $4$ & $817$ vs. $6$\\ 
     $45$ & $M_{24}$ & $M_{24}$ & $GL_{45}$ & $4$ & $42$ vs. $6$\\ \hline
     $52$ & $F_{4}(2)$$^{\;(\star)}$ & $2F_{4}(2) \cdot 2$ & $O_{52}$ & $6$ & $120$ vs. $105$\\ \hline
     $78$ & $Fi_{22}$ & $Fi_{22}$ & $O_{78}$ & $4$ & $21$ vs. $15$\\ \hline
     $133$ & $HN$ & $HN$ & $O_{133}$ & $4$ & $21$ vs. $15$\\ \hline
     $248$ & $Th$ & $Th$ & $O_{248}$ & $4$ & $20$ vs. $15$\\ \hline  
     $342$ & $O'N$ & $3O'N$ & $GL_{342}$ & $4$ & $3480$ vs. $6$\\ \hline
     $1333$ & $J_{4}$ & $J_{4}$ & $GL_{1333}$ & $4$ & $8$ vs. $6$\\ \hline
\end{tabular}
\end{center}

Notice that in the cases (B) and (C) of Theorem \ref{complex}, $M_{8}(G,V) > M_{8}(\GC,V)$; 
see Propositions \ref{ex-m8}, \ref{sp-odd}, \ref{sp-m6}, \ref{su-m8}, \ref{su-m6}, and Lemma 
\ref{sp-m4} for estimates of $M_{2k}(G,V)$ with $2 \leq k \leq 4$ for various groups occuring in 
these cases. In the case (D), Table I lists the largest $2k$ for which $M_{2k}(G,V) = M_{2k}(\GC,V)$ 
as well as the exact value of $M_{2k+2}(G,V)$. 

\medskip
Our strategy can be outlined as follows. Suppose $G$ is a closed subgroup of $\GC$ such that 
$G$ is irreducible on every $\GC$-composition factor of $V \otimes V^{*}$, where 
$\GC \in \{ GL(V), Sp(V), O(V)\}$. Then \S4 yields basic reductions to the case where $G$ is 
a finite group. In fact either $G$ is contained in the normalizer of certain subgroups of
symplectic type of $\GC$, or $\BS \lhd G/Z(G) \leq \Aut(\BS)$ for some finite simple group $\BS$. 
The former case is handled in \S5. In the treatment of the latter case, Proposition \ref{bound1}
plays a crucial role; in particular, it shows that $\dim(V)$ is bounded as in the key 
inequality (\ref{bound3}). Thus $V$ is a representation of small degree of $G$. Classification
results on low-dimensional representations of finite quasi-simple groups \cite{GMST}, \cite{GT1},
\cite{HM}, \cite{Lu1}, \cite{Lu2}, then allow us to identify $V$. 

Katz \cite{Ka2} defines the {\it projective drop} $\dr(g)$ of any element $g \in GL(V)$ to be 
the smallest codimension (in $V$) of $g$-eigenspaces on $V$. In \cite{Ka2},
Katz has formulated three conjectures on the projective drop, Conjectures (2.7.1), (2.7.4), and
(2.7.7), in increasing order of strength. We will show that the second strongest, Conjecture 
(2.7.4) of \cite{Ka2}, holds true.

\begin{theor}\label{main4}
{\sl Let $V = \CC^{d}$ with $d \geq 2$, $\GC = GL(V)$, $Sp(V)$, or $O(V)$. Assume $G$ is a finite 
subgroup of $\GC$ such that $M_{4}(G,V) = M_{4}(\GC,V)$. Then 
$$\min\left\{\dfrac{\dr(g)}{\dim(V)} \mid g \in G \setminus Z(\GC)\right\} \geq 1/8~.$$ 
Moreover, the equality occurs precisely when $G$ is the Weyl group $W(E_{8})$ of type $E_{8}$ on its
($8$-dimensional) reflection representation.}
\end{theor}

Observe that, without the condition $\MB$, the drop ratio $\dr(g)/\dim(G)$
can get as close to $0$ as one wishes (just look at the irreducible
complex representations of degree $n-1$ of $\SSS_{n}$). A slightly different version of Theorem 
\ref{main4} has also been proved by Gluck and Magaard (private communication).
   
\medskip
The {\it Weil representations} of finite symplectic and unitary groups 
provide the source for many of the examples listed in Theorem \ref{complex}.
The complex Weil modules for the symplectic groups $S = Sp_{2n}(q)$ with $q = p^{f}$ and 
$p$ an odd prime are constructed as follows. Let $E$ be an extraspecial group of order $pq^{2n}$
of exponent $p$ (i.e. $[E,E] = \Phi(E) = Z(E)$ has order $p$). For each nontrivial 
linear character $\chi$ of $Z(E)$, the group $E$ has a unique irreducible module $M$ 
of dimension $q^{n}$ over $\CC$ that affords the $Z(E)$-character $q^{n}\chi$. Now $S$ acts 
faithfully on $E$ and trivially on $Z(E)$, and one can extend $M$ to the semidirect product $ES$. If we 
restrict $M$ to $S$, then $M = [t,M] \oplus C_{M}(t)$ where $t$ is the central 
involution in $S$, and these two summands are irreducible modules. This construction gives 
two irreducible modules of each dimension $(p^{n} \pm 1)/2$. A similar but slightly more complicated 
construction \cite{S1} leads to the complex Weil modules of the special unitary groups $U := SU_{n}(q)$ 
(here $q$ may be even as well); there is one such a module of dimension $(q^{n}+q(-1)^{n})/(q+1)$ 
and $q$ such of dimension $(q^{n}-(-1)^{n})/(q+1)$. When $q$ is odd, among the latter modules there is 
exactly one self-dual module which we denote by $V_{0}$. The Weil modules over fields of positive 
characteristic $\ell$ are defined to be nontrivial irreducible constituents of reduction modulo 
$\ell$ of complex Weil modules. To determine the type of the classical group
$\GC$ on the Weil module $V$ that contains $S$ or $U$ acting on $V$, one needs to know the 
Frobenius-Schur indicator of the Weil module. We summarize this information in Table II. Weil modules
are studied in detail in \S\S \ref{weil-sp}, \ref{weil-su}.

\vskip10pt
\centerline
{{\sc Table} II. Types of complex Weil modules for finite symplectic and unitary groups}
\begin{center}
\begin{tabular}{|c||c|c|} \hline \skipa
     Group & $\dim(V)$ & $\GC$ \\ \skipa \hline \skipa
     $Sp_{2n}(q)$, $q \equiv 3 (\mod 4)$ & $(q^{n} \pm 1)/2$ & $GL(V)$ \\ \hline \skipa
     $Sp_{2n}(q)$, $q \equiv 1 (\mod 4)$ & $(q^{n} \pm 1)/2$ & 
       $\left\{ \begin{array}{c}Sp(V) \mbox{ if }2|\dim(V)\\ 
                  O(V) \mbox{ if }2 \not{|}\dim(V)\end{array}\right.$ \\ \hline \skipa
     $SU_{n}(q)$ & $(q^{n}+q(-1)^{n})/(q+1)$ & 
       $\left\{ \begin{array}{c}Sp(V) \mbox{ if }2 \not{|}n \\ 
                  O(V) \mbox{ if }2|n \end{array}\right.$  \\ \hline \skipa
     $SU_{n}(q)$, $q$ even & $(q^{n}-(-1)^{n})/(q+1)$ & $GL(V)$ \\ \hline \skipa
     $SU_{n}(q)$, $q$ odd, $V = V_{0}$ & $(q^{n} - (-1)^{n})/(q+1)$ & 
       $\left\{ \begin{array}{c}Sp(V) \mbox{ if }2|n \\ 
                  O(V) \mbox{ if }2 \not{|}n\end{array}\right.$  \\ 
     $SU_{n}(q)$, $q$ odd, $V \neq V_{0}$ & $(q^{n}-(-1)^{n})/(q+1)$ & $GL(V)$ \\ \hline
\end{tabular}
\end{center}

\bigskip
Notice that Larsen's conjecture for small rank $\leq 3$ groups is handled in Theorem \ref{small}. 
In fact we will study Conjecture \ref{larsen} and prove Theorems \ref{main1}, \ref{complex}
in the more general context of algebraic groups in any characteristic. This 
modular context, as well as basic notation used in the paper, is described in \S2 below.  

\section{Modular Analogue of Larsen's Conjecture}

Motivated by some other applications, we will consider a modular analogue of Larsen's
conjecture. Let us fix some notation first.
Throughout the paper, $\FF$ is an algebraically closed field of characteristic $\ell$.
The condition $\ell > b$ will mean that $\ell = 0$ or $\ell > b$. Let $V = \FF^{d}$ be
equipped with a nondegenerate bilinear symplectic, resp. orthogonal, form, and let $\XC$ be
$GL(V)$, or $Sp(V)$, resp. $O(V)$. The irreducible $\GC$-module with highest weight $\om$ will be 
denoted by $L(\om)$. We denote the fundamental weights of $\GC$ by $\om_{1}, \ldots ,\om_{n}$ with 
$L(\om_{1}) \simeq V$. If $\ell > 0$, then $V^{(\ell)}$ is the $\ell^{\mathrm {th}}$ Frobenius twist of 
$V$. If $V$ is a composition factor of an $\FF G$-module $U$, we will write $V \hra U$ and say by 
abusing language that $V$ enters $U$. ``Irreducible'' always means ``absolutely irreducible''. A 
finite group $G$ is {\it quasi-simple} if $G = [G,G]$ and $G/Z(G)$ is simple, and 
{\it nearly simple} if $\BS \lhd G/Z(G) \leq \Aut(\BS)$ for some finite simple non-abelian 
group $\BS$. For a finite group $G$, $F^{*}(G)$ is the generalized 
Fitting subgroup of $G$. Furthermore, a {\it component} of $G$ is a quasi-simple subnormal 
subgroup of $G$, and $E(G)$ is the subgroup generated by all components of $G$. For a prime $\ell$,
$O^{\ell}(G)$, resp. $\OL(G)$, is the smallest normal subgroup $N$ of $G$ such that
$G/N$ is an $\ell$-group, resp. an $\ell'$-group. We denote $PSL_{n}(q)$ by $L_{n}(q)$ and 
$PSU_{n}(q)$ by $U_{n}(q)$. In the paper, $G$ is said to have {\it symplectic type} if
$G = Z(G)E$ for some extraspecial $p$-subgroup $E$.    

\begin{defi}\label{m2k}
{\em Let $k \geq 1$. We say that a subgroup $G$ of $\XC$ satisfies the condition $\MK$,
if for any $j$ with $0 \leq j \leq k$ and for any two $\XC$-composition factors $Y$, $Y'$ of
the module $V^{\otimes (k-j)} \otimes (V^{*})^{\otimes j}$ the following holds:

(i) $G$ is irreducible on $Y$ and $Y'$; and

(ii) if the $\XC$-modules $Y$ and $Y'$ are nonisomorphic, then so are the $G$-modules
$Y$ and $Y'$.}
\end{defi}

The modular analogue of Larsen's conjecture we have in mind is that only very few specific
subgroups of $\GC$ can satisfy all the conditions $M_{2k}(\GC)$ for $k \leq 4$.

\begin{defi}\label{mp2k}
{\em We also introduce some more conditions for $G \leq \GC$:

$\MBA$: $G$ is irreducible on any $\GC$-composition factor of $V \otimes V^{*}$.

$\MCA$: $G$ is irreducible on the $\GC$-composition factor $L(3\om_{1})$ of $V^{\otimes 3}$.

$\MDA$: for any two $\GC$-composition factors $Y$, $Y'$ of the module $\VD$,
$G$ is irreducible on $Y$ and $Y'$, and moreover, if the $\GC$-modules $Y$ and $Y'$ are
nonisomorphic then so are the $G$-modules $Y$ and $Y'$.

$\MDAA$: $G$ is irreducible on any $\GC$-composition factor of $\VD$.}
\end{defi}

\begin{remar}
{\em (i) {\sl If $\FF = \CC$ and $\GCR$ is reductive, then the equality 
$M_{2k}(G,V) = M_{2k}(\XC,V)$ is equivalent to that $G$ satisfy the condition $\MK$.} Indeed,
for $0 \leq j \leq k$ we have
$$M_{2k}(G,V) = \dim\End_{G}(V^{\otimes k}) =
  \dim \Hom_{G}(V^{\otimes k} \otimes (V^{*})^{\otimes k},1_{G}) = $$
$$=  \dim \Hom_{G}(V^{\otimes (k-j)} \otimes (V^{*})^{\otimes j},
  V^{\otimes (k-j)} \otimes (V^{*})^{\otimes j}) =
  \dim \End_{G}(V^{\otimes (k-j)} \otimes (V^{*})^{\otimes j}).$$
Since any such a subgroup $G$ is completely reducible on finite dimensional $\GC$-modules, the 
condition $M_{2k}(G,V) = M_{2k}(\XC,V)$ now implies \ref{m2k}(i) and \ref{m2k}(ii) for any two 
$\XC$-composition factors $Y$, $Y'$ of $V^{\otimes (k-j)} \otimes (V^{*})^{\otimes j}$.

(ii) {\sl If $k \geq 3$ then $\MK$ implies $M_{2k-4}(\XC)$.} Indeed, for
$0 \leq j \leq k-2$, $V^{\otimes (k-2-j)} \otimes (V^{*})^{\otimes j}$ is a submodule of
$V^{\otimes (k-1-j)} \otimes (V^{*})^{\otimes (j+1)}$, since
the trivial module is a submodule of $V \otimes V^{*}$. Compare with Lemma \ref{m68}.

(iii) {\sl If $\dim(V) \geq 8$, then our results still hold true if we replace $\MDA$ by
$\MDAA$.} Indeed, if $\dim(V) \geq 8$ and $\GC = GL(V)$, then all distinct $\GC$-composition
factors are of distinct degree, cf. Proposition \ref{v4}. If $\GC = Sp(V)$ or $O(V)$, then in
all our arguments, we make use only of $\MDAA$.}
\end{remar}

\begin{remar}\label{choice}
{\em In this paper, we aim to find all subgroups $G$ of $\XC$ that satisfy the condition
$\MC \cap \MD$. Actually, our results (except for Theorem \ref{small}) still hold true 
if we replace $\MC \cap \MD$ by the weaker condition $\MBA \cap \MCA \cap \MDA$.}
\end{remar}

Before stating our main results, we display some examples clarifying the connections
between various conditions we formulated above.

\begin{examp}\label{clar1}
{\em (i) For each $\ell \in \{3,5,11,23\}$, $G = M_{23}$ has an irreducible $\FF G$-module of
dimension $d = 45$, where exactly one of $\AC(V)$, $\TSB(V)$, $\TWB(V)$ (see \S2.2 for
definition of these modules) is reducible over $G$ (and $V$ lifts to characteristic $0$.)
This phenomenon also happens to $\Aut(\ta B_{2}(8))$ in dimension $d = 14$ and characteristic
$\ell \in \{5,7,13\}$, and to $M_{22}$ in dimension $10$ and characteristic $2$. These
examples show in particular that in the case $\GC = GL(V)$ and $\ell > 0$, $\MB$ is strictly 
stronger than $\MBA$.

(ii) $G = 2M_{22}$ has a $10$-dimensional irreducible module $V$ in characteristic
$3$, where $G$ is irreducible on all $\GC$-composition factors of $V^{\otimes 3}$, but not
$V^{\otimes 2} \otimes V^{*}$. In particular, $\MC$ is strictly stronger than $\MCA$.

(iii) $G = 6_{1} \cdot U_{4}(3)$ has an irreducible module $V = \CC^{6}$,
where $G$ is irreducible on all but one $\GC$-composition factors of $V^{\otimes 4}$ (so $G$ almost
satisfies $\MD$ !).}
\end{examp}

\begin{examp}\label{clar2}
{\em (i) $S = Sp_{4}(4)$ has an irreducible $18$-dimensional module in characteristic
$\ell \neq 2,3$, where $\MB$ holds for some group $G = 2 \cdot \Aut(S)$ but not for any
other subgroup between $S$ and $G$, cf. Proposition \ref{sp44}. Thus $\MB$ may hold
for $G$ but fails for $F^{*}(G)$. This also happens to a $14$-dimensional complex
module of $\Aut(\ta B_{2}(8))$.

(ii) $S = 2 \cdot F_{4}(2)$ has an irreducible module $V = \CC^{52}$,
where $\MC$ holds for an extension $G = S \cdot 2$ (isoclinic to the one listed in
\cite{Atlas}) but not for $S$. Thus $\MC$ may hold for $G$ but fails for $F^{*}(G)$.
This also happens to an irreducible $6$-dimensional complex module of $U_{3}(3) \cdot 2$,
an irreducible $6$-dimensional module of $6 \cdot L_{3}(4) \cdot 2$ in characteristic $0$
or $3$, an irreducible $12$-dimensional module of $2 \cdot G_{2}(4) \cdot 2$
in characteristic $\neq 2,5$, as well as to an irreducible $\ell$-modular representation of
degree $10$ of $(2 \times U_{5}(2)) \cdot 2$ with $\ell \neq 2,3$.}
\end{examp}

\begin{examp}\label{clar3}
{\em $G = U_{4}(2)$ has an irreducible $5$-dimensional module $V$ in characteristic $3$.
This $V$ lifts to a complex module $\VC$. However, $\MC$ holds for $V$ but not for $\VC$.
A similar phenomenon happens to an irreducible $10$-dimensional $7$-modular representation
of $2M_{22}$.}
\end{examp}

\begin{examp}\label{clar4}
{\em Let $G = Sp_{4}(2) \simeq \SSS_{6}$ and $\ell = 2$. The natural representation of $G$ on
$V = \FF^{4}$ embeds $G$ in $\GC = Sp(V)$. All the $\GC$-composition factors on
$V^{\otimes 2}$ are $\Vbb$ and $\TWB(V) = L(\om_{2})$, which over $G$ is conjugate to $V$
via an outer automorphism of $G$. Hence $\MB$ holds for $G$. Furthermore, the
$\GC$-composition factors on $V^{\otimes 4}$ are $L(4\om_{1}) = \Vdd$,
$L(2\om_{1}+\om_{2}) = \Vbb \otimes L(\om_{2})$, $L(2\om_{2}) = L(\om_{2})^{(2)}$,
$L(2\om_{1}) = \Vbb$, $L(\om_{2})$, and $L(0)$, and all of them are irreducible over $G$. Thus
$\MD$ holds for $G$. However, $\MCA$ fails for $G$ as
$L(3\om_{1}) \da_{G} = (V \otimes \Vbb)\da_{G} =  (V \otimes V^{*})\da_{G}$ is reducible.
Consequently, $\MD$ does {\it not} imply $\MCA$ nor $\MC$. But see Remark \ref{m126}.}
\end{examp}

Define
$$\EC_{1} := \{\ta E_{6}(2), Ly, Th, Fi_{23}, Co_{1}, J_{4}, Fi'_{24}, BM, M\},~~~~
  \EC_{2} := \{F_{4}(2), Fi_{22}, HN\},$$
$$\EC_{3} := \{M_{11}, M_{12}, J_{1}, M_{22}, J_{2}, M_{22}, HS, J_{3}, M_{24}, McL, He,
  Suz, Co_{3}, Co_{2}, Ru, O'N,$$
$$L_{2}(7), L_{2}(11), L_{2}(13), L_{3}(4), U_{3}(3), U_{4}(2), U_{4}(3),
  U_{6}(2), Sp_{4}(4), Sp_{6}(2),$$
$$\Omega^{+}_{8}(2),\Omega_{7}(3),\ta B_{2}(8),\tb D_{4}(2),G_{2}(3),
  G_{2}(4),\ta F_{4}(2)',\AAA_{n} \mid 5 \leq n \leq 10\}.$$

One of our main results, which implies Theorem \ref{main1}, is

\begin{theor}\label{main2}
{\sl Let $\FF$ be an algebraically closed field of characteristic $\ell$,
$V = \FF^{d}$ with $d \geq 5$, $\GC = GL(V)$, $Sp(V)$, or $O(V)$. Assume $G$ is a
closed subgroup of $\GC$ that is irreducible on every $\GC$-composition factor of
$V \otimes V^{*}$ and of $\VD$, and on $L(3\om_{1})$. Then one of the following holds.

{\rm (A)} $G$ is of positive dimension, and one of the following holds.

\hspace{0.5cm}{\rm (i)} $G \geq [\GC,\GC]$;

\hspace{0.5cm}{\rm (ii)} $\GC = Sp(V)$, $\ell = 2$, and either $O(V) \geq G \geq \Omega(V)$ with 
$d \geq 10$ or $G = O(V)$ with $d = 8$.

{\rm (B)} $q \geq 4$ and one of the following holds.

\hspace{0.5cm}{\rm (i)} $\GC = GL(V)$, $\OL(G) = SL_{d}(q)$ or $SU_{d}(q)$, and $\ell | q$.

\hspace{0.5cm}{\rm (ii)} $\GC = Sp(V)$, $\OL(G) = Sp_{d}(q)$, and $\ell | q$.

\hspace{0.5cm}{\rm (iii)} $\GC = O(V)$, $\OL(G) = \Om^{\pm}_{d}(q)$, and $\ell | q$.

\hspace{0.5cm}{\rm (iv)} $\GC = Sp(V)$, $\ell = 2$, $\OL(G) = \Om^{\pm}_{d}(q)$, and
$2 | q$. Furthermore, either $d \geq 10$, or $d = 8$ and $G \geq O^{\pm}_{d}(q)$.

\hspace{0.5cm}{\rm (v)} $\GC = O(V)$, $d = 8$, $\OL(G) = \tb D_{4}(q)$, $\ell | q$.

{\rm (C)} $d = 6$, $\ell \neq 2,5$, $\GC = Sp(V)$, and $G = 2J_{2}$.

{\rm (D)} For $S = F^{*}(G)$, $S/Z(S) \in \EC_{1}$.\\
Conversely, all cases in {\rm (A), (B), (C)} give rise to examples.}
\end{theor}

Since none of the groups in Theorem \ref{main2}(D) can satisfy $\MF$
(see \S\ref{nearly3}), and since $\MD$ implies $\MB$ and $\ME$ implies $\MC$,
Theorem \ref{main2} immediately yields

\begin{corol}\label{m12}
{\sl Let $\FF$ be an algebraically closed field of characteristic $\ell$,
$V = \FF^{d}$ with $d \geq 5$, $\GC = GL(V)$, $Sp(V)$, or $O(V)$. Assume $G$ is a
closed subgroup of $\GC$ that satisfies both $\ME$ and $\MF$. Then one of the following holds.

{\rm (a)} $G$ is of positive dimension, and case {\rm (A)} of Theorem {\rm \ref{main2}} holds.

{\rm (b)} $G$ is finite, nearly simple, $E(G) \in Lie(\ell)$, and one of the conclusions 
{\rm (i) -- (v)} of Theorem {\rm \ref{main2}(B)} holds.
\hfill $\Box$}
\end{corol}

It is worthwhile to notice that, in the modular case, no condition of type
$\cap^{N}_{k=1}M_{2k}(\GC)$ would be strong enough to distinguish the algebraic group $\GC$
from its proper closed (or even finite) subgroups, cf. Lemma \ref{sp-o} for a glaring example.
Furthermore, even in the complex case and in big enough dimensions, the condition
$\MB \cap \MC$ is still not enough to tell apart $\GC = GL(V)$, $Sp(V)$, or $O(V)$ from its
finite subgroups, cf. Proposition \ref{ex-m6}.

Another main result is

\begin{theor}\label{main3}
{\sl Let $\FF$ be an algebraically closed field of characteristic $\ell$,
$V = \FF^{d}$ with $d \geq 5$, $\GC = GL(V)$, $Sp(V)$, or $O(V)$. Assume $G$ is a closed
subgroup of $\GC$ that is irreducible on every $\GC$-composition factor of
$V \otimes V^{*}$. Set $\BS = S/Z(S)$ for $S := F^{*}(G)$. Then one of the following holds.

{\rm (A)} (Positive dimensional case) $G$ is of positive dimension, and one of the following
holds.

\hspace{0.5cm}{\rm (i)} $G \geq [\GC,\GC]$.

\hspace{0.5cm}{\rm (ii)} $\GC = Sp(V)$, $\ell = 2$, and $O(V) \geq G \geq \Omega(V)$.

\hspace{0.5cm}{\rm (iii)} $\GC = Sp(V)$, $\ell = 2$, $d = 6$, and $G = G_{2}(\FF)$.

{\rm (B)} (Defining characteristic case) $\ell|q$ and one of the following holds.

\hspace{0.5cm}{\rm (i)} $\GC = GL(V)$, and $\OL(G) = SL_{d}(q)$ or $SU_{d}(q)$.

\hspace{0.5cm}{\rm (ii)} $\GC = Sp(V)$ and $\OL(G) = Sp_{d}(q)'$.

\hspace{0.5cm}{\rm (iii)} $\GC = Sp(V)$, $\ell = 2$, and $F^{*}(G) = \Om^{\pm}_{d}(q)$.

\hspace{0.5cm}{\rm (iv)} $\GC = O(V)$ and $F^{*}(G) = \Om^{\pm}_{d}(q)$.

\hspace{0.5cm}{\rm (v)} $\GC = O(V)$, $d = 8$, and $\OL(G) = \tb D_{4}(q)$.

\hspace{0.5cm}{\rm (vi)} $\GC = Sp(V)$, $\ell = 2$, $d = 6$, and $F^{*}(G) = G_{2}(q)'$.

\hspace{0.5cm}{\rm (vii)} $\GC = Sp(V)$, $\ell = 2$, $d = 8$, and $F^{*}(G) = \tb D_{4}(q)$.

{\rm (C)} (Alternating case) $S = \AAA_{n}$, and either $\ell|n$ and $d = n-2$, or $\ell = 2$ and 
$2|d = n-1$. Furthermore, $Z(G) \leq \ZZ_{2}$, $S \leq G/Z(G) \leq \Aut(S)$, and $V \da_{S}$ is 
labelled by $(n-1,1)$.

{\rm (D)} (Cross characteristic case) One of the following holds.

\hspace{0.5cm}{\rm (i)} $\BS = PSp_{2n}(q)$, $n \geq 2$, $\ell \neq q = 3,5$, $G = Z(G)S$, and
$V \da_{S}$ is a Weil module of dimension $(q^{n} \pm 1)/2$.

\hspace{0.5cm}{\rm (ii)} $\BS = U_{n}(2)$, $n \geq 4$, $\ell \neq 2$, and $V \da_{S}$ is a
Weil module of dimension $(2^{n} + 2(-1)^{n})/3$ or $(2^{n}-(-1)^{n})/3$.

{\rm (E)} (Extraspecial case) $F^{*}(G)$ is a group of symplectic type, and one of the
conclusions {\rm (i) -- (iii)} of Lemma \ref{ex-irr} holds.

{\rm (F)} (Small cases) If $\BS \in \EC_{3}$, then $V$ is as listed in Tables
{\rm III -- VI} located in \S8.

{\rm (G)} $\BS \in \EC_{1} \cup \EC_{2}$.\\
Conversely, all cases except possibly {\rm (G)} give rise to examples.}
\end{theor}

Notice that Theorems \ref{closed} and \ref{main3} immediately imply Theorem \ref{complex}.

Larsen's conjecture (and its modular analogue) for groups of small rank (i.e. of types 
$A_{1}$, $A_{2}$, $B_{2}$) is settled by the following theorem:
 
\begin{theor}\label{small}
{\sl Let $\FF$ be an algebraically closed field of characteristic $\ell$, $V = \FF^{d}$ with 
$d \leq 4$, $\GC \in \{SL(V), Sp(V)\}$ a simple simply connected algebraic group. Let $G$ be a 
proper closed subgroup of $\GC$. Then $G$ satisfies the condition $\MC \cap \MD$ if and only
if one of the following holds.

{\rm (i)} There is a Frobenius map $F$ on $\GC$ such that 
$\OL(G) = \GC^{F} \in \{SL_{d}(q), SU_{d}(q), Sp_{4}(q)\}$ with $\ell|q \geq 4$. Furthermore,
if $\OL(G) = SL_{3}(4)$ then $G = SL_{3}(4) \cdot 3$.

{\rm (ii)} $\GC = Sp_{4}(\FF)$, $\ell = 2$, and $G = \ta B_{2}(q)$ with $q \geq 32$.

{\rm (iii)} $\GC = SL_{2}(\FF)$, $\ell \neq 2,5$, and $G = SL_{2}(5)$.}
\end{theor}

Notice that, in the case (iii) of Theorem \ref{small}, $G$ satisfies $M_{2k}(\GC)$ for $k \leq 5$ (and 
also $M_{12}(\GC)$ if $\ell = 3$), but fails $M_{12}(\GC)$ if $\ell \neq 3$ and fails 
$M_{14}(\GC)$ if $\ell = 3$. Thus the groups $2J_{2}$ and $SL_{2}(5)$ are 
the only exceptions to Conjecture \ref{larsen} in the case of complex semisimple groups $\GC$. 

\medskip
Theorem \ref{main3} is related to the results proved in \cite{LiS, LST, Mag, MM, MMT, Mal}. Closed
subgroups of exceptional algebraic groups that are irreducible on a minimal or adjoint module have
been classified by Liebeck and Seitz \cite{LiS}. In the case $\ell = 0$, the finite nearly simple
subgroups of $\GC = GL(V)$ that are irreducible on the nontrivial $\GC$-composition factor of
$V \otimes V^{*}$ have been classified in \cite{Mal}. Finite quasisimple subgroups of $\GC$ that
are irreducible on the nontrivial $\GC$-composition factor of $V \otimes V^{*}$ if $V$ is not
self-dual, resp. of $\SB(V)$ or $\WB(V)$ if $V$ is self-dual, are considered in the other papers
mentioned above. Our methods are different from that of \cite{Mal} and similar to the ones
exploited in \cite{MMT}. Our assumption ($G$ is irreducible on every $\GC$-composition factor of
$V \otimes V^{*}$) is stronger than the ones used in the aforementioned papers in the case $V$ is
self-dual. On the other hand, we handle all closed (in particular finite) subgroups of $\GC$. Observe
that Example \ref{clar2} shows that the nearly simple case cannot be reduced to the
quasi-simple case: $\MB$ (or even $\MC$) holds for some $\FF G$-module $V$ but not for $V$
as a module over $F^{*}(G)$. Moreover, the modular case cannot be reduced to the complex
case even if the $\FF G$-module $V$ lifts to a complex module $\VC$: there are examples where $\MC$
holds for $V$ but not for $\VC$ even when $\ell = 7$, cf. Example \ref{clar3}. (But notice that 
if the $\FF G$-module $V$ lifts to a complex module $\VC$ and if $\ell \not{|} |G|$ then 
$M_{2k}(\GC)$ holds for $V$ if and only if it holds for $\VC$. Also see Remark \ref{v5} for
a partial reduction to the complex case.) 

\medskip
We expect Theorems \ref{main2} and \ref{main3} to have various applications. Below we point
out a particular consequence, in which $\overline{\langle X \rangle}$ denotes the closure
in $\GC$ of the subgroup generated by a subset $X$ of $\GC$. 

\begin{corol}\label{nice}
{\sl Let $\GC = GL(V)$, $Sp(V)$, or $O(V)$, where $V = \FF^{d}$ and $d \geq 5$. There is a
finite explicit list $\LF$ of irreducible $\FF\GC$-modules such that
$$\{(x,y) \in \GC \times \GC \mid \overline{\langle x,y \rangle}
  \mbox{ contains a subfield subgroup}\} = $$
$$ = \{(x,y) \in \GC \times \GC \mid \overline{\langle x,y \rangle} \mbox{ is irreducible on each }
  M \in \LF\}$$
is a non-empty open (in particular dense) subvariety of $\GC \times \GC$.
\hfill $\Box$}
\end{corol}

See \S\ref{generation} for more results in this direction. There is an analogous result for 
exceptional groups. Also, Theorem \ref{restr1} classifies all closed subgroups of a simple algebraic 
group $\GC$ that are irreducible on the adjoint module and the basic Steinberg module of $\GC$.

Corollary \ref{nice} is related to a conjecture of Dixon \cite{Di} stating that two random elements 
of a finite simple group $G$ generate $G$, with probability  $\to 1$ when $|G| \to \infty$. Dixon
himself proved the conjecture for alternating groups; the Lie type groups have been handled by 
Kantor and Lubotzky \cite{KaL} (classical groups), and Liebeck and Shalev \cite{LSh}) (exceptional 
groups). Random generation for semisimple algebraic groups over local fields has recently been addressed
in \cite{BL}. Corollary \ref{nice} gives another proof of random generation by a pair of elements for 
finite groups of Lie type of fixed type over $\FQ$ when $q \to \infty$. Roughly, the idea is as 
follows. Let $\GC$ be a simple algebraic group and consider $G := \GC^{F}$, the fixed point subgroup
under some Steinberg-Lang endomorphism $F$. The set of pairs in the finite group $G$ that generate an
irreducible subgroup on any finite collection of modules for $\GC$ are the $F$-fixed points of
an open subvariety of $\GC \times \GC$.  By results in arithmetic algebraic geometry, the proportion
of elements in $G \times G$ which are in the complement is at most $O(1/q)$ where $q$ is the size of
the field associated to $F$.  By the results above, the remaining elements must generate a subfield
subgroup. One then shows that the set of pairs of elements generating a proper subfield subgroup is
quite small, whence the result.  See \cite{gurnewton}, \cite{GLSS} and \cite{gurshalev} for
variations on this theme.

\medskip
The rest of the paper is organized as follows. In \S3 we prove some auxiliary results. Reduction 
theorems are proved in \S4 which reduce the problem to the cases where either $F^{*}(G)$ is of 
symplectic type, or $G$ is nearly simple. The case $F^{*}(G)$ is of symplectic type is completed in 
\S5. The nearly simple groups are dealt with in \S\S6 - 8. Main Theorems \ref{main1}, 
\ref{complex}, \ref{complex-m6}, \ref{main2}, \ref{main3}, and \ref{main4} are proved in \S9. 
Theorem \ref{small} is proved in \S\ref{smallrank}. 
Some results about generation of finite groups of Lie type; in particular, a proof of Corollary 
\ref{nice}, are exhibited in \S\ref{generation}. The labels for numbered formulae also include
the number of the section where the formula is located; for instance, formula (6.1) is the 
first numbered formula in \S6. Some calculation in the paper has been done using the package GAP
\cite{GAP}.

\section{Some Preliminary Results}

\subsection{Complex case}
In this subsection we gather some statements about representations of complex Lie groups.
In particular, $V = \CC^{d}$. We start with the following trivial observation.

\begin{lemma}\label{m68}
{\sl Suppose that $G_{1} \leq G_{2} \leq GL(V)$ and $M_{2k}(G_{1},V) > M_{2k}(G_{2},V)$. Then
$M_{2k+2}(G_{1},V) > M_{2k+2}(G_{2},V)$.}
\end{lemma}

\begin{proof}
Observe that the $\GC$-module $V \otimes V^{*} = A \oplus I$ contains the trivial submodule $I$ as
a direct summand. Let $W := V^{\otimes k}$. Then for any subgroup $X$ of $G_{2}$ we have
$$M_{2k+2}(X,V) = \dim\End_{X}(V \otimes W) =
  \dim \Hom_{X}(V \otimes V^{*}, W \otimes W^{*}) = $$
$$= \dim \Hom_{X}(I, W \otimes W^{*}) + \dim \Hom_{X}(A, W \otimes W^{*}) =
  M_{2k}(X,V) + f(X),$$
where $f(X) := \dim \Hom_{X}(A,W \otimes W^{*})$. Since
$f(G_{1}) \geq f(G_{2})$ we have $M_{2k+2}(G_{1},V) > M_{2k+2}(G_{2},V)$.
\end{proof}

Lemma \ref{m68} shows that when $\ell = 0$ the condition $M_{2k}(\GC)$ implies $M_{2j}(\GC)$
for any $j < k$.

The following result is well known in invariant theory (see \cite{Weyl}). Recall that
$(2k-1)!!$ is the product of all odd positive integers less than $2k$.

\begin{lemma}\label{dim}
{\sl Let $V$ be an $n$-dimensional vector space.

{\rm (i)} $M_{2k}(GL(V),V) = k!$  for $k = 1, 2, \ldots ,n$.

{\rm (ii)} $M_{2k}(Sp(V),V) = (2k-1)!!$ for $k = 1, 2, \ldots ,n$.

{\rm (iii)} $M_{2k}(SO(V),V) = (2k-1)!!$ for $2k < n$.
\hfill $\Box$}
\end{lemma}

We are thankful to N. Wallach for pointing out the following two results.

\begin{lemma}\label{nolan1}
{\sl Let $\AC$ be the adjoint module for $\XC = SL_{n}(\CC)$ with $n \geq 3$.  Then
$\dim \Hom_{\XC}(\AC, \AC \otimes \AC) = 2$.}
\end{lemma}

\begin{proof}
Let $V = \CC^{n}$ denote the natural $\XC$-module and $I$ the trivial module. It is well known that 
$W := V \otimes V^{*} = I \oplus \AC$, and 
$V^{\otimes 3} = \SC(V) \oplus \WC(V) \oplus 2 \cdot \SSS_{(2,1)}(V)$ (where the last summand is 
obtained by applying the Schur functor corresponding to the partition $(2,1)$ to $V$, cf. 
\cite[p. 78]{FH}). Hence 
$\dim \Hom_{\XC}(\AC, \AC \otimes \AC) = \dim \Hom_{\XC}(\AC^{\otimes 3},I)$ equals  
$$\dim \Hom_{\XC}(W^{\otimes 3},I) - 3 \cdot \dim \Hom_{\XC}(W \otimes W, I) + 
  3 \cdot \dim \Hom_{\XC}(W,I) - 1 = $$ 
$$  = \dim \Hom_{\XC}(V^{\otimes 3},V^{\otimes 3}) - 3 \cdot 2 + 3 - 1 = 2.$$     
\end{proof}

Indeed, Wallach has pointed out the following lovely formula (which we will not need):

\begin{lemma}
{\sl Let $\XC$ be a simple Lie group of rank $r$ and $\AC$ the adjoint module for $\XC$.  Then
$\dim \Hom_{\XC}(\AC,\AC \otimes \AC) = r - s$, where $s$ is the number of simple roots of $\XC$
that are perpendicular to the highest root.
\hfill $\Box$}
\end{lemma}

This leads to the following decomposition (and the answer can be easily verified).

\begin{lemma}\label{nolan2}
{\sl Let $\AC$ be the adjoint module for $\XC = SL_{n}(\CC)$ with $n \geq 4$. Let $\om_{i}$,
$1 \leq i < n$ be the fundamental weights for $\XC$. The composition factors of $\AC \otimes \AC$
are precisely the trivial module, $\AC$ (twice) and the modules

{\rm (i)} $ V(2\om_{1} + 2\om_{n-1})$, $\dim = n^{2}(n-1)(n+3)/4$;

{\rm (ii)} $ V(\om_{2} + 2\om_{n-1})$, $\dim =(n+2)(n+1)(n-1)(n-2)/4$;

{\rm (iii)} $ V(2\om_{1} + \om_{n-2})$, $\dim =(n+2)(n+1)(n-1)(n-2)/4$;

{\rm (iv)} $ V(\om_{2} + \om_{n-2})$, $\dim = n^{2}(n+1)(n-3)/4$.
\hfill $\Box$}
\end{lemma}

\subsection{Modular case}

First we introduce some more notation to be used throughout the paper.
Assume $\GC = GL(V)$ with $d := \dim(V) \geq 2$ and set $Y(V) := V \otimes V^{*}$. Since
$$\dim \Hom_{\GC}(Y(V),1_{\GC}) = \dim \Hom_{\GC}(1_{\GC},Y(V)) = 1,$$
$V \otimes V^{*}$ has a unique submodule $T(V)$ such that $Y(V)/T(V) \simeq 1_{\GC}$ and
a unique trivial submodule $I$. Then $\AC(V) := T(V)/(T(V) \cap I)$ is the only nontrivial
irreducible $\GC$-composition factor of $Y$.

Now assume $\GC = Sp(V)$ with $d := \dim(V) \geq 4$, or
$\GC = O(V)$ with $d := \dim(V) \geq 5$, and let $Y(V) \in\{ \SB(V), \WB(V)\}$.
If $\dim \Hom_{\GC}(Y(V),1_{\GC}) = 0$, set $I := 0$, $T(V) := Y(V)$. If
$\dim \Hom_{\GC}(Y(V),1_{\GC}) = \dim \Hom_{\GC}(1_{\GC},Y(V)) \geq 1$, $Y(V)$ has a
unique submodule $T(V)$ and a unique submodule $I$ such that
$Y(V)/T(V) \simeq I \simeq 1_{\GC}$. In both cases, $T(V)/(T(V) \cap I)$ is the only
nontrivial irreducible $\GC$-composition factor of $Y$, and we denote this subquotient by
$\TSB(V)$ if $Y = \SB$ and by $\TWB(V)$ if $Y = \WB$.

We also agree to talk about $\SB(V)$ only when $\ell \neq 2$, as $\SB(V)$ has a filtration
with quotients $\WB(V)$ and $\Vbb$ when $\ell = 2$.

\begin{lemma}\label{prim}
{\sl Assume $G < \GC$ satisfies $\MBA$. Assume $H \lhd G$ and $V \da_{H}$ is a direct sum of
$e \geq 2$ copies of an irreducible $\FF H$-module $U$. Then $O^{\ell}(H) \leq Z(G)$.}
\end{lemma}

\begin{proof}
We may assume $\dim(U) \geq 2$ as $H \leq Z(G)$ if $\dim(U) = 1$. Let $\chi$ denote the
Brauer character of $G$ on $V$. If $\GC = GL(V)$, then $(V \otimes V^{*})\da_{H}$ contains
$1_{H}$ with multiplicity $e^{2} \geq 4$. In the remaining cases, $\chi$ is real-valued and
so is the Brauer character of $U$, whence $U \simeq U^{*}$. Consider $Y = \SB$ if
$\ell \neq 2$ and $U$ is of type $+$, and $Y = \WB$ otherwise. This choice ensures that
$Y(U)$ contains the submodule $1_{H}$. Hence $Y(V)$ contains $1_{H}$ with multiplicity
$e(e+1)/2 \geq 3$. Thus in all cases $Y(V)\da_{H}$ contains $1_{H}$ with multiplicity $\geq 3$,
whence $X(V)$ contains $1_{H}$. But by assumption $G$ is irreducible on $X(V)$, so
$H$ is trivial on $X(V)$. If $\varphi$ denotes the Brauer character of $Y(V)$, then for
any $\ell'$-element $h \in H$ we see that $\varphi(h) = \varphi(1)$. Consider for instance
the case $Y = \WB$. We may assume that $h = \diag(\eps_{1}, \ldots, \eps_{d})$ for some
roots $\eps_{i}$ of unity. Then
$$d(d-1)/2 = \varphi(1) = \varphi(h) = (\chi(h)^{2} - \chi(h^{2}))/2 =
  \sum_{1 \leq i \neq j \leq d}\eps_{i}\eps_{j}.$$
Since $d = \chi(1) \geq 3$, this implies that $\eps_{1} = \ldots = \eps_{d} = \pm 1$,
whence $h \in Z(G)$. This is true for all $\ell'$-elements $h \in H$, so
$O^{\ell}(H) \leq Z(G)$. The other cases $Y = A$ or $Y = \SB$ can be dealt with similarly.
\end{proof}

\begin{propo}\label{irred1}
{\sl Let $\GC = GL(V)$, $Sp(V)$, or $O(V)$ and let $G \leq \GC$. Assume $G$ is irreducible on some
$X(V)$, with $X = \AC$ if $\GC = GL(V)$ and $X \in \{\TSB,\TWB\}$ if $\GC = Sp(V)$ or
$\GC = O(V)$. Assume in addition that $d \geq 3$ if $X = \AC$ or $X = \TSB$, and
$d \geq 5$ if $X = \TWB$. Then the following holds.

{\rm (i)} $G$ is irreducible on $V$.

{\rm (ii)} Assume $G$ is imprimitive on $V$. Then $X = \TWB$, $\ell \neq 2$, $\GC = O(V)$, and
$G$ is contained in the stabilizer $\ZZ_{2}^{d}:\SSS_{d}$ in $\GC$ of an orthonormal basis
$\{e_{1}, \ldots ,e_{d}\}$ of $V$. Moreover, $G$ acts $2$-homogeneously on
$\{\langle e_{1} \rangle_{\FF}, \ldots ,\langle e_{d}\rangle_{\FF}\}$.

{\rm (iii)} Assume in addition that $G$ satisfies the condition $\MBA$ and $d \geq 5$. Then
$G$ is primitive on $V$. If $H \lhd G$, then either $H$ is irreducible on $V$ or
$O^{\ell}(H) \leq Z(G)$.}
\end{propo}

\begin{proof}
1) Suppose $G$ is not irreducible on $V$. Let $U$ be a simple $G$-submodule of $V$ and let
$W := V/U$. First we consider the case $\GC = GL(V)$. By assumption, $G$ is irreducible on
$\AC(V)$, whence the $G$-module $V \otimes V^{*}$ has at most $3$ composition factors. But
the $G$-module $V \otimes V^{*}$ has a filtration with nonzero quotients $U \otimes W^{*}$,
$U \otimes U^{*}$, $W \otimes W^{*}$, and $W \otimes U^{*}$, a contradiction.

Next we assume $G$ is irreducible on $X(V)$ with $X = \TWB$. Then the $G$-module $\WB(V)$ has
at most $3$ composition factors: $X(V)$ and at most $2$ trivial ones. On the other hand,
the $G$-module $\WB(V)$ has a filtration with quotients $\WB(U)$, $U \otimes W$,
and $\WB(W)$, of dimensions $a(a-1)/2$, $ab$, and $b(b-1)/2$, where $a := \dim(U)$ and
$b := \dim(W)$. It follows that $a, b \leq 2$, and so $d \leq 4$, a contradiction.

Notice that the condition $d \geq 5$ is necessary, as shown by the example of
$G = SL_{2}(p^{2})$ embedded in $Sp(V)$ with $V = U \oplus U^{(p)}$ and $U$ being
the natural module for $G$.

The case $X = \TSB$ is similar.

2) Assume $G$ is imprimitive on $V$: $G$ preserves a decomposition
$V = \oplus^{t}_{i=1}V_{i}$ with $t > 1$. By irreducibility, $\dim(V_{i}) = a$ for $d = at$.

Suppose $X = \AC$, so $\GC = GL(V)$. Then $G$ also preserves the decomposition
$V^{*} = \oplus^{t}_{i=1}V^{*}_{i}$, where $V_{i}^{*} = \Ann_{V^{*}}(\oplus_{j \neq i}V_{j})$.
In particular, the permutation actions of $G$ on $\{V_{1}, \ldots ,V_{t}\}$ and on
$\{V^{*}_{1}, \ldots ,V^{*}_{t}\}$ are the same. Hence the $G$-module $V \otimes V^{*}$ has a
direct summand $\oplus^{t}_{i=1}V_{i} \otimes V_{i}^{*}$ of dimension $ta^{2} \geq 2$ and
codimension $(t^{2}-t)a^{2} \geq 2$. The irreducibility of $G$ on $\AC(V)$ forces
$2 \in \{ta^{2},(t^{2}-t)a^{2}\}$, i.e. $d = 2$, contrary to our assumption.

Suppose $X = \TSB$. Then the $G$-module $\SB(V)$ has a direct summand
$\oplus^{t}_{i=1}\SB(V_{i})$ of dimension $ta(a+1)/2 \geq 2$ and codimension
$t(t-1)a^{2}/2 \geq 1$. The irreducibility of $G$ on $\TSB(V)$ forces that either $ta(a+1)/2 = 2$
or $t(t-1)a^{2}/2 \leq 2$, i.e. $d = 2$, contrary to our assumption.

Finally, suppose $X = \TWB$. Then the $G$-module $\WB(V)$ has a direct summand
$U : = \oplus^{t}_{i=1}\WB(V_{i})$ of dimension $ta(a-1)/2 \geq 0$ and codimension
$t(t-1)a^{2}/2 \geq 1$. If $a \geq 3$, or if $a = 2$ but $t \geq 3$, then
$\dim(U), \codim(U) \geq 3$. Since $d = at \geq 5$, the irreducibility of $G$ on $\TWB(V)$
forces $a = 1$ and $t = d$. This also implies that $G$ acts primitively on
$\{V_{1}, \ldots ,V_{d}\}$.

3) We continue to study the situation of 2) with $X = \TWB$. Denote
$\Om := \{V_{1}, \ldots ,V_{d}\}$. We have shown that $V_{i} = \langle e_{i} \rangle_{\FF}$ for
some basis $\{e_{1}, \ldots ,e_{d}\}$ of $V$.

We claim that $G$ is $2$-homogeneous on $\Om$. Assume the contrary: $G$ has more than one
orbit on $\Om_{2} := \{ \{\al,\beta\} \mid \al \neq \beta \in \Om\}$. Let $\OC$ be such an orbit.
If $\OC$ has length $1$, say $\OC = \{\{V_{1},V_{2}\}\}$, then $G$ preserves $\{V_{1},V_{2}\}$. 
If $\OC$ has length $2$, say $\OC = \{\{V_{1},V_{2}\},\{V_{1},V_{3}\}\}$ or
$\OC = \{\{V_{1},V_{2}\},\{V_{3},V_{4}\}\}$, then $G$ preserves
$\{V_{1},V_{2},V_{3}\}$ or $\{V_{1},V_{2},V_{3},V_{4}\}$. In both cases $G$ is intransitive on
$\Om$, contrary to the irreducibility of $G$ on $V$. So the length of any $G$-orbit on $\Om_{2}$
is at least $3$. It follows that $G$ has a submodule of dimension and codimension $\geq 3$ in
$\WB(V)$, contrary to the irreducibility of $G$ on $\TWB(V)$.

Let $(\cdot,\cdot)$ be the $\GC$-invariant bilinear form on $V$. Suppose that
$(e_{1},e_{2}) = 0$. Since $G$ is $2$-homogeneous on $\Om$, we have $(e_{i},e_{j}) = 0$
whenever $i \neq j$. This can happen only when $\ell \neq 2$, $\GC = O(V)$, and clearly
$(e_{1},e_{1}) = \ldots = (e_{d},e_{d})$, so we may assume that $\{e_{1}, \ldots ,e_{d}\}$ is
an orthonormal basis of $V$, as recorded in (ii).

From now on we will assume that $(e_{i},e_{j}) \neq 0$ whenever $i \neq j$. Observe that
if $g \in G$ preserves say $V_{1}$ and $V_{2}$, then $g$ fixes $e_{1} \wedge e_{2}$. Indeed, by
assumption $g(e_{i}) = a_{i}e_{i}$ for some $a_{i} \in \FF$ and $i = 1,2$. Then
$0 \neq (e_{1},e_{2}) = (g(e_{1}),g(e_{2})) = a_{1}a_{2}(e_{1},e_{2})$, so
$a_{1}a_{2} = 1$, whence $g$ fixes $e_{1} \wedge e_{2}$.

Similarly, if $(\cdot,\cdot)$ is alternating and $g \in G$ fixes
$\langle e_{1} \wedge e_{2} \rangle_{\FF}$, then $g$ fixes $e_{1} \wedge e_{2}$. Beacause of the
above discussion, it suffices to consider the case
$g~:~e_{1} \mapsto b_{1}e_{2},~e_{2} \mapsto b_{2}e_{1}$ for some $b_{i} \in \FF$. Then
$$0 \neq (e_{1},e_{2}) = (g(e_{1}),g(e_{2})) = b_{1}b_{2}(e_{2},e_{1}) =
  -b_{1}b_{2}(e_{1},e_{2}),$$
so $b_{1}b_{2} = -1$, whence $g$ fixes $e_{1} \wedge e_{2}$ as stated.

4) We continue to study the situation of 3). For $i \neq j$, let
$G_{ij} := Stab_{G}(\langle e_{i} \wedge e_{j} \rangle) = Stab_{G}(\{V_{i},V_{j}\})$, and let
$\lam_{ij}$ be the character of $G_{ij}$ acting on $e_{i} \wedge e_{j}$. Since $G$ is
$2$-homogeneous on $\Om$, $\WB(V) = \lam_{12} \ua G$. We claim that
$\dim \End_{G}(\WB(V))$ is at least $2$ if $\ell \neq 2$, and at least $3$ if
$(\cdot,\cdot)$ is alternating. Indeed, by Mackey's theorem we have
$$\dim \End_{G}(\WB(V)) \geq 1 + \dim \Hom_{G_{12} \cap G_{13}}(\lam_{12},\lam_{13}) +
  \dim \Hom_{G_{12} \cap G_{34}}(\lam_{12},\lam_{34}).$$
Let $g \in G_{12} \cap G_{13}$. Then $g$ fixes $V_{1}$, $V_{2}$, and $V_{3}$, whence $g$
fixes $e_{1} \wedge e_{2}$ and $e_{1} \wedge e_{3}$ by 3). Thus
$\lam_{12} = \lam_{13} = 1$ on $G_{12} \cap G_{13}$. This implies that
$\dim \End_{G}(\WB(V)) \geq 2$. Now assume that $(\cdot,\cdot)$ is alternating. Consider
any $g \in G_{12} \cap G_{34}$. Then $g$ fixes $\langle e_{1} \wedge e_{2}\rangle_{\FF}$, so $g$
fixes $e_{1} \wedge e_{2}$ by 3). Similarly, $g$ fixes $e_{3} \wedge e_{4}$. Thus
$\lam_{12} = \lam_{34} = 1$ on $G_{12} \cap G_{34}$. This implies that
$\dim \End_{G}(\WB(V)) \geq 3$.

5) Here we complete the analysis of the situation of 3). Suppose $\ell \neq 2$ and $(\cdot,\cdot)$
is symmetric. Then $G$ is irreducible on $\TWB(V) = \WB(V)$, but this contradicts the inequality
$\dim \End_{G}(\WB(V)) \geq 2$ proved in 4). So we may assume that $(\cdot,\cdot)$ is
alternating, whence $\dim \End_{G}(\WB(V)) \geq 3$ as proved in 4). Fix three linearly
independent $f,g,h \in \End_{G}(\WB(V))$. Since $(\cdot,\cdot)$ is alternating, the
$\GC$-module $\WB(V)$ contains a trivial submodule $I$.

Assume $f(I), g(I) \neq 0$. Then $f(I), g(I) \simeq 1_{G}$. If $f(I) \neq g(I)$, then
$\dim\Hom_{G}(1_{G},\WB(V)) \geq 2$, contradicting the irreducibility of $G$ on $V$. So
$f(I) = g(I)$. Replacing $g$ by $g -\al f$ for a suitable $\al \in \FF$, we get $g(I) = 0$.
Thus we may assume $g(I) = h(I) = 0$. Denoting $\TWB(V)$ by $X$ for short. Then $X$ is a
submodule of the $\GC$-module $\WB(V)/I$. If $g(X) = h(X) = 0$, then
$g,h \in \Hom_{G}(1_{G},\WB(V))$ which is of dimension $1$, whence $g,h$ are dependent, a
contradiction. So we may assume that $h(X) \neq 0$, hence $X$ is a submodule of the
$G$-module $\WB(V)$ as $G$ is irreducible on $X$. Now if the $\GC$-module $\WB(V)$ has only
one trivial composition factor, then $\WB(V) \simeq X \oplus 1_{G}$ as $G$-modules and
$\dim \End_{G}(\WB(V)) = 2$, again a contradiction. So the $\GC$-module $\WB(V)$ has $I$ as
a composition factor of multiplicity $2$.

The action of $G$ on $\WB(V)/X$ induces a homomorphism $\pi$ from $G$ into the $\ell$-group
$\left\{ \begin{pmatrix}1 & * \\0 & 1 \end{pmatrix} \right\}$. Then
$K := \Ker(\pi) \geq O^{\ell}(G)$ acts trivially on the $2$-dimensional module $\WB(V)/X$,
so we get
$$2 \leq \dim \Hom_{K}(\WB(V),1_{K}) = \dim \Hom_{K}(1_{K},\WB(V)) \leq
  \dim \Hom_{K}(1_{K}, V \otimes V^{*}),$$
whence $V_{K}$ is reducible by Schur's lemma. Write $V\da_{K} = e(\oplus^{s}_{j=1}U_{j})$
where $U_{1}, \ldots ,U_{s}$ are distinct $K$-irreducibles. Then $es > 1$ and $G$ permutes the
$s$ isotypic components of $K$ on $V$. The results of 2) imply that either $s = 1$ or $s = d$.
First we consider the former case, that is $s = 1$. Since $V$ is self-dual, the $K$-module
$U := U_{1}$ is also self-dual. Observe that $\WB(V) \da_{K}$ contains $1_{K}$ with multiplicity
$\geq 3$. (For it is so if $e \geq 3$, or if $e = 2$ and $U$ is of type $-$. If $e = 2$
and $U$ supports no alternating form, then $\ell \neq 2$ and
$\WB(V) \da_{K} \simeq \WB(U) \oplus U \otimes U \oplus \WB(U)$ contains $1_{K}$ with
multiplicity $1$, contrary to the inequality $2 \leq \dim \Hom_{K}(\WB(V),1_{K})$.) As in the
proof of Lemma \ref{prim}, we can conclude that $Z(G) \geq O^{\ell}(K) = O^{\ell}(G)$. It follows
that $G = Z(G)P$ for some Sylow $\ell$-subgroup $P$ of $G$. But in this case the dimension of
any $\ell$-modular absolutely irreducible representation of $G$ is $1$. In particular,
$\dim(V) = 1$, a contradiction.

So we must have $s = d$ and $e = 1$. By 3), $G$, and so $G/K$, acts $2$-homogeneously on
$\{U_{1}, \ldots ,U_{d}\}$. Thus $d(d-1)/2$ divides $G/K$, which is an $\ell$-power as
$K \geq O^{\ell}(G)$. This is impossible as $d \geq 5$. This contradiction finishes the proof
of (ii).

6) Assume $G$ satisfies $\MBA$ and $d \geq 5$. If $\ell \neq 2$ then $\MBA$ implies that $G$ is
irreducible on $\AC(V)$ or $\TSB(V)$. Hence (ii) implies that $G$ is primitive on $V$.

Assume $H \lhd G$ and $H$ is not irreducible on $V$. Then Clifford's theorem and the primitivity
of $G$ implies that $V\da_{H}$ is the direct sum of $e \geq 2$ copies of an irreducible
$H$-module $U$. By Lemma \ref{prim}, $O^{\ell}(H) \leq Z(G)$.
\end{proof}

\begin{corol}\label{irred2}
{\sl Assume $d \geq 5$, $G < \GC$ and $G$ is reducible on $V$. Then for all
$k \geq 2$, $\MK$ fails for $G$.}
\end{corol}

\begin{proof}
By Proposition \ref{irred1}(i), $\MB$ fails for $G$. Next, $V$ is a $\GC$-composition
factor inside $V^{\otimes 2} \otimes V^{*}$, so $\MC$ fails for $G$. Since
$\MK$ implies $M_{2k-4}(\GC)$, $\MK$ fails for $G$.
\end{proof}

\begin{lemma}\label{indec1}
{\sl Assume either $\GC = GL(V)$ with $d \geq 3$, or $\GC = Sp(V)$, $O(V)$ with $d \geq 5$.
Assume $G \leq \GC$ satisfies the condition $\MBA$. Consider $X = \AC$ if $\GC = GL(V)$ and
$X \in \{\TSB,\TWB\}$ otherwise. Then the $G$-module $T(V)$ is indecomposable.}
\end{lemma}

\begin{proof}
Assume the contrary. Then $T(V) \simeq X(V) \oplus 1_{G}$ as $G$-modules. In particular,
$1 = \dim \Hom_{G}(1_{G},T(V)) \leq \dim \Hom_{G}(1_{G},Y(V)) = \dim \Hom_{G}(Y(V),1_{G})$,
whence $Y(V)/T(V) \simeq 1_{G}$. The action of $G$ on $Y(V)/X(V)$ induces a homomorphism
$\pi$ from $G$ into the $\ell$-group
$\left\{ \begin{pmatrix}1 & * \\0 & 1 \end{pmatrix} \right\}$.
Clearly $\Ker(\pi) \geq H := O^{\ell}(G)$. Since $H$ acts trivially on the
$2$-dimensional module $Y(V)/X(V)$, we get
$$2 \leq \dim \Hom_{H}(Y(V),1_{H}) = \dim \Hom_{H}(1_{H},Y(V)) \leq
  \dim \Hom_{H}(1_{H}, V \otimes V^{*}),$$
whence $V_{H}$ is reducible by Schur's lemma. By Proposition \ref{irred1},
$Z(G) \geq O^{\ell}(H) = H$. It follows that $G = Z(G)P$ for some Sylow $\ell$-subgroup
$P$ of $G$. But in this case the dimension of any $\ell$-modular absolutely irreducible
representation of $G$ is $1$. In particular, $\dim(V) = 1$, a contradiction.
\end{proof}

\begin{lemma}\label{indec2}
{\sl Let $T$ be an $\FF G$-module with $\soc(T) = 1_{G}$ and $X := T/\soc(T)$ is nontrivial
irreducible. Assume $\dim \Hom_{C}(T,1_{C}) \geq s \geq 1$ for some $C \leq G$. Then

{\rm (i)} $X$ is a composition factor of $(1_{C}) \ua G$ with multiplicity at least $s$.

{\rm (ii)} Assume $C \lhd N \leq G$ and $N/C$ is abelian. Then $X \hra \lam \ua G$ for some
$1$-dimensional $\FF N$-representation $\lam$.}
\end{lemma}

\begin{proof}
(i) is just \cite[Lemma 2.1]{MMT}.

(ii) Let $M$ be the smallest $C$-submodule of $T$ such that all $C$-composition factors of
$T/M$ are trivial. By assumption, $M \neq T$, and clearly $M$ is $N$-stable. Next,
$T/M$ has a simple quotient $1_{C}$ and clearly $1_{C}$ extends to $N$ since $N/C$ is
abelian. By \cite[Lemma 2.2]{MT1} the $N$-module $T/M$, and so $T$, has a $1$-dimensional
quotient say $\lam$, i.e. $0 \neq \Hom_{N}(T,\lam) \simeq \Hom_{G}(T,\lam \ua G)$. Pick a nonzero
$f \in \Hom_{G}(T,\lam \ua G)$. If $\Ker(f) \supseteq \soc(T)$, then $\Ker(f) = \soc(T)$ and so
$0 \neq f \in \Hom_{G}(X,\lam \ua G)$, whence $X \hra \lam \ua G$. Otherwise $\Ker(f) = 0$,
$f$ embeds $T$ in $\lam \ua G$ and we are again done.
\end{proof}

A key ingredient of our arguments is the following statement (see also \cite[Prop. 2.3]{MMT}):

\begin{propo}\label{bound1}
{\sl Let $\GC = GL(V)$ with $d \geq 3$ or $\GC \in \{Sp(V),O(V)\}$ with $d \geq 5$.
Assume $G$ is a subgroup of $\GC$ satisfying the condition $\MBA$.
Assume $C \lhd N \leq G$, $N/C$ is abelian, and $V\da_{C} = V_{1} \oplus \ldots \oplus V_{t}$
for some $t \geq 2$ nonzero $C$-submodules $V_{i}$. Then the following hold.

{\rm (i)} If $\GC = GL(V)$ then $\AC(V)$ enters $(1_{C}) \ua G$ with multiplicity $\geq t-1$,
and $\AC(V) \hra \lam \ua G$ for some $1$-dimensional $\FF N$-module $\lam$.

{\rm (ii)} Assume $V = V^{*}$ as $\GC$-modules and $\Hom_{C}(V_{i},V^{*}_{j}) = 0$ whenever
$i \neq j$. Then there is an $X \in \{\TSB,\TWB\}$ such that $X(V)$ enters $(1_{C}) \ua G$ with
multiplicity $\geq t-1$, and $X(V) \hra \lam \ua G$ for some $1$-dimensional $\FF N$-module $\lam$.

{\rm (iii)} Assume $V = V^{*}$ as $\GC$-modules. If $\ell \neq 2$, assume that there are some
$i \neq j$ such that $\Hom_{C}(V_{i},V^{*}_{j}) \neq 0$. If $\ell = 2$, assume that $V\da_{C}$ is
the sum of at least $3$ indecomposable summands. Then there is an $X \in \{\TSB,\TWB\}$ such that
$X(V)$ enters $(1_{C}) \ua G$ and $\lam \ua G$ for some $1$-dimensional $\FF N$-module $\lam$.}
\end{propo}

\begin{proof}
By Lemma \ref{indec1}, the $G$-module $T(V)$ is indecomposable for $X = \AC$ if $\GC = GL(V)$ and
$X \in \{\TSB,\TWB\}$ otherwise. By assumption, $G$ is irreducible on $T(V)/\soc(T(V))$.

1) First we consider the case $\GC = GL(V)$ and set $Y(V) = V \otimes V^{*}$. Since
$(V^{*})\da_{C} \simeq \oplus^{t}_{i=1}V^{*}_{i}$, we see that $Y(V)$ contains the submodule
$\oplus^{t}_{i=1}V_{i} \otimes V^{*}_{i}$, whence $\dim \Hom_{C}(Y(V),1_{C}) \geq t$. Since
$\dim (Y(V)/T(V)) \leq 1$, the restriction map
$$\Res_{T(V)}~:~\Hom_{C}(Y(V),1_{C}) \to \Hom_{C}(T(V),1_{C})$$
has kernel of dimension $\leq 1$.
Hence $\dim \Hom_{C}(T(V),1_{C}) \geq t-1 \geq 1$. Now (i) follows from Lemma \ref{indec2}
applied to $T = T(V)$.

From now on we may assume that $V = V^{*}$ as $\GC$-modules and let $\bb$ be a $\GC$-invariant
nondegenerate bilinear form on $V$.

2) Here we consider the case of (ii). For $B := \oplus_{j \neq i}V_{i}$ we have
$B/(B \cap V_{i}^{\perp}) \simeq (B + V_{i}^{\perp})/V_{i}^{\perp} \leq V/V_{i}^{\perp}
 \simeq V_{i}^{*}$. If $B/(B \cap V_{i}^{\perp}) \neq 0$, then we get a nonzero
$C$-homomorphism $B \to V^{*}_{i}$, contradicting the assumption that
$\Hom_{C}(V_{j},V^{*}_{i}) = 0$ for all $j \neq i$. Hence $B \subseteq V_{i}^{\perp}$.
Comparing dimension we see that $B = V_{i}^{\perp}$. Thus $V_{i} \cap V_{i}^{\perp} = 0$,
whence $\bb$ is nondegenerate on $V_{i}$.

Choose $Y = \SB$ if $\ell \neq 2$ and $\bb$ is symmetric, and $Y = \WB$ otherwise. Then
$\Hom_{C}(Y(V_{i}),1_{C}) \neq 0$, whence $\dim \Hom_{C}(Y(V),1_{C}) \geq t$. Now we can
proceed as in 1).

3) Assume $V = V^{*}$, $\ell \neq 2$, and $\Hom_{C}(V_{i},V^{*}_{j}) \neq 0$ for some $i \neq j$. 
Choose $Y = \SB$ if $\bb$ is alternating and $Y = \WB$ if $\bb$ is symmetric (recall that
$\ell \neq 2$ here). This choice ensures that $X(V) = Y(V)$. Since $d \geq 5$ and
$G$ satisfies $\MBA$, $G$ is irreducible on $X(V)$. Next, $Y(V)\da_{C}$ contains the
direct summand $V_{i} \otimes V_{j}$, hence
$$\dim \Hom_{C}(X(V),1_{C}) \geq
  \dim \Hom_{C}(V_{i} \otimes V_{j},1_{C}) = \dim \Hom_{C}(V_{i},V^{*}_{j}) \neq 0,$$
and so $X(V)$ is a submodule of $(1_{C}) \ua G$ by irreducibility of $X(V)$. Arguing as in
the proof of Lemma \ref{indec2} we see that $X(V)$ is a submodule of $\lam \ua G$ for some
$1$-dimensional $\FF N$-module $\lam$.

4) Finally, assume that $V = V^{*}$, $\ell = 2$ and $V\da_{C} = \oplus^{s}_{i=1}U_{i}$ is the
sum of $s \geq 3$ indecomposable direct summands. Choose $Y = \WB$. As above, it suffices to
show that $\dim \Hom_{C}(Y(V),1_{C}) \geq 2$. We have
$\oplus^{s}_{i=1}U_{i} = V\da_{C} \simeq V^{*}\da_{C} \simeq \oplus^{s}_{i=1}U^{*}_{i}$. By the
Krull-Schmidt theorem, each $U^{*}_{i}$ is isomorphic to some $U_{j}$. Thus we get an
involution permutation $\tau~:~i \mapsto j$ of $\{1,2, \ldots ,s\}$.
If $\tau$ contains $\geq 2$ transpositions, we may assume that $U_{2} \simeq U^{*}_{1}$ and
$U_{4} \simeq U^{*}_{3}$, whence $Y(V)\da_{C}$ has a submodule
$$U_{1} \otimes U_{2} \oplus U_{3} \otimes U_{4} \simeq
  U^{*}_{2} \otimes U_{2} \oplus U^{*}_{4} \otimes U_{4}$$
and so $\dim \Hom_{C}(Y(V),1_{C}) \geq 2$. Assume $\tau$ has $\geq 2$ fixed points. Then
we may assume that $U_{1} \simeq U^{*}_{1}$ and $U_{2} \simeq U^{*}_{2}$. Since $U_{j}$,
$j = 1,2$, is self-dual and $\ell = 2$, $Y(U_{i})$ has $1_{C}$ as a submodule. Now $Y(V)\da_{C}$
has a submodule $Y(U_{1}) \oplus Y(U_{2})$, so $\dim \Hom_{C}(Y(V),1_{C}) \geq 2$. Finally,
assume $\tau$ has $< 2$ transpositions and $<2$ fixed points. Since $s \geq 3$, this implies
that $s = 3$ and we may assume $U_{1} \simeq U^{*}_{1}$ and $U_{3} \simeq U^{*}_{2}$. In
this case, $Y(V)\da_{C}$ has a submodule $Y(U_{1}) \oplus U_{2} \otimes U_{3}$, so again
$\dim \Hom_{C}(Y(V),1_{C}) \geq 2$.
\end{proof}

The following lemma is useful in checking $\MK$ with large $k$:

\begin{lemma}\label{weight}
{\sl Assume $d = \dim(V) \geq \max\{4,2k\}$ and let $e := \lfloor d/2 \rfloor$.

{\rm (i)} The $\GC$-module $\wedge^{k}(V)$ inside $V^{\otimes k}$ contains a composition
factor $W$ of dimension at least $2^{k} \begin{pmatrix}e\\k\end{pmatrix}$. In fact
$\wedge^{k}(V)$ is irreducible if $\GC = GL(V)$.

{\rm (ii)} Assume $\ell \neq 2$, $\GC = SO(V)$, and $k <d/2$. Then $\wedge^{k}(V)$ is
$\GC$-irreducible. }
\end{lemma}

\begin{proof}
(i) Consider the highest weight $\om$ of $\wedge^{k}(V)$. If $\GC \neq O(V)$, or
if $\GC = O(V)$ but $d \geq 2k+3$, then $\om = \om_{k}$. Assume $\GC = O(V)$.
If $d = 2k+1$, then $\om = 2\om_{k}$. If $d = 2k$ then $\om = 2\om_{k-1}+2\om_{k}$,
and if $d = 2k+2$ then $\om = \om_{k}+\om_{k+1}$, cf. \cite{Dyn}. In all cases,
the orbit of $\om$ under the action of the Weyl group of $\GC$ has length at least
$2^{k} \begin{pmatrix}e\\k\end{pmatrix}$ (and exactly
$\begin{pmatrix}d\\k\end{pmatrix} = \dim(\wedge^{k}(V))$ if $\GC = GL(V)$). Let
$W$ be the irreducible $\GC$-module $L(\om)$. Since $L(\om)$ has to afford the whole
orbit of $\om$, the statements follow.

(ii) See \cite[(8.1)]{S2}.
\end{proof}

We will need some statements about $\GC$-composition factors of $\VE$ with $e \leq 4$.

\begin{propo}\label{v4}
{\sl Let $\ell \neq 2,3$, $1 \leq e \leq 4$, $V = \FF^{d}$, and $\VC = \CC^{d}$.

{\rm (i)} Assume $\XC = GL(V)$ and $d \geq 4$. Then the $\XC$-module $\VE$ is
semisimple. Moreover, all $GL(\VC)$-composition factors of $\VC^{\otimes e}$ remain
irreducible modulo $\ell$. If $d \geq 8$, then all distinct $\XC$-composition factors of
$V^{\otimes e}$ are of distinct degrees.

{\rm (ii)} Assume $\XC = Sp(V)$. Then $S^{e}(V)$ is an irreducible direct summand of
the $\XC$-module $\VE$.

{\rm (iii)} Assume $\XC = SO(V)$ and $d \geq 9$. Then $\wedge^{e}(V)$ is an irreducible
direct summand of the $\XC$-module $\VE$.}
\end{propo}

\begin{proof}
1) We may assume that $V$ is obtained from $\VC$ by reducing modulo $\ell$.
Consider the natural action of $S := \Sym_{e}$ on $\VE$ by permuting the $e$ components.
Clearly, this action commutes with the natural action of $\GC := GL(V)$ on $\VE$. Hence we get
a homomorphism $\pi~:~\FF S \to \End_{\GC}(\VE)$ of algebras. Since $\ell \neq 2,3$,
we can decompose $\FF S = \oplus_{\lam \vdash e}A_{\lam}$ as direct sum of matrix algebras.
Claim that $\pi$ is injective. Assume the contrary. Then $\Ker(\pi)$ contains $A_{\lam}$ for
some $\lam \vdash e$. Observe that the Young symmetrizer $c_{\lam}$ annihilates all
$A_{\mu}$ with $\mu \neq \lam$. Hence $c_{\lam}$ annihilates the action of $\FF S$ on $\VE$,
whence $c_{\lam}(\VE) = 0$. On the other hand, it is well known that $c_{\lam}(\VC^{\otimes e})$ 
is the direct sum of $\dim(S^{\lam})$ copies of $\SSS_{\lam}(\VC)$, where $S^{\lam}$ is the Specht 
module of $S$ labelled by $\lam$ (cf. \cite{FH} for the definition of $\SSS_{\lam}$). Reducing 
modulo $\ell$, we see that $M_{\lam} := c_{\lam}(\VE)$ has the same
dimension as of $c_{\lam}(\VC^{\otimes e})$ and so it is nonzero, a contradiction.

The above argument also shows that $\VE = \oplus_{\lam \vdash e}M_{\lam}$ as
$\GC$-modules. Setting $m = \dim(S^{\lam})$, we see that $A_{\lam} \simeq Mat_{m}(\FF)$ acts
on $M_{\lam}$. We can decompose $c_{\lam}$ as the sum of $m$ idempotents $e_{ii}$,
$1 \leq i \leq m$, of $Mat_{m}(\FF)$. This allows us to write
$M_{\lam} = \oplus_{1 \leq i \leq m}e_{ii}M_{\lam}$ as $\GC$-modules. Moreover,
the $e_{ii}M_{\lam}$'s are transitively permuted by $Mat_{m}(\FF)$, so they are isomorphic
$\GC$-modules. Denoting any of them by $\SSS_{\lam}(V)$, we see that it can be obtained
by reducing $\SSS_{\lam}(\VC)$ modulo $\ell$. Moreover, $\SSS_{\lam}(V)$ is a Weyl module
$W(\om)$ with highest weight say $\om$ and simple head $L(\om)$. Similarly,
$(V^{*})^{\otimes e} = \oplus_{\lam \vdash e}c_{\lam}((V^{*})^{\otimes e})$, and each
$c_{\lam}((V^{*})^{\otimes e})$ is the direct sum of $\dim(S^{\lam})$ copies of
$\SSS_{\lam}(V^{*})$.

The crucial observation is that the $\GC$-module $V$ is a Weyl module and at the same time
a dual Weyl module. In other words, $V^{*} \simeq V^{\#}$ as $\GC$-modules, where
$V^{\#}$ is $V$ as an $\FF$-space, but with $g \in \GC$ acting as $\tau(g)$ on $V$,
and $\tau~:~g \mapsto \tn g^{-1}$ is an outer automorphism of $\GC$. Since
$c_{\lam}(\VE)$ and $c_{\lam}((V^{*})^{\otimes e})$ are annililated by all $c_{\mu}$ with
$\mu \neq \lam$, applying $\#$ we get $c_{\lam}(\VE) = (c_{\lam}((V^{*})^{\otimes e}))^{\#}$.
The $\GC$-modules $\SSS_{\lam}(V)$ and $\SSS_{\lam}(V^{*})$ are indecomposable as Weyl modules,
so by the Krull-Schmidt theorem, $\SSS_{\lam}(V) = (\SSS_{\lam}(V^{*}))^{\#}$. The left-hand
side has $L(\om)$ as simple head, meanwhile the right-hand side has
$(L(\om)^{*})^{\#} \simeq L(\om)$ as simple socle. Thus $L(\om)$ is simultaneously simple
head and simple socle of $W(\om)$, and it is a composition factor of multiplicity $1$ of
$W(\om)$. It follows that $W(\om) = L(\om)$. We have shown that $\VE$ is semisimple and
$\SSS_{\lam}(\VC) (\mod \ell)$ is irreducible.

Finally, we indicate the dimension of $\SSS_{\lam}(V)$. This is just the multiplicity of
$S^{\lam}$ in the $\Sym_{e}$-module $V^{\otimes e}$. If $e = 1$ we get $\lam = (1)$ and
$\dim(\SSS_{\lam}(V)) = d$. If $e = 2$ we get $\begin{pmatrix}d+1\\2\end{pmatrix}$, resp.
$\begin{pmatrix}d\\2\end{pmatrix}$, for $\lam = (2)$, resp. $(1^{2})$. If $e = 3$ we
get $\begin{pmatrix}d+2\\3\end{pmatrix}$, resp. $d(d^{2}-1)/3$,
$\begin{pmatrix}d\\3\end{pmatrix}$, for $\lam = (3)$, resp. $(2,1)$, $(1^{3})$.
If $e = 4$ we get $\begin{pmatrix}d+3\\4\end{pmatrix}$, resp. $d(d+2)(d^{2}-1)/8$,
$d(d-2)(d^{2}-1)/8$, $d^{2}(d^{2}-1)/12$, $\begin{pmatrix}d\\4\end{pmatrix}$, for
$\lam = (4)$, resp. $(3,1)$, $(2,1^{2})$, $(2,2)$, $(1^{4})$ (and the highest weights are
$4\om_{1}$, $2\om_{1}+\om_{2}$, $\om_{1}+\om_{3}$, $2\om_{2}$, and $\om_{4}$). In particular,
if $d \geq 8$ then all the dimensions in $V^{\otimes e}$ are distinct.

2) Now (i) follows from the results proved in 1). For (ii), resp. (iii), observe that
$Y := S^{e}(V)$, resp. $\wedge^{e}(V)$, is a direct summand of the $\XC$-module $\VE$ by 1),
and $V \simeq V^{*}$ as $\XC$-modules. Moreover, $Y$ is a Weyl module for $\XC$. Arguing as
in 1) with $\tau$ replaced by the trivial map, we conclude that $Y$ is irreducible.
\end{proof}

\begin{remar}\label{v4e}
{\em (i) {\sl An analogue of Proposition \ref{v4} holds for $V^{\otimes e}$ if
$\ell > e$ and $\dim(V) > 2e$.} Indeed, we need only the fact that the algebra $\FF\Sym_{e}$
is semisimple, and $\wedge^{e}(\VC)$ is irreducible over $O(V)$.

(ii) {\sl Proposition \ref{v4} fails for $\ell = 2,3$.} For, the Weyl module $S^{4}(V)$
has highest weight $4\om_{1}$ and so by Steinberg's tensor product theorem it has a proper
simple head $V \otimes \Vcc$ when $\ell = 3$ and $\Vdd$ when $\ell = 2$. Furthermore,
the $O(V)$-module $\wedge^{4}(V)$ is reducible when $\ell = 2$ (since $O(V)$ embeds in
$Sp(V)$).}
\end{remar}

\begin{remar}\label{m126}
{\em {\sl Assume $\ell \geq 7$, $\GC = GL(V)$, and $d > 14$. Then $\MF$ implies $\MCA$.} Indeed,
consider the corresponding module $\VC$ for $\GCC$. Then
$S^{2}(S^{3}(\VC)) = \SSS_{(6,0)}(\VC) \oplus \SSS_{(4,2)}(\VC)$ is a direct sum of two
$\GCC$-irreducibles (with highest weights $6\om_{1}$ and $2(\om_{1}+\om_{2})$). By Remark
\ref{v4e}, the same is true for $\GC$. Now if $G \leq \GC$ satisfies $\MF$, then $G$ also has
two composition factors on $S^{2}(S^{3}(V))$, whence $G$ is irreducible on $S^{3}(V)$.}
\end{remar}

\begin{propo}\label{v3}
{\sl Let $\GC = SO(V)$ with $\ell > 3$, $d \geq 10$, and $\ell|(d+2)$. Then all
$\GC$-composition factors of $\SC(V)$ are $L(3\om_{1})$ (with multiplicity $1$), and
$V$. Furthermore, $\dim(L(3\om_{1}))$ is at least $(2m+1)(2m^{2}+2m+3)/3$ if $d = 2m+1$, resp.
$4m(m^{2}+2)/3$ if $d = 2m$. Moreover, all $O(V)$-composition factors of $\SC(V)$ have dimension
equal to $\dim(L(3\om_{1}))$ (with multiplicity $1$), or $\dim(V)$.}
\end{propo}

\begin{proof}
Let $\VC = \CC^{d}$ be endowed with a nondegenerate symmetric bilinear form and let
$\GCC = SO(\VC)$. Furthermore, $\LC(\om)$ denotes the irreducible $\GCC$-module with highest
weight $\om$. Since $\ell > 3$, the weight $3\om_{1}$ is restricted. Hence, by Premet's
theorem \cite{Pre}, $L(3\om_{1})$ affords the same weights as of $\LC(3\om_{1})$. Observe that
$\SC(\VC) = \LC(3\om_{1}) \oplus \VC$. We may assume that $V$ is obtained from $\VC$ by
reducing modulo $\ell$.

1) First we consider the case $d = 2m+1$. One can check that the dominant weights
afforded by $\LC(3\om_{1})$ are $3\om_{1}$, $\om_{1} + \om_{2}$, $\om_{3}$, $2\om_{1}$,
$\om_{2}$ (all with multiplicity $1$), $\om_{1}$ (with multiplicity $m$), and $0$ (with
multiplicity $m$). The length of the orbit of those weights under the Weyl group is
$2m$, $4m(m-1)$, $4m(m-1)(m-2)/3$, $2m$, $2m(m-1)$, $2m$, and $1$, respectively. Hence Premet's
theorem implies that $\dim(L(3\om_{1}))$ is at least the sum of the lengths of all these orbits,
which is $(2m+1)(2m^{2}+2m+3)/3$. It is well-known that the Weyl module $V(3\om_{1})$ with
highest weight has the same composition factors (with the same multiplicities) as of the
$\GC$-module $\LC(3\om_{1}) (\mod \ell)$. In particular, the multiplicity of $L(3\om_{1})$ in
$\SC(V)$ is $1$. Let $L(\mu)$ be any $\GC$-composition factor of $\SC(V)$ other than
$L(3\om_{1})$. Then $\mu$ cannot be any of the weights $\om_{1} + \om_{2}$, $\om_{3}$,
$2\om_{1}$, or $\om_{2}$, since each of them occurs in $\LC(3\om_{1})$ with multiplicity $1$.
Observe that $\mu \neq 0$. Assume the contrary. Then the weights $3\om_{1}$ and
$0$ are linked by the linkage principle \cite{Jan}. This means that $3\om_{1}$ and $0$ belong
to the same orbit under the action of the affine Weyl group, which is generated by the maps
$s_{\al,k}~:~\lam \mapsto s_{\al}(\lam+\rho) - \rho + k\ell\al$, $\al$ any simple root and
$k$ any integer. Here $s_{\al}$ is the reflection corresponding to $\al$, and
$\rho = \sum^{m}_{i=1}\om_{i}$. Direct computation shows that
$|s_{\al,k}(\lam) + \rho|^{2} \equiv |\lam + \rho|^{2} (\mod \ell)$ for any weight $\lam$.
So the assumption that $3\om_{1}$ and $0$ are linked implies that $\ell$ divides
$|3\om_{1}+\rho|^{2} - |\rho|^{2} = 6m+3$. But this is impossible, since $3 < \ell|(2m+3)$.
Thus $\mu = \om_{1}$, as stated.

2) Next we consider the case $d = 2m$. One can check that the dominant weights
afforded by $\LC(3\om_{1})$ are $3\om_{1}$, $\om_{1} + \om_{2}$, $\om_{3}$ (all with
multiplicity $1$), and $\om_{1}$ (with multiplicity $m-1$). The length of the orbit of those
weights under the Weyl group is $2m$, $4m(m-1)$, $4m(m-1)(m-2)/3$, and $2m$, respectively.
Hence Premet's theorem implies that $\dim(L(3\om_{1}))$ is at least the sum of the lengths of
all these orbits, which is $4m(m^{2}+2)/3$. As above, the multiplicity of $L(3\om_{1})$ in
$\SC(V)$ is $1$. Let $L(\mu)$ be any $\GC$-composition factor of $\SC(V)$ other than
$L(3\om_{1})$. Then $\mu$ cannot be any of $\om_{1} + \om_{2}$ or $\om_{3}$, since each of them
occurs in $\LC(3\om_{1})$ with multiplicity $1$. Hence $\mu = \om_{1}$, as desired.

3) For $O(V)$-composition factors of $\SC(V)$, just observe that $\om_{1}$ is $O(V)$-stable, so
$V$ and $L(3\om_{1})$ both extend to $O(V)$-modules.
\end{proof}

Sometimes we will have to use the condition $\MDA$ to exclude the examples that survive
the condition $\MBA$:

\begin{propo}\label{bound2}
{\sl Assume $G < \GC$ is a finite subgroup satisfying the condition $\MDA$.

{\rm (i)} Assume that $\GC = GL(V)$ with $d \geq 5$, $\ell \neq 2,3$, and that the $G$-module
$V$ lifts to a complex module $\VC$. In addition, assume $C \leq G$ and
$\VC\da_{C} = A \oplus B$ with $A, B \neq 0$. Then either $\AC(\VC)$ enters $1_{C} \ua G$ with
multiplicity $\geq 2$, or $\begin{pmatrix}d+3\\4\end{pmatrix} \leq (G:C)$.

{\rm (ii)} Assume $V\da_{C} = A \oplus B \oplus B^{*}$ for some $C \leq G$ and $B \neq 0$.
Assume $\ell \neq 2$ if $\GC = O(V)$ and $\ell \neq 2,3$ if $\GC = Sp(V)$.
Assume in addition that $d \geq 4$ if $\GC = Sp(V)$, and $d \geq 9$ and
$\dim(B) > 1$ if $\GC = O(V)$. Then one of the following holds.

\hspace{0.5cm}{\rm (a)} $\begin{pmatrix}d+3\\4\end{pmatrix} < (G:C)$, $\ell \neq 2,3$, and
$\GC = GL(V)$ or $Sp(V)$.

\hspace{0.5cm}{\rm (b)} $\begin{pmatrix}d\\4\end{pmatrix} < (G:C)$. Furthermore, either
$\GC = GL(V)$ and $\ell = 2,3$, or $\GC = O(V)$.}
\end{propo}

\begin{proof}
1) Consider the case of (i) and set $\GCC = GL(\VC)$. By Proposition \ref{v4}(i),
the $\GC$-module $\VD$ is semisimple and $\dim \End_{\GC}(\VD) = 24$
according to Lemma \ref{dim}. Hence the condition $\MDA$ implies that
$\dim \End_{G}(\VD) = 24$ as well. Since the $G$-module $\VD$ is
semisimple and $V$ lifts to $\VC$, this in turn implies that
$\dim \End_{G}(\VC^{\otimes 4}) = 24$. Thus $G$ satisfies the condition $M_{8}(\GCC)$.
In particular, $G$ is irreducible on $\SD(\VC)$. Assume $\dim(A) = 1$. Then
$\SD(\VC)\da_{C}$ contains the submodule $\SD(A)$ of dimension $1$, so by Frobenius'
reciprocity $\SD(\VC)$ is a quotient of $\SD(A) \ua G$, whence
$\begin{pmatrix}d+3\\4\end{pmatrix} \leq (G:C)$. So we may assume $\dim(A), \dim(B) \geq 2$.

We also have $\dim \End_{G}(W) = 24$ for $W := \VC^{\otimes 2} \otimes (\VC^{*})^{\otimes 2}$.
Clearly, $(\VC \otimes \VC^{*})\da_{C} = M \oplus N \oplus K \oplus K^{*} \oplus 2 \cdot 1_{C}$,
where $A \otimes A^{*} = M \oplus 1_{C}$, $B \otimes B^{*} = N \oplus 1_{C}$,
and $K = A \otimes B^{*}$. Since $\dim(A), \dim(B) \geq 2$, we have
$\dim(M), \dim(N), \dim(K) \geq 3$ and $M$, $N$ are self-dual. Next, $W\da_{C}$ contains
$4 \cdot 1_{C} \oplus M \otimes M \oplus N \otimes N \oplus 2 \cdot K \otimes K^{*}$, whence
$W\da_{C}$ contains $8 \cdot 1_{C}$. According to Lemma \ref{nolan2},
$W\da_{\GC} = 2 \cdot 1_{\GC} \oplus 4 \cdot \AC(\VC) \oplus (\oplus^{4}_{i=1}B_{i})$ for some
$\GC$-modules $B_{i}$ of dimension larger than $\begin{pmatrix}d+3\\4\end{pmatrix}$. Since
$\dim \End_{G}(W) = 24$, $G$ is irreducible on $\AC(\VC)$ and on each $B_{i}$. If $1_{C}$ enters
$B_{i}$ for some $i$, then by Frobenius' reciprocity $B_{i}$ enters $1_{C} \ua G$ and so
$\dim(B_{i}) \leq (G:C)$. Otherwise $1_{C}$ enters $\AC(\VC)$ with multiplicity $> 1$, and so
we are done.

2) Next we consider the case of (ii). Set $Z = S$ if $(\GC,\ell) = (GL(V),\neq 2,3)$ or
$\GC = Sp(V)$, and $Z = \wedge$ otherwise. By Lemma \ref{weight} and Proposition \ref{v4},
$\GC$ is irreducible on $Z^{4}(V)$. Next, $\MDA$ implies $G$ is irreducible on $Z^{4}(V)$. On
the other hand, $Z^{4}(V)\da_{C}$ contains the nonzero submodule
$Z^{2}(B) \otimes Z^{2}(B^{*})$, so it contains $1_{C}$. By Frobenius'
reciprocity, $Z^{4}(V)$ is a (nontrivial) irreducible quotient of $1_{C} \ua G$, whence
$\dim(Z^{4}(V)) < (G:C)$.
\end{proof}

\begin{remar}\label{v5}
{\em Statement (i) of Remark \ref{v4e} and the argument in p.1) of the proof of Proposition
\ref{bound2} also show the following partial reduction to complex case. {\sl Assume $\ell > k$,
$G$ a finite group, the $\FF G$-module $V$ is not self-dual, $\dim(V) > k$, and $V$ lifts to a 
complex module $\VC$. Let $\GC = GL(V)$ and $\GCC = GL(\VC)$. Then $G$ is irreducible on all 
$\GC$-composition factors of $V^{\otimes k}$ if and only if $G$ is irreducible on all 
$\GCC$-composition factors of $\VC^{\otimes k}$. In particular, if $M_{2k}(\GC)$ holds for $V$ then 
$M_{2k}(\GCC)$ holds for $\VC$.}}
\end{remar}

\section{Reduction Theorems}\label{red-thms}
In this section we provide the reduction to the case where a finite group $G$ is either
nearly simple or the normalizer of a group of symplectic type. We begin with the following
proposition.

\begin{propo}\label{tensor}
{\sl Let $\FF$ be an algebraically closed field of characteristic $\ell > 0$, $V = \FF^{d}$ with
$d \geq 5$, and let $\GC$ be a classical group on $V$. Let $X = \AC$ if $\GC = GL(V)$,
$X = \TSB$ if $\GC = Sp(V)$ and $\ell \neq 2$, and $X = \TWB$ otherwise. Assume that
$G < \GC$ is a finite subgroup such that $G$ is irreducible on $X(V)$. Assume $G$ has 
a normal quasi-simple subgroup $S \in Lie(\ell)$. Then

{\rm (i)} $S$ is not of types $\ta B_{2}$, $\ta G_{2}$, or $\ta F_{4}$;

{\rm (ii)} $V\da_{S}$ is restricted (up to a Frobenius twist).}
\end{propo}

\begin{proof}
1) First we show that $V$ is irreducible over $S$. Assume the contrary. Then p. 2) of the
proof of Proposition \ref{irred1} shows that either $V\da_{S} = eU$ for some $\FF S$-irreducible
module $U$ and $e \geq 2$ or all irreducible constituents of $V\da_{S}$ are of dimension $1$. The
latter case is impossible as $S$ is perfect. Assume we are in the former case. If $X = \AC$ or if
$\ell = 2$, then the proof of Lemma \ref{prim} shows that $X(V)\da_{S}$ contains $1_{S}$.
Assume $\ell \neq 2$ and $X \neq A$. Then $V \simeq V^{*}$ and so $U \simeq U^{*}$. Observe that
$X(V) = \SB(V)$ or $\WB(V)$ here, so again $X(V)\da_{S}$ contains $1_{S}$. Thus in all cases
$X(V)\da_{S}$ contains $1_{S}$, whence $S$ must act trivially on $V$. The proof of Lemma \ref{prim}
now yields that $S = O^{\ell}(S) \leq Z(G)$, a contradiction.

As a consequence, we see that $O_{\ell}(Z(S)) = 1$, so $S$ is a quotient of the finite group of Lie
type $\hat{S}$ of simply connected type, defined over $\FF_{q}$ with $q = \ell^{f}$. Let $\HC$ be
the algebraic group corresponding to $\hat{S}$. We can view $V$ as an $\hat{S}$-module and extend
it to an $\HC$-module.

2) Here we show that $S$ cannot be of types $\ta B_{2}$, $\ta F_{4}$. Assume the contrary.
Then $\ell = 2$. Decompose $V\da_{S} = V_{1} \otimes V_{2} \otimes \ldots \otimes V_{t}$
into nontrivial {\it tensor indecomposable} irreducible $S$-modules. Then all $V_{i}$ are
self-dual, and so $X(V)\da_{S}$ contains $\TWB(V)\da_{S}$. Observe that
$\WB(V)\da_{S}$ contains $\SB(V_{1}) \otimes \SB(V_{2}) \otimes \SB(V_{t-1}) \otimes \WB(V_{t})$.
But $\SB(V_{i})$ and $\WB(V_{i})$ contains $1_{S}$, and $\SB(V_{1})$ contains $V_{1}^{(2)}$.
It follows that $V_{1}^{(2)} \hra X(V)\da_{S}$. Similarly, $V_{i}^{(2)} \hra X(V)\da_{S}$ for all
$i$. Since $G$ is irreducible on $V$ and $S \lhd G$, we get $\dim(V_{i}) = e$ for all $i$, and
$e$ divides $\dim(X(V))$. Now $e = 4$ if $S = \ta B_{2}(q)$, and $e = 26$, $246$, $2048$, or
$4096$ if $S = \ta F_{4}(q)'$. Since $\dim(X(V)) = e^{2t} - \alpha$ or $e^{t}(e^{t}-1)/2 - \alpha$
with $\alpha = 1$ or $2$, we see that $e$ cannot divide $\dim(X(V))$ if $e^{t} > 4$. Since
$e^{t} = d > 4$ by assumption, we get a contradiction.

3) Next we show that $S$ cannot be of type $\ta G_{2}$.
Assume the contrary: $S = \ta G_{2}(3^{2a+1})$.
Then $\ell = 3$. Decompose $V = V_{1} \otimes V_{2} \otimes \ldots \otimes V_{t}$
into nontrivial tensor indecomposable irreducible $\HC$-modules. Then all $V_{i}$ are self-dual of
type $+$. So $V\da_{S}$ is of type $+$, whence $X = \AC$ or $\TWB$, and $X(V)\da_{S}$ contains
$\TWB(V)$. Moreover, each $V_{i}$ is a Frobenius twist of $L(\om)$ with the highest weight
$\om = b_{i}\om_{1}$ or $b_{i}\om_{2}$, and $b_{1} = 1$ or $2$. (As in 2), this is a special
feature of the algebraic group of type $G_{2}$ in characteristic $3$, and types $B_{2}$ and $F_{4}$
in characteristic $2$.) We will denote $L(b\om_{1}+b'\om_{2})$ by $L(b,b')$. We also assume that
$\om_{1}$ corresponds to a short root.

Suppose that some $b_{i}$ is equal to $2$. Without loss we may assume that $V_{1} = L(2,0)$.
Writing $V = V_{1} \otimes A$, we see that $\WB(V)$ contains $\WB(V_{1}) \otimes \SB(A)$. But $A$ is
of type $+$, so $X(V)$ contains $\TWB(V_{1})$. One can checks that the $\HC$-module $\TWB(V_{1})$
contains $L(2,1)$ (of degree $189$) and $L(1,1)$ (of degree $49$), both irreducible over $S$. It
follows that $V$ is reducible over $G$, a contradiction.

Thus all $b_{i}$ are equal to $1$. Suppose $t > 1$. Without loss we may assume that
$V_{1} = L(1,0)$, $V_{2} = L(3^{c},0)$ with $0 < c \leq 2a$. Notice that
$\WB(L(1,0)) = 2L(1,0) + L(0,1)$ and $\SB(L(1,0)) = L(0,0) + L(2,0)$ as $\HC$-modules. As above,
$\WB(V)$ contains $\WB(V_{1})$, and so $X(V)$ contains $L(1,0)$, which is of degree $7$ and
irreducible over $S$. Next, $\WB(V)$ contains $\SB(V_{1}) \otimes \WB(V_{2})$ and so it
contains $L(2,0) \otimes (L(3^{c},0) + L(0,3^{c}))$. If $0 < c \leq a$, then
$L(2,0) \otimes L(3^{c},0) \simeq L(3^{c}+2,0)$ is an $S$-irreducible of degree $189$.
If $a < c \leq 2a$, then $0 < c-a \leq a$, and
$L(2,0) \otimes L(0,3^{c}) \simeq L(2,0) \otimes L(3^{c-a},0) \simeq L(3^{c-a}+2,0)$ is an
$S$-irreducible of degree $189$. In all cases, $X(V)$ contains $S$-irreducibles of distinct degrees,
so it is again reducible over $G$.

Thus $t = 1$. Without loss we may assume that $V = L(1,0)$. In this case,
$\AC(V)$ contains $L(2,0)$ (of degree $27$) and $L(1,0)$ (of degree $7$), both irreducible over $S$,
so $\AC(V)$ is $G$-reducible. Furthermore, $\TWB(V) = 2L(1,0) + L(0,1)$ contains distinct
$S$-irreducibles with distinct multiplicities, so again it is reducible over $G$, a contradiction.

4) Now we can assume that $S$ is not of types $\ta B_{2}$, $\ta G_{2}$, or $\ta F_{4}$.
According to Steinberg's tensor product theorem, cf. \cite[Thm. 5.4.1]{KL},
$V = V_{1} \otimes V_{2} \otimes \ldots \otimes V_{t}$ as $\hat{S}$-modules, where
$V_{i} = U_{i}^{(\ell^{n_{i}})}$, $0 \leq n_{1} < n_{2} < \ldots < n_{t} < f$, and all $U_{i}$ are
nontrivial and restricted. Moreover, this decomposition is unique up to $S$-isomorphism.
Observe that any automorphism of $S$ transforms this decomposition into another tensor decomposition
satisfying the same constraints on the factors. Hence $G$ acts on the set
$\{V_{1}, V_{2}, \ldots ,V_{t}\}$ (up to $S$-isomorphism).

Next assume $X \in \{\TWB,\TSB\}$. Then $V = V^{*}$ as $S$-modules, so
$\{V_{1}, V_{2}, \ldots ,V_{t}\} = \{V_{1}^{*}, V_{2}^{*}, \ldots ,V_{t}^{*}\}$ by uniqueness of
the Steinberg decomposition, whence $V_{i} \simeq V_{i}^{*}$ as $S$-modules by the irreducibility
of $V\da_{S}$.

We need to show that $t = 1$. Assume the contrary: $t > 1$.

5) Let $e := \min\{\dim(V_{i}) \mid 1 \leq i \leq t\}$. Also, for each $i$ let $t_{i}$ be the number
of indices $j$ with $\dim(V_{j}) = \dim(V_{i})$.
Observe that $|\Out(S)| = af$, with $a \leq 3e$ (cf. \cite[\S5.1]{KL}). Let $M := Z(G)S$. Since
$V\da_{S}$ is irreducible, $G/M \hra \Out(S)$. Here we show that $|G/M| \leq at_{i}$. Indeed, suppose
$|G/M| > at_{1}$. Consider the subgroup $F$ of outer field automorphisms of $S$.
Then $|F| = f$. Since $|G/M| > at_{1}$ and $|\Out(S)| = af$, it follows that
$F' := F \cap G/M$ has order $> t_{1}$. On the other hand, since $U_{1}$ is restricted, no
nontrivial element in $F$ can fix $V_{1}$. Thus the $G$-orbit of $V_{1}$ (cf. 4)) has length $>t_{1}$,
contrary to the choice of $t_{1}$.

6) Assume that $\dim(V_{i}) \geq 3$ for some $i$. It is easy to see that $\AC(V_{i})$, resp.
$\WB(V_{i})$, $\SB(V_{i})$, contains some nontrivial $S$-irreducible say $A_{i}$, resp. $B_{i}$,
$C_{i}$. Here we show that $X(V)$ contains a nontrivial $S$-irreducible
$B \in \{A_{i},B_{i},C_{i}\}$; moreover, one can choose $B \in \{B_{i},C_{i}\}$ if all $S$-modules
$V_{j}$ are self-dual. Indeed, $(V \otimes V^{*})\da_{S}$ contains
$V_{i} \otimes V_{i}^{*}$, and $\WB(V_{i}) \hra V_{i} \otimes V_{i}^{*}$ if $V_{i} \simeq V_{i}^{*}$,
so the claims follows for $X = \AC$. Next assume all $S$-modules $V_{j}$ are self-dual. Then for
each $j \neq i$ we can choose $Y_{j} \in \{\SB,\WB\}$ such that $1_{S} \hra Y_{j}(V_{j})$.
Given any $Y \in \{\SB,\WB\}$, we can choose $Y_{i} \in \{\SB,\WB\}$ such that
$Y(V)\da_{S}$ contains $\otimes_{j \neq i}Y_{j}(V_{j}) \otimes Y_{i}(V_{i})$, and so
$Y(V)\da_{S}$ contains $B_{i}$ or $C_{i}$. Hence $X(V)\da_{S}$ contains $B_{i}$ or $C_{i}$.

7) Suppose $X = \TWB$ for instance. Then according to 4) all the $S$-modules $V_{j}$ are self-dual.
We may assume $\dim(V_{1}) = e$.

First we consider the case $e \geq 3$ and $t \geq 3$. According to 6), $X(V)$ contains an
$S$-irreducible $B \in \{B_{1},C_{1}\}$; in particular, $\dim(B) \leq e(e+1)/2$.  Clearly, $B$
extends to an $M$-irreducible which we also denote by $B$, and then $X(V) \hra (B \ua G)$. It follows
that
$$e^{t}(e^{t}-1)/2 -2 \leq \dim(X(V)) \leq (G:M) \cdot \dim(B) \leq 3et \cdot e(e+1)/2,$$
a contradiction since $e,t \geq 3$.

The same argument applies to the case $e \geq 7$ and $t = 2$. If $e = 6$, then notice that
$|\Out(S)| = af$ with $a \leq 2e$, so $(G:M) \leq 2et$, and we can use the same argument. Assume
$3 \leq e \leq 5$, $t = 2$, and $\dim(V_{2}) > e$. Since the $G$-orbit of $V_{1}$
has length $1$, $(G:M) \leq 3e$ by 5). On the other hand, $\dim(V) \geq e(e+1)$, so
$\dim(X(V)) > (G:M) \cdot \dim(B)$, a contradiction. Next assume $3 \leq e \leq 5$, $t = 2$,
and $\dim(V_{1}) = \dim(V_{2})$. Since $S$ has a self-dual irreducible module $V_{1}$ of dimension
$e$ with $3 \leq e \leq 5$, $S = PSp_{4}(q)$, and we may assume that both $V_{1}$ and $V_{2}$ are
(different) Frobenius twists of $L(\om_{i})$ for $i = 1$ or $2$. Direct calculation shows that in
these subcases $X(V)$ contains $S$-irreducibles of distinct dimension, so $X(V)$ is $G$-reducible.

Now we may assume $e = 2$; in particular, $S = L_{2}(q)$ and $|\Out(S)| = af$ with $a \leq 2$.
Assume that some $V_{i}$ has dimension $g > 2$ and set $s := t_{i}$. Then $\dim(V) \geq 2g^{s}$.
According to 6), $X(V)$ contains an $S$-irreducible $B \in \{B_{i},C_{i}\}$, and
$\dim(B) \leq g(g+1)/2$. According to 5), $(G:M) \leq 2s$. As above, $B$ extends to
$M$ and $X(V) \hra (B \ua G)$, so we get
$$\dim(X(V)) \geq g^{s}(2g^{s}-1) -2 > 2s \cdot g(g+1)/2 \geq (G:M) \cdot \dim(B),$$
a contradiction.

So now we can assume $\dim(V_{1}) = \ldots = \dim(V_{t}) = 2$, and $t \geq 3$ as $d = \dim(V) > 4$.
Also, $(G:M) \leq 2t$. Assume $t$ is even. Then $\WB(V)$ contains
$\SB(V_{1}) \otimes \WB(V_{2}) \otimes \ldots \otimes \WB(V_{t}) \simeq \SB(V_{1})$, whence
$X(V)$ contains an $S$-irreducible $B$ of dimension $2$ or $3$. It follows that 
$$2^{t-1}(2^{t}-1) -2 \leq \dim(X(V)) \leq (G:M) \cdot \dim(B) \leq 2t \cdot 3,$$
a contradiction since $t \geq 3$. Assume $t$ is odd. Since $X = \TWB$ is not considered when
$\GC = Sp(V)$ and $\ell \neq 2$, we must have $\ell = 2$. In this case
$\SB(V_{1}) = 1_{S} + B$ for some $S$-irreducible $B$ of dimension $2$. Furthermore, $\SB(V_{i})$
and $\WB(V_{i})$ contain $1_{S}$. Since $\WB(V)$ contains
$\SB(V_{1}) \otimes \SB(V_{2}) \otimes \WB(V_{3}) \otimes \ldots \otimes \WB(V_{t})$, we see
that $X(V)\da_{S}$ contains $B$, and again the inequality
$\dim(X(V)) \leq (G:M) \cdot \dim(B)$ yields a contradiction.

The cases $X = \AC$ or $X = \TSB$ can be handled similarly.
\end{proof}

\begin{remar}
{\em (i) Proposition \ref{tensor} fails if $d = 4$. For consider $S = L_{2}(q^{2})$ with $q$ odd,
and $G$ is $S$ extended by the field automorphism of order $2$. If $A$ is the natural
$2$-dimensional module for $SL_{2}(q^{2})$, then $V := A \otimes A^{(q)}$ is a $G$-stable
$S$-module, so $V$ can be extended to $G$ in such a way that $G < \GC = O(V)$. One checks
that $\WB(V)$ is irreducible over $G$. However, $V\da_{S}$ is tensor decomposable and not
restricted.

(ii) An analogue of Proposition \ref{tensor} fails for $\GC = Sp(V)$ and $X = \TWB$. For consider
$S = SL_{2}(q^{3})$ with $q$ odd, and $G$ is $S$ extended by the field automorphism of order $3$.
If $A$ is the natural $2$-dimensional module for $S$, then
$V := A \otimes A^{(q)} \otimes A^{(q^{2})}$ is $G$-stable, so $V$ can be extended to $G$ in such
a way that $G < \GC = Sp(V)$. One checks that $\TWB(V)$ is irreducible over $G$. However, $V\da_{S}$
is tensor decomposable and not restricted.

(iii) The above examples also show that the $\FF G$-module $V$ may be tensor indecomposable
over $G$ but yet tensor decomposable over $S = F^{*}(G)$. Here is a similar example but in
cross characteristic. Consider $S = Sp_{2n}(q)$ with $q = 3,5$ and $G$ is $S$ extended by
the diagonal outer automorphism of order $2$. Then $S$ has complex Weil representations
$A$, $B$ of degree $(q^{n}-1)/2$ and $(q^{n}+1)/2$, respectively, such that
$A \otimes B$ is irreducible and $G$-stable, cf. \cite[Prop. 5.4]{MT1}. Now the $\CC G$-module
$A \otimes B$ is tensor indecomposable over $G$ but not over $S$. Thus it can happen that 
a (nearly simple) group $G$ does not preserve a tensor decomposition of a space $W$ but $F^{*}(G)$ 
does; in particular, $F^{*}(G)$ is contained in a positive dimensional subgroup of $GL(W)$ (the
stabilizer of the tensor decomposition) but $G$ is not.}
\end{remar}

\begin{theor}\label{red-sl}
{\sl Let $\FF$ be an algebraically closed field of characteristic $\ell \geq 0$. Let
$\XC = GL(V) = GL_{d}(\FF)$, $d > 2$, and let $G$ be a closed subgroup of $\XC$. Assume
that $G$ is irreducible on the nontrivial composition factor $\AC(V)$ of $V \otimes V^{*}$.
Then one of the following holds:

{\rm (i)} $SL(V) \lhd G$;

{\rm (ii)} $\ell > 0$, $SL_{d}(q)' \lhd G$ or $SU_{d}(q)' \lhd G$ for a power $q$ of $\ell$;

{\rm (iii)}  $F^{*}(G) = Z(G)E$, where $E$ is extraspecial of order $p^{1+2a}$ for some
prime $p \neq \ell$ and $d=p^{a}$;

{\rm (iv)} $G$ is finite and nearly simple with the nonabelian composition factor not a
Chevalley group in characteristic $\ell$.}
\end{theor}

\begin{proof}
First suppose that $G$ is contained in $\HC$, a proper closed subgroup of $\XC$ of positive
dimension. Then any nontrivial irreducible composition factor of the adjoint module for the
connected component $\HC$ embeds in the irreducible part $\AC(V)$ of the adjoint module for $\XC$.
A straightforward dimension computation shows that $\HC$ is not irreducible on $\AC(V)$,
(since $\dim(\AC(V)) = d^{2} - 1 - \al$, where $\al =1$ if $\ell$ divides $d$ and $0$ otherwise).

By Aschbacher's classification of subgroups of $GL(V)$ \cite{A}, the only possibilities for $G$
not contained in a positive dimensional closed (proper) subgroup are normalizers of groups of
symplectic type or nearly simple groups.

All that remains to show is to consider the case $G$ is nearly simple and its composition factor
is a finite simple group $\BS \in Lie(\ell)$. Moreover, the representation
must be tensor indecomposable (more strongly, preserves no tensor structure) and primitive. Notice
that $S = E(G)$ is a quotient of the finite group of Lie type $\hat{S}$ of simply connected
type. Let $\HC$ be the algebraic group corresponding to $\hat{S}$. We can view $V$ as an
$\hat{S}$-module and extend it to an $\HC$-module. Using \cite[Prop. 5.4.3]{KL}, one can see
that the self-duality of the $S$-module $V$ implies that $V \simeq V^{*}$ as $\HC$-modules.
It follows that if $S$ supports some nondegenerate bilinear form $\bb$ on $V$ then $\bb$ is
also $\HC$-invariant (we will need this observation for proving further theorems in this section).
By Proposition \ref{tensor}, $S$ is not a Suzuki or a Ree group, and we may assume that $V\da_{S}$
is restricted. In particular, $G$ cannot induce any (nontrivial) field automorphism on $S$.

We claim that $G$ preserves the adjoint module for $\HC$ which is contained in the adjoint module for
$\GC$. Assuming this claim, the irreducibility of $G$ on $\AC(V)$ implies that $\HC$ must be $SL(V)$.
Thus $G$ is essentially either $SL_{d}(q)$ or $SU_{d}(q)$; more precisely, (ii) holds. So we will
prove the claim, and we give a proof which also works in the cases $\GC = Sp(V)$
or $O(V)$.

The representation of $G$ and $\HC$ on $V$ induces a homomorphism $\Phi~:~G \to \GC$ and
$\Psi~:~\HC \to \GC$, with $\Phi = \Psi$ on $S$. Assume some $g \in G$ induces an outer diagonal
automorphism on $S$. Then we can find $h \in \HC$ such that $g$ and $h$ induce the same automorphism
on $S$ (via conjugation). Since $V\da_{S}$ is irreducible, $\Phi(g) = \lambda \Psi(h)$ for some
$\lambda \in \FF^{*}$ by Schur's lemma. It follows that $\Phi(g)$ normalizes $\Psi(\HC)$.
Now assume $g \in G$ induces a nontrivial graph automorphism on $S$. Then there is a graph
automorphism $\tau$ of $\HC$ such that $\tau(x) = gxg^{-1}$ for all $x \in S$. Since the $S$-module
$V$ is $g$-stable, its highest weight is also $\tau$-stable, whence the $\HC$-module $V$ is
$\tau$-stable. Thus $\Psi(\tau(y)) = M\Psi(y)M^{-1}$ for some $M \in GL(V)$ and for all $y \in \HC$.
Also, for all $x \in S$ we have
$$\Phi(g)\Phi(x)\Phi(g)^{-1} = \Phi(gxg^{-1}) = \Psi(gxg^{-1}) = \Psi(\tau(x)) = M\Psi(x)M^{-1}
  = M \Phi(x) M^{-1}$$
so again $\Phi(g) = \gamma M$ for some $\gamma \in \FF^{*}$ by Schur's lemma. Since
$\Phi(g)\Psi(y)\Phi(g)^{-1} = M\Psi(y)M^{-1} = \Psi(\tau(y))$ for all $y \in \HC$, we see that
$\Phi(g) \in N_{\GC}(\Psi(\HC))$. We have shown that $\Phi(G) < N_{\GC}(\Psi(\HC))$. Since
$N_{\GC}(\Psi(\HC))$ preserves the adjoint module for $\HC$, the same is true for $\Phi(G)$, and
the claim is established.
\end{proof}

Next consider $\XC = O(V) = O_{d}(\FF)$ with $d \geq 5$. We assume that $\XC$ is irreducible on
$V$, i.e. if $\ell = 2$, $n$ is even. The adjoint module for $\XC$ is essentially $\WB(V)$ (if
$\ell = 2$, then there is a trivial composition factor when $n \equiv 2 (\mod 4)$ and two
trivial composition factors when $4|n$). The proof of the previous theorem yields the following
result.

\begin{theor}\label{red-so}
{\sl Let $\FF$ be an algebraically closed field of characteristic $\ell \geq 0$.  Let
$\XC = O(V) = O_{d}(\FF)$, $d \geq 5$, and let $G$ be a closed subgroup of $\XC$.
If $\ell = 2$, assume that $n$ is even. Assume that $G$ is irreducible on $\TWB(V)$.
Then one of the following holds:

{\rm (i)} $\Om(V) \lhd G$;

{\rm (ii)} $\ell > 0$, and $\Om^{\pm}_{d}(q) \lhd G$ with $q$ a power of $\ell$;

{\rm (iii)} $\ell > 0$, $d = 8$, and $\tb D_{4}(q) \lhd G$ with $q$ a power of $\ell$;

{\rm (iv)} $\ell \neq 2$, and $G$ is contained in the stabilizer $\ZZ_{2}^{d}:\SSS_{d}$ of an
orthonormal basis of $V$ in $\XC$;

{\rm (v)} $\ell \neq 2$, $F^{*}(G) = E$, where $E = 2^{1+2a}_{+}$ is extraspecial of type $+$
and $d = 2^{a}$;

{\rm (vi)} $G$ is finite and nearly simple with the nonabelian composition factor not a
Chevalley group in characteristic $\ell$.
\hfill $\Box$}
\end{theor}

Note that the possibility (iv) arises exactly when $G$ is imprimitive on $V$, cf.
Proposition \ref{irred1}. Also, the assumption $d \geq 5$ is necessary, see the remark in p.1) of
the proof of Proposition \ref{irred1}.

We next consider the symplectic case.  If $\ell \neq 2$, then the same argument as above yields:

\begin{theor}\label{red-sp}
{\sl Let $\FF$ be an algebraically closed field of characteristic $\ell \neq 2$. Let
$\XC = Sp(V) = Sp_{2m}(\FF)$, $m > 2$, and let $G$ be a proper closed subgroup of $\XC$. Assume
that $G$ is irreducible on $\TSB(V)$.  Then one of the following holds:

{\rm (i)} $\ell > 0$, and $Sp_{2m}(q) \lhd G$ with $q$ a power of $\ell$;

{\rm (ii)} $F^{*}(G) = E$, where $E = 2^{1+2a}_{-}$ is extraspecial of type $-$, and $2m = 2^{a}$;

{\rm (iii)} $G$ is finite and nearly simple with the nonabelian composition factor not a
Chevalley group in characteristic $\ell$.
\hfill $\Box$}
\end{theor}

Finally, assume $\ell = 2$ and $\XC = Sp(V) = Sp_{2m}(\FF)$. Then the largest composition factor
of the adjoint module of $\XC$ is $\TWB(V)$, other composition factors being the Frobenius
twist $\Vbb$ and the trivial ones. If $\HC$ is a positive dimensional closed subgroup
of $\XC$, then the heart of the adjoint module for $\HC$ is a proper subquotient for $\TWB(V)$
(again compute dimensions -- the only case where the dimension is sufficiently large is the case
of $D_{n}$ contained in $C_{n}$ or $G_{2}$ inside $C_{3}$). Now if $G$ is contained in some
orthogonal group $O(V)$, Theorem \ref{red-so} applies. Assume $G$ is contained in
$\HC = G_{2}(\FF)$ and irreducible on $\TWB(V)$. Again, if $G$ normalizes some proper positive
dimensional closed subgroup $\HC_{1}$ of $\HC$, then $G$ stabilizes the adjoint module for $\HC_{1}$
which is a proper submodule of $\TWB(V)$. Thus $G$ is a Lie primitive finite subgroup of $\HC$.
In this case one can show (cf. \cite{LiS}) that $G$ is as in cases (vi) or (vii) of Theorem
\ref{red-sp2}. Finally, if $G$ is a finite subgroup of $\GC$ then one easily reduces to the case
of nearly simple groups.

\begin{theor}\label{red-sp2}
{\sl Let $\FF$ be an algebraically closed field of characteristic $2$. Let
$\XC = Sp(V) = Sp_{2m}(\FF)$, $m > 2$, and let $G$ be a proper closed subgroup of $\XC$. Assume
that $G$ is irreducible on $\TWB(V)$.  Then one of the following holds:

{\rm (i)} $Sp_{2m}(q) \lhd G$ with $q$ a power of $2$;

{\rm (ii)} $\Om(V) \lhd G$ for some quadratic form on $V$;

{\rm (iii)} $\Om^{\pm}_{2m}(q) \lhd G$ with $q$ a power of $2$;

{\rm (iv)} $m = 4$, and $\tb D_{4}(q) \lhd G$ with $q$ a power of $2$;

{\rm (v)} $m = 3$, and $G = G_{2}(\FF)$;

{\rm (vi)} $m = 3$, and $F^{*}(G) = G_{2}(q)'$ with $q$ a power of $2$;

{\rm (vii)} $G$ is finite and nearly simple with the nonabelian composition factor not a
Chevalley group in characteristic $2$.
\hfill $\Box$}
\end{theor}

Notice $\ta B_{2}(q)$ is excluded as we assume $d \geq 5$.

\section{Normalizers of Symplectic Type Subgroups}\label{extra}

Here we consider the case $F^{*}(G) = Z(G)E$ is of symplectic type, i.e. $E$ is either 
extraspecial of odd exponent $p$, an extraspecial $2$-group or a central product of an 
extraspecial $2$-group with a cyclic group of order $4$ (with the central involutions identified).
If $E$ is an extraspecial $2$-group, we say $E$ is of $+$ type if $E$ is a central product of
$D_{8}$'s and of $-$ type if $E$ is a central product of $Q_{8}$ and some number of $D_{8}$'s.

If $E$ is extraspecial of order $p^{1+2a}$, then an irreducible faithful module for $E$ has 
dimension $p^{a}$ and is unique once the character of $Z(E)$ is fixed. Moreover, if 
$E \lhd G \leq \GC$ then $(\ell,|E|) = 1$ and one of the following holds, where we denote 
$Z := Z(\GC)$:

(a) $p$ is odd, $G \leq N := (EZ):Sp_{2a}(p)$ and $\XC=GL(V)$;

(b) $p=2$, $G \leq N := (EZ) \cdot Sp_{2a}(p)$ and $\XC=GL(V)$;

(c) $p=2$, $G \leq  N := E \cdot O^{+}_{2a}(2)$ and $\XC=O(V)$;

(d) $p=2$, $G \leq  N := E \cdot O^{-}_{2a}(2)$ and $\XC=Sp(V)$.

We now assume that $E \lhd G \leq N$, $|E| = p^{1+2a}$, $d = \dim(V) = p^{a} > 4$,
$W := \FF_{p}^{2a}$ the natural module for $N/EZ$.

\begin{lemma}\label{ex-irr}
{\sl $G$ satisfies $\MBA$ if and only if one of the following holds.

{\rm (i)} $\GC = GL(V)$ and $G/E(Z \cap G)$ is a subgroup of $Sp_{2a}(p)$ that is transitive on
$W^{\#} := W \setminus \{0\}$ (such subgroups are classified by Hering's theorem, cf. 
\cite{He}, \cite{Li}).

{\rm (ii)} $\GC = O(V)$ and $G/E$ is a subgroup of $O^{+}_{2a}(2)$ that is transitive on
the isotropic vectors and on the non-isotropic vectors of $W^{\#}$ (such subgroups are classified 
in \cite{Li}).

{\rm (iii)} $\GC = Sp(V)$ and $G/E$ is a subgroup of $O^{-}_{2a}(2)$ that is transitive on
the isotropic vectors and on the non-isotropic vectors of $W^{\#}$ (such subgroups are classified 
in \cite{Li}).}
\end{lemma}

\begin{proof}
(i) Let $\varphi$ be the Brauer character of $G$ on $\AC(V)$. Since $\ell$ is coprime to
$\dim(V) = p^{a}$, $\varphi \da_{E} = \sum_{\lam \in \Irr(E/Z(E))}\lam - 1_{E}$. By Clifford's
theorem, $G$ is irreducible on $\AC(V)$ if and only if $G/E$ is transitive on
$\Irr(E/Z(E)) \setminus \{1_{E}\}$, equivalently, on $W^{\#}$. Observe that $H := G/ZE \leq Sp(W)$ 
is such that $WH$ is doubly transitive on $W$; such subgroups are classified by 
Hering's theorem, cf. \cite{He}, \cite{Li}. 

(ii) Using the alternating form on $W$ and identifying $E/Z(E)$ with $W$, we can identify
$\Irr(E/Z(E))$ with $W$. Then one can check that $\TSB(V)$ and $\TWB(V)$ afford the
$E$-characters
$$\sum_{v \in W^{\#},~v \mbox{ {\tiny isotropic}}}v,~~~~\mbox{resp.}~~~
  \sum_{v \in W^{\#},~v \mbox{ {\tiny nonisotropic}}}v.$$
Now the claim follows by Clifford's theorem. Also observe that $H := G/E \leq O(W)$ is such that 
$WH$ is an affine permutation group of rank $3$ on $W$; such subgroups are classified by Liebeck 
\cite{Li}. 

(iii) is similar to (ii).
\end{proof}

From now on we assume that $G \leq N$ satisfies $\MBA$ and aim to show that $G$ fails
$\MCA \cap \MDA$. So there is no loss in assuming $G = N$. The following two statements deal 
with the complex case.

\begin{propo}\label{ex-m8}
{\sl Assume $\ell = 0$ and $p^{a} > 4$. 

{\rm (i)} Assume $\GC = GL(V)$. Then $M_{6}(N,V) - M_{6}(\GC,V) \geq 2p-5$. If $p = 2$ then 
$M_{8}(N,V) - M_{8}(\GC,V) = 6$.

{\rm (ii)} Assume $p = 2$, $a \geq 4$, and $\GC = Sp(V)$ or $O(V)$. Then $M_{6}(N,V) = M_{6}(\GC,V)$
and $M_{8}(N,V) - M_{8}(\GC,V) = 20$.}
\end{propo}

\begin{proof}
(i) Observe that $M := V \otimes V^{*}$ is trivial on $Z$, and, considered as a module over $N/Z$, it
is the permutation module on $W$ with $E/Z(E)$ acting by translations and $N/EZ \simeq Sp(W)$ acting
naturally. So $M$ affords the character $\rho := \sum_{v \in W}v$, where we again identify
$\Irr(E/Z(E))$ with $W$ as in the proof of Lemma \ref{ex-irr}. It follows that $M^{\otimes 3}$ is 
the permutation module on $W \times W \times W$, and that the fixed point
subspace for $ZE$ inside $M^{\otimes 3}$ affords the $E/Z(E)$-character 
$1_{E} \cdot \left(\sum_{u,v,w \in W, u + v + w = 0}1\right)$. On the triples
$(u,w,w)$, $u,v,w \in W$, $u+v+w = 0$, $Sp(W)$ acts with exactly $2p+2$ orbits if $a > 1$ 
and $2p+1$ orbits if $a = 1$. The orbit representatives are 
$(0,0,0)$; $(u,0,-u)$, $(u,-u,0)$, $(0,u,-u)$ with $u \neq 0$; 
$(u,\lam u,(-1-\lam) u)$ with $u \neq 0$, $\lam \in \FF_{p} \setminus \{0,-1\}$; and  
$p$ orbits $(u,v,-(u+v))$ with $u,v \in W$ linearly independent and the inner product 
$(u|v) = \mu \in \FF_{p}$ (and $\mu \neq 0$ if $a = 1$). Each orbit gives rise to a 
permutation module for $Sp(W)$, so it yields a trivial character $1_{N}$. It follows that 
$M^{\otimes 3}$ contains $1_{N}$ with multiplicity at least $2p+1$, while 
$M_{6}(\GC,V) = 6$.

Next we consider the case $p = 2$ and $\GC = GL(V)$; in particular $a \geq 3$. Then $M^{\otimes 4}$ 
is the permutation module on $W \times W \times W \times W$, and the fixed point
subspace for $ZE$ inside $M^{\otimes 4}$ affords the $E/Z(E)$-character 
$1_{E} \cdot \left(\sum_{u,v,w,t \in W, u + v + w +t = 0}1\right)$. On the quadruples
$(u,w,w,t)$, $u,v,w,t \in W$, $u+v+w+t = 0$, $Sp(W)$ acts with exactly 
$$1 + 6 + 4 + 4 + 1 + 3 + 3 + 8 = 30$$ 
orbits. The orbit representatives are 
$(0,0,0,0)$; $(u,u,0,0)$ with $u \neq 0$ and $5$ other permutations of this;
$(u,v,u+v,0)$ with $(u|v) = 1$ and $3$ other permutations of this;
$(u,v,u+v,0)$ with $(u|v) = 0$ but $u \neq v$ and $3$ other permutations of this;
$(u,u,u,u)$ with $u \neq 0$; $(u,v,u,v)$ with $(u|v) = 1$ and $2$ other permutations of this;
$(u,v,u,v)$ with $(u|v) = 0$ but $u \neq v$ and $2$ other permutations of this;
and $(u,v,w,u+v+w)$ with $u,v,w$ linearly independent and the Gram matrix of $(\cdot|\cdot)$ on 
$\langle u,v,w \rangle_{\FF_{2}}$ being $\begin{pmatrix}0 & x & y\\x & 0 & z \\y & z & 0\end{pmatrix}$,
$x,y,z \in \FF_{2}$. Again, each orbit gives rise to a trivial character $1_{N}$. It follows that 
$M^{\otimes 4}$ contains $1_{N}$ with multiplicity $30$, while $M_{8}(\GC,V) = 24$. 
 
(ii) Finally, we consider the case where $p = 2$, $a \geq 4$, and $\GC = O(V)$ or $Sp(V)$. Then we can 
repeat the computations in (i), except that now we have to take into account the values of the 
quadratic form $Q$ on $W = \FF_{2}^{2a}$. With this in mind, one can show that the total number of 
$O(W)$-orbits on the the triples $(u,w,w)$, $u,v,w \in W$, $u+v+w = 0$, is 
$1 + 2 \cdot 3 + 2 \cdot 4 = 15$, yielding $M_{6}(N,V) = 15$. Similarly, the total number of 
$O(W)$-orbits on the the quadruples $(u,w,w,t)$, $u,v,w,t \in W$, $u+v+w+t = 0$, is  
$$1 + 2 \cdot 6 + 4 \cdot 4 + 4 \cdot 4 + 2 \cdot 1 + 4 \cdot 3 + 4 \cdot 3 + 8 \cdot 8 = 135,$$
yielding $M_{8}(N,V) = 135$.  
\end{proof}

\begin{propo}\label{ex-m62}
{\sl Assume $\ell = 0$ and $p^{a} > 4$. Then $M_{6}(\GC,V) = M_{6}(G,V)$ if and only if $G$ is 
as described in the case {\rm (B)} of Theorem \ref{complex-m6}.}
\end{propo}

\begin{proof}
By Proposition \ref{ex-m8} we may assume $p = 2$. In fact the proof of Proposition \ref{ex-m8}
establishes the ``if'' part of our claim. It also shows that, for the ``only if'' part, it suffices
to prove that if $H := G/Z(G)E$ has the same orbits on $W \times W$ as of $M := N/Z(N)E$ then 
$H \geq [M,M]$. Since $\MC$ implies $\MB$, $G$ satisfies the conclusions of Lemma \ref{ex-irr}.  

First we consider the case $\GC = GL(V)$. Then $M = Sp(W)$ and $M$ has an orbit 
$\Delta = \{(u,v) \mid 0 \neq u,v \in W, u \neq v, (u|v) = 0\}$ of length $(d^{2}-1)(d^{2}-4)/2$ 
on $W \times W$, and the longest orbit of $M$ on $W \times W$ has length $d^{2}(d^{2}-1)/2$. By 
Hering's theorem \cite{He}, \cite{Li}, $H$ is one of the following four possibilities. 
The transitivity of $H$ on $\Delta$ eliminates the possibilities $H \rhd G_{2}(q)'$ with 
$q^{6} = 2^{2a}$ and $H \leq \Gamma L_{1}(2^{2a})$. Assume $H \rhd H_{0}:= Sp_{2b}(q)$ with 
$q = 2^{c} > 2$ and $a = bc$. If $b > 1$, then $H_{0}$ has orbits of different lengths 
$d^{2}-1$ and $(d^{2}-1)(q^{2b-1}-q)$ on $\Delta$, whence $H_{0}$ is not semi-transitive on $\Delta$
and $H$ cannot be transitive on $\Delta$, a contradiction. If $b = 1$ then $H$ is too small to 
be transitive on $\Delta$. Assume $H \rhd H_{0}:= SL_{b}(q)$ with $q = 2^{c}$, $2a = bc$ and $b > 2$. 
Then $H_{0}$ has an orbit of length $(d^{2}-1)(d^{2}-q) > d^{2}(d^{2}-1)/2$ on $W \times W$, again a 
contradiction. Thus $H = Sp_{2a}(2)$ as stated. 

Next we consider the cases $(\GC,\eps) = (Sp(V),-)$ or $(O(V),+)$. Then $M = O(W)$. Assume the 
contrary: $H \not\geq [M,M]$. Since $H$ has the same orbits on $W$ as of $M$, the main result of 
\cite{Li} implies that one of the following four possibilities occurs. Notice that $\MC$ implies 
that $G$ has an irreducible constituent of degree $d(d^{2}-4)/3$ on $V^{\otimes 3}$. This condition
excludes the first possibility $H \leq \Gamma L_{1}(2^{2a})$. Assume $H \rhd SU_{a}(2)$. If $2|a$, 
then we can find a primitive prime divisor of $2^{a-1}-1$ \cite{Zs} which divides $d(d^{2}-4)/3$ 
but not $|G|$, a contradiction.  If $a$ is odd, then a primitive prime divisor of $2^{2a-2}-1$ again 
divides $d(d^{2}-4)/3$ but not $|G|$. Assume $H \rhd Sp_{6}(2)$; in particular $\eps = +$. Since $H$ 
is an irreducible subgroup of $M = O^{+}_{8}(2)$, $H = Sp_{6}(2)$. Now $H$ acts on $W$ with orbits
of length $1$, $120$ and $135$, so the permutation character $\rho$ of $H$ on $W$ is uniquely
determined by \cite{Atlas}. Hence we can check that $(\rho,\rho) = 16$, meaning $H$ has $16$ orbits
on $W \times W$, one more than $M$ does, a contradiction. Finally, assume 
$H \leq N_{\Gamma L_{3}(4)}(3^{1+2})$; in particular, $\eps = -$. Here the permutation character
of $M = O^{-}_{6}(2)$ on $W$ is 
$3 \cdot \overline{1} + \overline{6} + \overline{15} + 2 \cdot \overline{20}$ (\cite{Atlas}), where
$\overline{m}$ denotes an irreducible character of degree $m$. One can show that $(|H|,5) = 1$, 
whence $\overline{20}|_{H}$ is reducible and so $H$ has at least $19$ orbits on $W \times W$, again 
a contradiction.         
\end{proof}

\subsection{Case $\GC = GL(V)$.}

We begin with the case $\GC = GL(V)$. 
\begin{lemma}\label{ex-p2}
{\sl Assume $p=2$ and $a > 2$. Then $G$ fails $\MDA$.}
\end{lemma}

\begin{proof}
Assume the contrary. Then $G$ is irreducible on the $\GC$-composition factor $U$, where
we choose $U$ to be $L(4\om_{1}) \simeq V \otimes \Vcc$ if $\ell = 3$, and $\WD(V)$ if
$\ell > 3$, cf. Proposition \ref{v4}.

Recall that $G = ZE \cdot S$. As $E$ is the subgroup of elements of order $\leq 4$ in $ZE$,
$E$ is normal in $G$, and $G$ has only one orbit on nontrivial linear characters of $E/Z(E)$.
Also, $Z(E)$ acts trivially on $U$, but $E$ does not. It follows that $2^{2a}-1$ divides
$\dim(U)$. If $\ell = 3$, this gives a contradiction as $\dim(U) = 2^{2a}$. Assume $\ell > 3$.
Then  $\dim(U) = 2^{a}(2^{a}-1)(2^{a}-2)(2^{a}-3)/24$. If $a > 3$, there exists a primitive
prime divisor $\pi$ of $2^{2a}-1$, and $\pi$ does not divide $\dim(U)$, again a contradiction.
If $a=3$, then we compute directly that $2^{2a}-1$ does not divide $\dim(U)$.
\end{proof}

\begin{lemma}\label{ex-p3}
{\sl Assume $p=3$ and $a > 1$. Then $G$ is not irreducible on $L(3\om_{1})$.}
\end{lemma}

\begin{proof}
Note that $Z(E)$ is trivial on $V^{\otimes 3}$ but that $E$ is not trivial on $\SC(V)$. The
dimension of any irreducible $G/Z(E)$-module with $E$ acting nontrivially is a multiple of
$d^{2}-1$ (since $G$ acts transitively on the nontrivial linear characters of $E/Z(E)$, cf.
Lemma \ref{ex-irr}).
Assume $\ell = 2$. Then $L(3\om_{1}) \simeq V \otimes \Vbb$ has dimension $d^{2}$, whence
$G$ is reducible on $L(3\om_{1})$.

Assume $\ell > 3$. Then $L(3\om_{1}) = \SC(V)$ by Proposition \ref{v4}. Let $\pi$ be a primitive
prime divisor of $d-1$, cf. \cite{Zs}. Then $\pi > 3$ and $d \equiv 1 (\mod \pi)$, hence
$\pi$ does not divide the dimension of $\SC(V)$ and so $\SC(V)$ is not $G$-irreducible.
\end{proof}

\begin{lemma}\label{ex-p4}
{\sl Assume $\ell > 3$, and $p^{a} \geq 5$. Then $G$ fails $\MDA$.}
\end{lemma}

\begin{proof}
Assume the contrary: $\ell > 3$ but $G$ satisfies $\MDA$. Let $H = E:S$ with $S := Sp_{2a}(p)$
(notice $S$ embeds in $G$ as $p > 2$).
Since $G = HZ$, $\MDA$ holds for $H$ as well. Observe that the $H$-module $V$ lifts to
a complex module $\VC$. Remark \ref{v5} implies that $H$ satisfies
the condition $M_{8}(\GCC)$ for $\GCC := GL(\VC)$. But this is impossible by Proposition \ref{ex-m8}.
\end{proof}

Finally, we show that $G$ fails $\MCA$ or $\MDA$ if $\ell = 2,3$ and $p^{a} > 4$.
Assume the contrary. By Lemmas \ref{ex-p2} and \ref{ex-p3}, we may assume that $p \geq 5$.
Clearly, $G$ is irreducible on a $\GC$-composition factor $U$ of $V^{\otimes (\ell+1)}$. We
choose $U$ to be $L((\ell+1)\om_{1}) \simeq V \otimes V^{(\ell)}$ if $a > 1$ or if
$(a,\ell) = (1,2)$, and $L(2\om_{1}+\om_{2})$ if $(a,\ell) = (1,3)$.
As in the proof of Lemma \ref{ex-p4}, we may replace $G$ by $E:S$. In the former case, since
$p > \ell+1$, $U \da_{E}$ is a direct sum of $d = p^{a}$ copies of a faithful irreducible
representation say $T$ of $E$. Clifford's theory now implies that $U \simeq T' \otimes R$,
where $T'$ is a projective $G$-representation with $T' \da_{E} = T$, and $R$ is a projective
$S$-representation. But $\Mult(S) = 1$, whence we may assume that both $T'$ and $R$ are linear
representations. Since $G$ is irreducible on $U$, $R$ is irreducible over $S$. However, if
$a > 1$ then $S$ has no irreducible representation of degree $\dim(R) = p^{a}$, cf. \cite{GMST}.
Assume $a = 1$. If $\ell = 2$, then $\dim(R) = p$ by our choice of $U$, and again
$S = SL_{2}(p)$ has no irreducible $2$-modular representation of degree $p$. Assume
$\ell = 3$. Then $U = L(2\om_{1}+\om_{2})$ by our choice, and so $\dim(U)$ is
$p(p+2)(p^{2}-1)/8$ if $p \leq 11$ and at least $(p-1)^{3}/8$ if $p \geq 13$, cf. \cite{Lu2}.
In particular, $\dim(U) > p(p+1)$, $\dim(R) > p+1$ and so $R$ cannot be irreducible over
$S = SL_{2}(p)$, a contradiction.

\subsection{Case $\GC = Sp(V)$ or $O(V)$.}
Now we deal with the case $\GC = Sp(V)$ or $O(V)$. In particular $p = 2$, $a \geq 3$, and
$G = E \cdot S$ (and the extension is nonsplit). We claim that $G$ fails $\MDA$. Assume the
contrary. Then $G$ is irreducible on the $\GC$-composition factor $U$ of $\VD$, where we
choose $U$ to be $L(4\om_{1}) \simeq V \otimes \Vcc$ if $\ell = 3$, $\SD(V)$ if $\ell > 3$
and $\GC = Sp(V)$, $\WD(V)$ if $\ell > 3$ and $\GC = O(V)$ but $\dim(V) > 8$, and
$L(4\om_{1})$ if $\ell > 3$, $\GC = O(V)$ and $\dim(V) = 8$. This choice is possible because of
Proposition \ref{v4}(ii), (iii). In all cases, $Z(E)$ acts trivially on $U$ but $E$ does not.

First we consider the case $\GC = Sp(V)$. By Lemma \ref{ex-irr}, $G$ has two orbits on
$\Irr(E/Z(E)) \setminus \{1_{E}\}$, of length $n_{1} := d(d+1)/2$ and $n_{2} := (d-2)(d+1)/2$,
where $d = \dim(V) = 2^{a}$ as usual. By Clifford's theorem, $\dim(U)$ is divisible by
$n_{1}$ or $n_{2}$. However, this is impossible if $\ell = 3$ as $\dim(U) = d^{2}$. Assume
$\ell > 3$. Then $\dim(U) = d(d+1)(d+2)(d+3)/24$, and again $\dim(U)$ is not divisible by $n_{1}$
nor by $n_{2}$.

Next let $\GC = O(V)$. By Lemma \ref{ex-irr}, $G$ has two orbits on
$\Irr(E/Z(E)) \setminus \{1_{E}\}$, of length $n_{1} := d(d-1)/2$ and $n_{2} := (d+2)(d-1)/2$.
By Clifford's theorem, $\dim(U)$ is divisible by $n_{1}$ or $n_{2}$. However, this is impossible
if $\ell = 3$ as $\dim(U) = d^{2}$. Assume $\ell > 3$ and $d > 8$. Then
$\dim(U) = d(d-1)(d-2)(d-3)/24$, and again $\dim(U)$ is not divisible by $n_{1}$
and by $n_{2}$. Finally, assume $\ell > 3$ and $d = 8$. Then $\dim(U) = 293$ or $294$ (cf.
\cite{Lu2}), and $n_{1} = 28$, $n_{2} = 35$, again a contradidiction.

We have proved

\begin{theor}\label{ext}
{\sl Assume $G < \GC$ satisfies $\MBA$ and $F^{*}(G) = Z(G)E$ is a group of symplectic type with
$|E|= p^{1+2a}$. Assume in addition that $p^{a} \geq 5$. Then $G$ satisfies one of the
conclusions {\rm (i) -- (iii)} of Lemma \ref{ex-irr}. Moreover, $G$ cannot satisfy
$\MCA \cap \MDA$.
\hfill $\Box$}
\end{theor}

We finish this section with the following statement

\begin{propo}\label{ex-m6}
{\sl {\rm (i)} Let $G = (\ZZ_{4} * 2^{1+2n}) \cdot Sp_{2n}(2)$ be embedded in $\GC = GL(V)$ with
$V = \FF^{2^{n}}$ and $n > 3$. Assume in addition that $(\ell,2^{2n}-1) = 1$. Then $G$
satisfies $\MB \cap \MC$.

{\rm (ii)} Let $H = 2^{1+2n}_{\eps} \cdot O^{\eps}_{2n}(2)$ be embedded in $\HC = O(U)$ for
$\eps = +$ and in $\HC = Sp(U)$ for $\eps = -$, where $U = \CC^{2^{n}}$. Assume $n \geq 3$.
Then $H$ satisfies $M_{6}(\HC)$.}
\end{propo}

\begin{proof}
(i) Notice that $V$ lifts to a complex module $U$. Our strategy is to decompose the modules
$U^{\otimes 3}$ and $U^{\otimes 2} \otimes U^{*}$ over $G$, and then reduce modulo $\ell$.
Let $d = 2^{n}$, $E := O_{2}(G)$ and $\GCC = GL(U)$ as usual. By Lemma \ref{ex-irr}, $G$ satisfies
$\MB$ and so it is irreducible on $\AC(U)$, $\SB(U)$, and $\WB(U)$.

1) Let $X$ be any irreducible constituent of the $G$-module $U^{\otimes 3}$. Inspecting the
action of $Z(E)$ on $X$, we see that $X \da_{E}$ is a direct sum of some copies of
$U^{*} \da_{E}$. Clearly, $U^{*} \da_{E}$ is irreducible and extends to $G$. Hence by Clifford's
theory $X \simeq U^{*} \otimes T$ for some $S$-module $T$, where $S := G/E$.
Since we are considering complex modules, the same conclusion holds for any $G$-submodule of
$U^{\otimes 3}$. Similarly, any $G$-submodule of $U^{\otimes 2} \otimes U^{*}$ can be written
as $U \otimes T$ for some $S$-submodule $T$.

2) Over $\GCC$ we have $U^{\otimes 3} = \SC(U) \oplus \WC(U) \oplus 2L$, where
$L := L(\om_{1}+\om_{2})$. Next,
$(U^{\otimes 2} \otimes U^{*}, U) = (U \otimes U^{*}, U \otimes U^{*}) = 2$, both over $G$ and
$\GCC$. Also, $U^{\otimes 2} \otimes U^{*}$ has $\GC$-submodules $M := L(2\om_{1} + \om_{d-1})$,
$N := L(\om_{2} + \om_{d-1})$. Hence over $\GC$ we have
$U^{\otimes 2} \otimes U^{*} = M \oplus N \oplus 2U$.

Applying 1) to the summands of these two decompositions,
we get $\SC(U) \da_{G} = U^{*} \otimes T_{1}$, $\WC(U) \da_{G} = U^{*} \otimes T_{2}$,
$L \da_{G} = U^{*} \otimes T_{3}$, $M \da_{G} = U \otimes T_{4}$, and
$N \da_{G} = U \otimes T_{5}$ for some $S$-modules $T_{i}$. Observe that
\begin{equation}\label{ex5}
  (U^{\otimes 3},U^{*})_{G} = (U \otimes U,U^{*} \otimes U^{*})_{G} =
    (\SB(U) \oplus \WB(U),\SB(U^{*}) \oplus \WB(U^{*}))_{G} \leq 2.
\end{equation}
It is easy to see that $\dim(T_{i}) \leq 2^{2n}-2^{n+1}$ for all $i$. According to \cite{GT1},
$S$ has exactly $6$ irreducible complex modules of dimension $\leq 2^{2n}-2^{n+1}$, namely,
the trivial one, the so-called {\it unitary Weil} modules $\an$, $\bn$, $\zn$, and {\it linear
Weil} modules $\ran$, $\rbn$. Here,
$\dim(\bn) = \dim(T_{1}) = (2^{n}+1)(2^{n}+2)/6$,
$\dim(\an) = \dim(T_{2}) = (2^{n}-1)(2^{n}-2)/6$,
$\dim(\zn) = \dim(T_{3}) = (2^{2n}-1)/3$, $\dim(\ran) = (2^{n}+1)(2^{n-1}-1)$, and
$\dim(\rbn) = (2^{n}-1)(2^{n-1}+1)$. In particular,
\begin{equation}\label{ex6}
  \dim(\ran) + \dim(\rbn) = \dim(T_{4}) + \dim(T_{5}).
\end{equation}
Keeping (\ref{ex5}) in mind, we can easily deduce that $T_{2} = \an$, $T_{1} = \bn$ and
$T_{3} = \zn$. Thus we have shown that $G$ satisfies the condition $M_{6}(\GCC)$, i.e.
$M_{6}(G,U) = 6$. This in turn implies that $T_{4}$ and $T_{5}$ are irreducible. Also notice that
$$(U^{\otimes 3}, U^{\otimes 2} \otimes U^{*})_{G} = (U^{\otimes 4}, U^{\otimes 2})_{G} = 0$$
since $Z(E)$ acts trivially on $U^{\otimes 4}$ but not on $U^{\otimes 2}$. In particular,
$T_{4}$ and $T_{5}$ are not isomorphic to $\an$, $\bn$, and $\zn$. This observation and
(\ref{ex6}) imply that $\{T_{4},T_{5}\} = \{\ran,\rbn\}$.

3) Finally, assume $(\ell,2^{2n}-1) = 1$. Then all $T_{i}$'s are irreducible modulo $\ell$,
according to \cite[Cor. 7.5, 7.10]{GT1}. It follows that $G$ is irreducible on $L$, $M$, and
$N$, whence $G$ satisfies $\MC$.

(ii) This has been proved in Propositions \ref{ex-m8} and \ref{ex-m62}.
\end{proof}

Proposition \ref{ex-m6}(ii) has also been proved in \cite{NRS} by different means.

\section{Nearly Simple Groups. I}\label{nearly1}

In this and the two subsequent sections we consider the finite subgroups $G < \GC$ that
satisfy the condition $\MBA$ and have a unique component $S$. Then $S$ is quasi-simple and
$\BS := S/Z(S)$ is simple. Let $\chi$ denote the Brauer character of $G$ on $V$ and let
$M := Z(G)S$. Assume $d := \dim(V) \geq 3$ if $\GC = GL(V)$ and $d \geq 5$ if $\GC = Sp(V)$
or $O(V)$.

We begin with a simple observation which reduces some computations on $G$ to those on $M$:

\begin{lemma}\label{g-m}
{\sl Under the above assumptions,

{\rm (i)} $(G:M) \leq |\Out(S)| \leq |\Out(\BS)|$;

{\rm (ii)} $(G:C_{G}(Z)) \leq |\Out(\BS)| \cdot (S:C_{S}(Z))$ and
$(G:N_{G}(Z)) \leq |\Out(\BS)| \cdot (S:N_{S}(Z))$ for any $Z \leq S$.

{\rm (iii)} For any nontrivial composition factor $X$ of the $\GC$-module
$V \otimes V^{*}$, $X \da S$ is a direct sum of at most $(G:M)$ irreducible constituents
of same dimension. In particular, if $G = M$ then $S$ also satisfies $\MBA$.}
\end{lemma}

\begin{proof}
(i) Clearly $S \lhd G$, so we get a homomorphism $\pi~:~G \to \Aut(S)$ with
$\Ker(\pi) = C_{G}(S)$. By Proposition \ref{irred1}, $S$ is irreducible on $V$, whence
$C_{G}(S)$ acts scalarly on $V$ by Schur's lemma. Since $G$ is faithful on $V$, it follows
that $C_{G}(S) = Z(G)$ and $\pi$ embeds $G/Z(G)$ in $\Aut(S)$. Also, $Z(G) \cap S = Z(S)$,
so $\pi$ maps $Z(G)S/Z(G) \simeq S/Z(S)$ isomorphically onto $\Inn(S)$. Thus
$$G/M \simeq (G/Z(G))/(Z(G)S/Z(G)) \leq \Aut(S)/\Inn(S) \leq \Out(S).$$
A standard argument shows $\Out(S)$ embeds in $\Out(\BS)$.

(ii) Since $C_{M}(Z) = Z(G)C_{S}(Z)$ and $Z(G) \cap S = Z(G) \cap C_{S}(Z) = Z(S)$, we have
$$\frac{|G|}{|C_{G}(Z)|} \leq \frac{|G|}{|C_{M}(Z)|} =
  \frac{(G:M) \cdot |Z(G)S|}{|Z(G)C_{S}(Z)|} =
  \frac{(G:M) \cdot |S|}{|C_{S}(Z)|} \leq \frac{|\Out(\BS)| \cdot |S|}{|C_{S}(Z)|}.$$
Similarly for $N_{G}(Z)$.

(iii) is clear, since $M = Z(G)S$ and $Z(G)$ acts scalarly on $X$.
\end{proof}

The next observations are useful in dealing with nearly simple groups:

\begin{lemma}\label{change1}
{\sl Let $H < GL(V)$ such that $\FF^{*} \cdot G = \FF^{*} \cdot H$, where
$\FF^{*} = Z(GL(V))$. Let $W$ be any $\GC$-composition factor of
$V^{\otimes k} \otimes (V^{*})^{\otimes l}$. Then $G$ and $H$ both act on $W$ (although
$H$ may not necessarily be contained in $\GC$). Moreover, if $G$ is irreducible on $W$ then
so is $H$.}
\end{lemma}

\begin{proof}
$G$ acts on $W$ as $G < \GC$. Next, $\FF^{*}$ acts (scalarly) on $V$ and on $V^{*}$, so it
also acts scalarly on $V^{\otimes k} \otimes (V^{*})^{\otimes l}$. In particular, $\FF^{*}$
acts scalarly on $W$. This implies that $H$ acts on $W$, and moreover $H$ is irreducible
on $W$ if $G$ is irreducible on $W$.
\end{proof}

\begin{lemma}\label{change2}
{\sl Assume $V = \FF^{d}$ with $d \geq 5$ and $G \leq \GC \leq GL(V)$. Assume there is a
subgroup $H \leq GL(V)$ such that $\FF^{*} \cdot G = \FF^{*} \cdot H$ (where
$\FF^{*} = Z(GL(V))$ and $V\da_{H}$ is self-dual. Then one of the following holds:

{\rm (i)} $V\da_{G}$ is self-dual and $\GC \neq GL(V)$;

{\rm (ii)} Both $\MBA$ and $\MC$ fail for $G$.}
\end{lemma}

\begin{proof}
By way of contradiction, assume (ii) fails but $\GC = GL(V)$.

1) First suppose that $\MBA$ holds for $G$. The $H$-module $V \otimes V^{*}$ has
subquotients $\SB(V)$, $\WB(V)$ if $\ell \neq 2$, and $\WB(V)$ (twice), $\Vbb$ if
$\ell = 2$. It follows that $H$ is reducible on $\AC(V)$, whence $G$ is reducible on
$\AC(V)$ by Lemma \ref{change1}, a contradiction.

2) Next we suppose that $\MC$ holds for $G$.

By Lemma \ref{weight}, $\WC(V)$ is $\GC$-irreducible. If $\ell \neq 2$ and $V\da_{H}$ is
of type $-$, or if $\ell = 2$, then the contraction map projects the $H$-module
$\WC(V)$ onto $V$, hence $\WC(V)$ is reducible over $H$ and so over $G$ as well, contrary
to the condition $\MC$ for $G$.

Assume $\ell \neq 2$ and $V\da_{H}$ is of type $+$. Then $H \leq \HC := O(V)$. Let
$W$ be the irreducible $\GC$-module of highest weight $\om_{1}+\om_{2}$ inside $V^{\otimes 3}$.
By \cite[(8.1)]{S2}, $W$ is reducible over $\HC$ (if $d = 6$ we consider $\HC$ as having
type $A_{3}$). Therefore $W$ is reducible over $H$ and over $G$, again a contradiction.
\end{proof}

Notice that, in the modular case $M_{2k}(\GC)$ does not imply $M_{2k-2}(\GC)$; in particular
$\MB$ may fail for some groups satisfying $\MC$.

\medskip
Next we show that the situation where the $G$-module $V$ is not self-dual but $V \da_{S}$ is
self-dual can rarely happen.

\begin{lemma}\label{rare}
{\sl Assume $G < \GC$ satisfies $\MBA$, $d \geq 2$, and $S = E(G)$ is quasi-simple. Assume
$V\da_{S}$ is self-dual, but either the $G$-module $V$ is not self-dual or $\GC = GL(V)$. Then

{\rm (i)} Any $S$-composition factor of $\AC(V)$ has dimension $\leq d+1$;

{\rm (ii)} $|\Out(\BS)| \geq d-1$.}
\end{lemma}

\begin{proof}
(i) By the assumptions $\GC = GL(V)$, and $\MBA$ implies that $G$ is
irreducible on $\AC(V)$. Let $e$ be the common degree of $S$-composition factors of $\AC(V)$.
Observe that $(V \otimes V^{*}) \da_{S} \simeq V^{\otimes 2} \da_{S}$.

Assume $\ell = 2$. Then $\Vbb$ is a subquotient of $V^{\otimes 2}$ and it is irreducible
over $S$ as $S$ is irreducible on $V$ by Proposition \ref{irred1}. But $d \geq 2$, so $\Vbb$ is an
$S$-composition factor of $\AC(V)$. Thus $e = d$.

Assume $\ell \neq 2$. Then $V^{\otimes 2} = \SB(V) \oplus \WB(V)$. Counting the dimensions of
$S$-composition factors in $\SB(V)$, $\WB(V)$, and $\AC(V)$, we see that
$d(d+1)/2 = me+x$ and $d(d-1)/2 = ne + y$ for some nonnegative integers $m,n,x,y$ and
$x + y \leq 2$. In particular, $d+y-x  = (m-n)e$. As $d-1 \leq d+y-x \leq d+1$, we get
$m > n$ and $e \leq d+1$.

(ii) Clearly, $M$ fixes every $S$-composition factor of $\AC(V)$. Hence by Clifford's theorem
and Lemma \ref{g-m}(i), $|\Out(\BS)| \geq (G:M) \geq \dim(\AC(V))/e \geq (d^{2}-2)/(d+1) > d-2$,
as stated.
\end{proof}

Most of the times, we will apply Propositions \ref{bound1} and \ref{bound2} to $C = C_{G}(Z)$
and $N = N_{G}(Z)$ for a suitable cyclic $\ell'$-subgroup $Z < G$ (so $N/C$ is abelian). In
this setup, $V\da_{C} = \oplus^{t}_{i=1}V'_{i}$, where $V'_{i}$ is the $\lam_{i}$-eigenspace
for $Z$ and the spectrum $\Spec(Z,V)$ of $Z$ on $V$ contains exactly $t$ pairwise
distinct characters $\lam_{1}, \ldots ,\lam_{t}$. Once we are able to apply Proposition
\ref{bound1}, we will have $X(V) \hra \lam \ua G$, and so the following key inequality holds:
\setcounter{equation}{0}
\begin{equation}\label{bound3}
  d \leq \frac{1}{2} + \sqrt{2(G:N) + \frac{17}{4}}\;.
\end{equation}

In the rest of this section we treat the case $\BS = \AAA_{n}$ with $n \geq 8$. The cases
$8 \leq n \leq 13$ can be checked directly using \cite{Atlas} and \cite{JLPW}, so
we will assume $n \geq 14$ in the following analysis.

\medskip
{\bf Case 1: $\ell \neq 3$.} Choose $Z = \langle c \rangle \simeq \ZZ_{3}$, where $c$ is
an inverse image of order $3$ of a $3$-cycle in $S$. Since $Z \not\leq Z(G)$, 
$|\Spec(Z,V)| \geq 2$. If $|\Spec(Z,V)| = 3$, then we label the $\lam_{i}$'s such that
$\lam_{1} = 1_{Z}$ and set $V_{1} = V'_{1}$, $V_{2} = V'_{2} \oplus V'_{3}$. This ensures that 
$\Hom_{C}(V_{1},V_{2}) = 0$ (as it is so over $Z$). Assume $|\Spec(Z,V)| = 2$. Then, according 
to \cite{Wa1}, $\ell \neq 2$, $S = \HA_{n}$ is the double cover of $\AAA_{n}$, $V \da S$ is a 
basic spin module, and $\lam_{2} = \bar{\lam}_{1} \neq 1_{Z}$. Setting $V_{i} = V'_{i}$, we see that
$V_{2} \simeq V_{1}^{*}$ as $C$-modules, provided that $V = V^{*}$ as $\GC$-modules. Thus in
all cases the assumptions of either (i), or (ii), or (iii) of Proposition \ref{bound1} are
fulfilled. Hence there is a $1$-dimensional $\FF N$-module $\lam$ and $X \in \{\AC,\TSB, \TWB\}$
such that $X(V) \hra \lam \ua G$. By Lemma \ref{g-m}, $(G:N) \leq n(n-1)(n-2)/3$. Hence
(\ref{bound3}) implies
\begin{equation}\label{an-1}
  \dim(V) < (n-2)(n-3)/2
\end{equation}
if $n \geq 9$, and
\begin{equation}\label{an-2}
  \dim(V) < n(n-3)/4 < 2^{[(n-3)/2]}
\end{equation}
if $n \geq 15$.

Assume $n \geq 15$. According to \cite{KT}, the dimension of any faithful $\HA_{n}$-module
is at least $2^{[(n-3)/2]}$. Hence (\ref{an-2}) implies that $S = \AAA_{n}$. Applying
(\ref{an-2}) and \cite[Thm. 7]{J} we see that $V \da S$ is a composition factor of the
irreducible $\FF \SSS_{n}$-module $D^{\mu}$ labelled by an $\ell$-regular partition
$\mu \in \{(n-1,1), (n-2,1^{2}), (n-2,2)\}$. Since $\dim(V) \geq (n^{2}-5n+2)/2$ for
$\mu = (n-2,2)$ or $(n-2,1^{2})$, we get $\mu = (n-1,1)$. If $(\ell,n) = 1$ and
$\ell \neq 2$, then $\TSB(V) \da S$ contains irreducible constituents
$D^{(n-1,1)}$ and $D^{(n-2,2)}$ by \cite[Lemma 2.1]{MM} which are not $G$-conjugate,
contrary to Lemma \ref{g-m}(iii). So $\ell|n$ or $\ell = 2$. This case
can indeed happen, see \cite[Table 2.1]{MM}.

Assume $n = 14$. Then (\ref{bound3}) implies $d \leq 38$. Using \cite{HM} one can show that either
$S = \AAA_{14}$ and $V\da_{S} = D^{(13,1)}$, or $S = \HA_{14}$, $\ell = 7$, and $d = 32$. Suppose
we are in the latter case. Then $G = Z(G)S$ by \cite{KT}, so we may replace $G$ by $S$ by Lemma
\ref{g-m}. Consider an inverse image $g$ of order $3$ in $S$ of a $3$-cycle. Then
$C := C_{S}(g) = \HA_{11} \times \langle g \rangle$. By Proposition \ref{bound1},
there is an $X \in \{\AC,\TWB, \TSB\}$ such that $X(V) \hra 1_{C} \ua S$. Observe that
$1_{C} \ua S$ is trivial on $Z(S)$, and as an $S/Z(S)$-module it equals
$1_{\AAA_{11} \times \AAA_{3}} \ua \AAA_{14} = W \da_{\AAA_{14}}$, where
$W := 1_{\SSS_{11} \times \AAA_{3}} \ua \SSS_{14}$. Observe that
$W$ has a subquotient $1_{\SSS_{11} \times \SSS_{3}} \ua \SSS_{14}$ of dimension
$364$ which is $\dim(W)/2$. It follows that any composition factor of the
$\SSS_{14}$-module $W$ has dimension $\leq 364$. On the other hand,
$\dim(X(V)) \geq (32 \cdot 31)/2-2 = 494$, a contradiction.

\medskip
{\bf Case 2: $\ell = 3$.} Choose $Z = \langle c \rangle$, where $c$ is an inverse image of
a double transposition in $S$, and set $V_{i} = V'_{i}$. If $S = \AAA_{n}$, then
$c^{2} = 1$ and $c \notin Z(G)$, whence $\Spec(c,V) = \{1,-1\}$ and
$\Hom_{C}(V_{1},V^{*}_{2}) = 0$. If $S = \HA_{n}$, then $c^{2} = -1$ and $c \notin Z(G)$,
whence $\Spec(c,V) = \{i,-i\}$, and $V_{1} \simeq V^{*}_{2}$ if $V = V^{*}$ as $\GC$-modules.
Thus in all cases the assumptions of either (i), or (ii), or (iii) of Proposition \ref{bound1} are
fulfilled. Hence there is a $1$-dimensional $\FF N$-module $\lam$ and $X \in \{\AC,\TSB, \TWB\}$
such that $X(V) \hra \lam \ua G$ and $X(V) \hra 1_{C} \ua G$. Observe that if $G = M$ then
we may assume $G = S$ by Lemma \ref{g-m}(iii), and $C/Z(G)$ contains $\AAA_{n-4}$, so
by Proposition \ref{bound1} $X(V) \hra 1_{\AAA_{n-4}} \ua \AAA_{n}$, i.e. the conclusion of
\cite[Lemma 2.4]{MM} holds.

Assume $n \geq 15$ and $S = \HA_{n}$. One can check that $(G:N) = n(n-1)(n-2)(n-3)/8$.
Therefore, (\ref{bound3}) implies that $d < 2^{[(n-2-\kappa_{n})/2]}$ unless $n = 15,18$,
where $\kappa_{n} = 1$ if $3|n$ and $0$ otherwise. On the other hand, the dimension of
any faithful $\FF\HA_{n}$-module is at least $2^{[(n-2-\kappa_{n})/2]}$ by \cite{KT}.
Hence $n = 15$ or $n = 18$. If $G = M$, then we may assume $G = S$ by Lemma \ref{g-m}(iii)
and get a contradiction by \cite[Prop. 2.5]{MM}. Assume $G > M$. If $n = 18$, then
(\ref{bound3}) implies $d \leq 136$, meanwhile the dimension of any faithful $\FF G$-module
is at least $2^{[(18-1-\kappa_{18})/2]} = 256$ by \cite{KT}, a contradiction. Assume
$n = 15$. Then (\ref{bound3}) implies that $d \leq 91$. Restricting to the subgroup
$B = \HA_{13}$ inside $A$ and using \cite{JLPW}, we see that all irreducible constituents
of $V \da_{B}$ are basic spin, whence $V \da_{S}$ is basic spin by \cite{KT}, and $d = 64$.
Since $C/Z(G)$ contains $\SSS_{n-4}$, $X(V) \hra 1_{\SSS_{n-4}} \ua \SSS_{n}$ by
Proposition \ref{bound1}. Now we can get a contradiction as in the proof of
\cite[Prop. 2.5]{MM}.

Assume $n = 14$ and $S = \HA_{14}$. Then (\ref{bound3}) implies $d \leq 78$, whence
$d = 64$ and $V \da_{S}$ is the reduction modulo $3$ of a (unique) basic spin complex module by
\cite{KT}. In particular, $V \da_{S}$ is self-dual and rational-valued. This implies by
\cite{JLPW} that $V \da_{\HA_{13}}$ has exactly two composition factors, nonisomorphic and
both of type $-$. The latter in turn implies that $V\da_{S}$ is of type $-$.
By Lemma \ref{rare}, $V$ is also self-dual as a $G$-module, whence $V$ has type $-$ and
$\GC = Sp(V)$. So $\SB(V)$ is irreducible over $\GC$, $Z(G) = Z(S)$, and $G = S$ or $\HS_{14}$.
The proof of Proposition \ref{bound1}(iii) shows that $\SB(V) \hra W := 1_{C} \ua G$. Since
$C > Z(G)$ and $C/Z(G) \geq G/Z(G) \cap \AAA_{10}$, $W$ is trivial on $Z(G)$, and
$W \hra (1_{\SSS_{10}} \ua \SSS_{14}) \da G/Z(G)$ as $G/Z(G)$-modules. The largest (complex)
constituent of $1_{\SSS_{10}} \ua \SSS_{14}$ has degree
$(14 \cdot 12 \cdot 11 \cdot 9)/8 = 2079$ (see \cite[p. 174]{MM}). On the other hand,
$\dim(\SB(V)) = (64 \cdot 65)/2 = 2080$, a contradiction.

Assume $n \geq 14$ and $S = \AAA_{n}$. By Proposition \ref{bound1}, $X(V) \hra 1_{C} \ua G$.
In particular, all $S$-composition factors of $X(V)$ are inside $1_{\AAA_{n-4}} \ua \AAA_{n}$,
(as $(G:C) = (S:C_{S}(Z))$ here), i.e. the conclusion of \cite[Lemma 2.4]{MM} holds. Now we
appeal to \cite[Prop. 2.5]{MM} to get $3|n$ and $V \da_{S}$ is labelled by $(n-1,1)$. Since
$V \da_{S}$ is of type $+$ and $|\Out(S)| = 2 < d-1$, $V$ is of type $+$ by Lemma \ref{rare},
whence $Z(G) \leq \ZZ_{2}$.

We have proved

\begin{theor}\label{an}
{\sl Assume $G \leq \GC$, $G$ satisfies $\MBA$ and $S := E(G) \simeq \AAA_{n}$ or $\HA_{n}$.
Assume in addition that $n \geq 8$. Then one of the following holds.

{\rm (i)} Either $\ell|n$ and $d = n-2$, or $\ell = 2$ and
$d = n-1$. Furthermore, $S = \AAA_{n}$, $Z(G) \leq \ZZ_{2}$, $S \leq G/Z(G) \leq \Aut(S)$,
and $V \da_{S}$ is labelled by $(n-1,1)$.

{\rm (ii)} $\BS = \AAA_{9}$, $\ell \neq 3$, and $d = 8$. Furthermore, $S = \HA_{9}$ if
$\ell \neq 2$ and $S = \AAA_{9}$ if $\ell = 2$.
\hfill $\Box$}
\end{theor}

\begin{propo}\label{an-m6}
{\sl Assume $G \leq \GC$ satisfies $\MBA$ with $S := E(G) \simeq \AAA_{n}$, $n \geq 5$.
Assume $d \geq 4$ and $V \da_{S}$ is labelled by $(n-1,1)$. Assume furthermore that either
$\ell|n$ or $\ell = 2$. Then $G$ cannot satisfy $\MCA$ when $\ell \neq 3$, and $\MDA$ when
$\ell = 3$.}
\end{propo}

\begin{proof}
1) We construct a group $H \in \{\AAA_{n},\SSS_{n}\}$ as follows. If $G = Z(G)S$, then set $H := S$.
Otherwise $G$ induces all automorphisms of $S$ by Lemma \ref{g-m}(i). In particular, there is
$t \in G$ whose conjugation on $S = \AAA_{n}$ is the same as the conjugation by a transposition in
$\SSS_{n}$. Then $t^{2} \in C_{G}(S) = Z(G) \leq \FF^{*}$ (see the proof of Lemma \ref{g-m}).
But $\FF$ is algebraically closed, so there is $\nu \in \FF^{*}$ such that $\nu^{2} = t^{2}$,
whence $\tau^{2} = 1$ for $\tau := t\nu^{-1}$. Setting $H = \langle S,\tau \rangle$ we see that
$H \simeq \SSS_{n}$ and $\FF^{*} \cdot G = \FF^{*} \cdot H$. Also denote 
$H' = [H,H] \simeq \AAA_{n}$.

2) Assume $\ell = 2$. Then $W := L(3\om_{1}) \simeq V \otimes \Vbb$. In this case,
$W\da_{H'} \simeq U \otimes U^{(2)} \simeq U \otimes U$ contains $1_{H'}$, and so $W\da_{H}$ is
reducible, whence $G$ is reducible on $W$.

Similarly, let $\ell = 3$. Then $G$ and so $H$ is irreducible on the $\GC$-composition
factor $W' := L(4\om_{1}) \simeq V \otimes \Vcc$. But
$W'\da_{H'} \simeq U \otimes U^{(3)} \simeq U \otimes U$ contains $1_{H'}$, and so $W'\da_{H}$
is reducible, whence $G$ is reducible on $W'$.

3) So we may now assume $\ell > 3$. The cases $5 \leq n \leq 11$ can be checked directly using
\cite{Atlas} and \cite{JLPW}, so we assume $n \geq 12$.
Assume $G$ satisfies $\MCA$. Then $G$ is irreducible on the $\GC$-composition
factor $W$ of $V^{\otimes 3}$ with $W \da_{\GC'} = L(3\om_{1})$.  By Lemma \ref{change1}, $H$ is
also irreducible on $W$.

We are given that $V \da_{S}$ is labelled by $(n-1,1)$. One can show that there are (at most)
two extensions of $V \da_{S}$ to $H$, the deleted permutation module which we will denote by
$U$, and $U \otimes \sign$, where $\sign$ is the sign representation of $\SSS_{n}$. Since
$\SC(U \otimes \sign) \simeq \SC(U) \otimes \sign$, there is no loss to assume that
$V \da_{H} \simeq U$. 

We may consider $U$ as the deleted permutation module for $K := \SSS_{n}$. Consider the
subgroup $N = \SSS_{3} \times \SSS_{n-3}$ of $K$. It is easy to see that
$U \da_{N} = A \oplus B$, where $A$ is the deleted permutation module for $\SSS_{3}$ (and
$\SSS_{n-3}$ acts trivially on $A$), and $B$ is the deleted permutation module for $\SSS_{n-3}$
(and $\SSS_{3}$ acts trivially on $B$). Notice that
$\SC(A) = A \oplus 1_{\SSS_{3}} \oplus \sign$. Hence
$\SC(U) \da_{N}$ contains the direct summand $1_{N}$. It follows that
$0 \neq \Hom_{N}(\SC(U),1_{N}) \simeq \Hom_{K}(\SC(U),1_{N} \ua K)$. Consider
a nonzero map $f \in \Hom_{K}(\SC(U),1_{N} \ua K)$.

If $\GC = GL(V)$ or $Sp(V)$, then $W = \SC(V)$ by Proposition \ref{v4}. If $\GC = O(V)$
then $\dim(W)$ is bounded below by Proposition \ref{v3}. Thus in all cases we have
$\dim(L(3\om_{1})) > n(n-1)(n-5)/6$, where the right-hand side equals $\dim(S^{n-3,3})$ and so
it is an upper bound on the dimension of irreducible constituents of $1_{N} \ua K$. It follows
that $W \hra \Ker(f)$ as $H$ is irreducible on $W$. But $f \neq 0$, hence $\GC = O(V)$, and all 
composition factors of $\SC(U)/\Ker(f)$ have dimension equal to $\dim(U)$. This implies that either 
$U$ or $U \otimes \sign$ embeds in $1_{N} \ua K$, whence $1_{N}$ enters $U \da_{N}$ or
$(U \otimes \sign)\da_{N}$. This conclusion contradicts the decomposition
$U \da_{N} = A \oplus B$ mentioned above.
\end{proof}

\section{Nearly simple groups. II}\label{nearly2}
In this section we continue to treat the nearly simple groups and assume that
$S$ is a finite quasi-simple group of Lie type,
defined over a field $\FQ$ of characteristic $p \neq \ell$, and that $\BS$ is not isomorphic
to any of the following groups: $L_{2}(q)$ with $q = 5,7,9,11,13$, $L_{3}(4)$,
$SL_{4}(2) \simeq \AAA_{8}$, $U_{3}(q)$ and $U_{4}(q)$ with $q = 2,3$, $U_{6}(2)$,
$Sp_{4}(2)' \simeq \AAA_{6}$, $Sp_{4}(4)$, $Sp_{6}(2)$, $\Om_{7}(3)$, $\Omp_{8}(2)$,
$\ta B_{2}(8)$, $G_{2}(q)$ with $q = 3,4$, $F_{4}(2)$, $\ta F_{4}(2)'$, $\tb D_{4}(2)$, and
$\ta E_{6}(2)$ (they are included in $\cup^{3}_{i=1}\EC_{i}$).
Hence $S$ is a quotient of $\HS$, the corresponding finite Lie-type group of
simply connected type. We will write $q = p^{f}$ and denote by $\dd(\BS)$ the smallest degree
of faithful projective $\FF \BS$-representations. Lower bounds for $\dd(\BS)$ were given in
\cite{LS,SZ}.

\subsection{Some generalities.}

We will apply Proposition \ref{bound1} to $C := C_{G}(Z)$ and $N := N_{G}(Z)$, where
$Z$ is a long-root subgroup. This choice is justified by the following statement:

\begin{propo}{\rm \cite[Cor. 2.10]{MMT}}\label{longroot}
{\sl Let $G$ be a universal-type quasi-simple finite group of Lie type defined over $\FF_{q}$, 
$q = p^{f}$, and $Z$ a long-root subgroup of $G$. Let $\FF$ be an algebraically closed field 
of characteristic $\ell \neq p$ and $V$ a nontrivial irreducible $\FF G$-module. Then either

{\rm (i)} $\Spec(Z,V)$ contains at least two distinct characters which are not dual to each 
other, or

{\rm (ii)} $G \in \{SL_{2}(5),~SU_{3}(3),~Sp_{4}(3)\}$.
\hfill $\Box$}
\end{propo}

Proposition \ref{longroot} implies that $\Spec(Z,V)$ contains at least two distinct characters, 
say $\lam_{1}$ and $\lam_{2}$, which are not dual to each other. 
Write $V|_{Z} = V_{1} \oplus V_{2}$, where
$V_{1}$ can afford only $Z$-characters $\lam_{1}$ and $\overline{\lam}_{1}$, and $V_{2}$ affords
only the remaining $Z$-characters. Then clearly $V_{1}$ and $V_{2}$ are $C$-stable, and they
satisfy the assumptions of (i) and (ii) of Proposition \ref{bound1}. Thus there is
$X \in \{\AC,\TSB,\TWB\}$ and an $1$-dimensional $\FF N$-module $\lam$ such that
$X(V) \hra \lam \ua G$; in particular, (\ref{bound3}) holds.

Next we record the following refinement to Lemma \ref{g-m}(ii):

\begin{lemma}\label{root}
{\sl Under the above assumptions, assume that $S$ is not of types $B_{2}$ or $F_{4}$
if $p = 2$, and type $G_{2}$ if $p = 3$. If $Z < S$ is a long-root subgroup, then
$(G:N_{G}(Z)) = (S:N_{S}(Z))$.}
\end{lemma}

\begin{proof}
Let $\Om$ be the set of long-root subgroups in $S$. Since $S$ is transitive on $\Om$,
$|\Om| = (S:N_{S}(Z))$. Hence it suffices to show that $Z^{g} \in \Om$ for any $g \in G$.
Since $S \lhd G$, it suffices to show that $\sigma(Z) \in \Om$ for any $\sigma \in \Aut(S)$.
Any such a $\sigma$ is a product of inner, diagonal, graph, and field automorphisms.

The case $S$ is of type $\ta A_{2}$ may be checked directly. So we may assume that $S$ is not
of type $\ta A_{2}$, and let $Z = \{x_{\al}(t) \mid t \in \FQ\}$ for some long root $\al$. Also, if
$S$ is of type $\ta B_{2}$, $\ta G_{2}$, or $\ta F_{4}$, then $\Out(S)$ is generated by a
field automorphism.

The claim is obvious for inner automorphisms. If $\sigma$ is diagonal, then we may assume
that $\sigma~:~x_{\al}(t) \mapsto x_{\al}(f(t))$ for some $f~:~\FQ \to \FQ$ and so
$\sigma(Z) = Z$. If $\sigma$ is a field automorphism, then we may assume
$\sigma~:~x_{\al}(t) \mapsto x_{\al}(t^{\gam})$ for some $\gam \in \Aut(\FQ)$, so
$\sigma(Z) = Z$. Finally, assume that $\sigma$ is a nontrivial graph automorphism. Then $S$ is of
type $A_{n}$, $D_{n}$, or $E_{6}$. For the last two types, we may choose $\al$ to be fixed
by $\sigma$, whence $\sigma(Z) = Z$. If $S$ is of type $A_{n}$, then we may assume that
$\sigma$ is induced by the map $X \mapsto \tn X^{-1}$ on $SL_{n+1}(q)$, and so it is easy to
check that $\sigma(Z) = Z$ for some $Z \in \Om$.
\end{proof}

\subsection{Non-generic case.}

In this subsection we consider the non-generic case, that is, $S$ is of types $B_{2}$,
$F_{4}$ if $p = 2$, or type $G_{2}$ if $p = 3$.

Assume $S = Sp_{4}(q)$ with $q = 2^{f}$. Then
$(G:N) \leq |\Out(\BS)| \cdot (S:N_{S}(Z)) = 2f(q^{4}-1)/(q-1)$ by Lemma \ref{g-m}. If
$q \geq 8$, then (\ref{bound3}) implies $d < q(q-1)^{2} = \dd(\BS)$, a contradiction. Assume
$q = 4$. Then (\ref{bound3}) implies $d \leq 26$. According to \cite{Atlas, JLPW},
$\dim(V) = 18$ (as we assume $\ell \neq 2$). This case is completed in the following
Proposition, which shows in particular that $\MBA$ may hold for nearly simple groups but not
for their generalized Fitting subgroups:

\begin{propo}\label{sp44}
{\sl Assume $G \leq \GC$ satisfies $\MBA$ with $S := E(G) \simeq Sp_{4}(4)$. Then $G$ is the
unique group determined by

{\rm (i)} $G$ is a non-split extension  of $\ZZ_{2} \times S$ by $\ZZ_{4}$, and

{\rm (ii)} $G/Z(G) \simeq \Aut(S)$.\\
Furthermore, $\ell \neq 2,3$, and $V$ is any of two self-dual irreducible representations of
$G$ of degree $18$. Conversely, if $G,V,\ell$ are as described, then $G < \GC = O(V)$ satisfies
$\MBA$, but it fails both $\MCA$ and $\MDA$. Finally, if $E(G) = Sp_{4}(4)$ then
$G$ cannot satisfy $\MC$ nor $\MDA$.}
\end{propo}

\begin{proof}
1) Assume $E(G) = S = Sp_{4}(4)$ and $G$ satisfies $\MBA$. We have shown above that $d = 18$.
In particular, $V$ lifts to a complex module $\VC$, and $V \da_{S}$ is uniquely determined.
Inspecting \cite{Atlas, JLPW}, we see that $\WB(V)$ is irreducible over $S$. On the other hand,
if $\ell = 3$ then $\TSB(V)\da_{S}$ contains irreducible constituents of degree $34$ and $50$,
contrary to Lemma \ref{g-m}(iii). So $\ell \neq 3$. In this case, $\TSB(V)\da_{S}$ has
exactly two irreducible constituents say $A$ and $B$, both of degree $85$, which extend to
$S \cdot \ZZ_{2}$ but fuse under $\Aut(S)$. Hence $G/Z(G) \simeq \Aut(S)$. Since
$V\da_{S}$ is of type $+$ and $|\Out(S)| = 4 < d-1$, $V$ is of type $+$ by Lemma \ref{rare}.
In particular, $|Z(G)| \leq 2$. Observe that $\Aut(S)$ cannot be embedded in $G$, since
$\Aut(S)$ has no self-dual irreducible representations of degree $18$. It follows that
$Z(G) = \ZZ_{2}$ and $G$ is a non-split extension  of $\ZZ_{2} \times S$ by $\ZZ_{4}$.

Conversely, suppose $G$ satisfies (i) and (ii). Then we can construct $G$ as follows.
Let $H = \Aut(S) = \langle H,x\rangle$ where $x$ is an element of class $4F$ in the notation of
\cite{Atlas}. Embed $H$ in $GL(\VC)$ with $\VC = \CC^{18}$, and set $G = \langle H,y\rangle$
with $y = \exp(\pi i/4)x$. If $\chi$ is the character of $G$ on $\VC$, then direct computation
shows that $\chi = \overline{\chi}$. Since $\chi\da_{S}$ is of type $+$, $\chi$ itself is of
type $+$. As we mentioned above, $\WB(\chi)$ is irreducible over $S$. Furthermore,
$\TSB(\chi)$ is irreducible, since $G$ fuses the two irreducible constituents of
$\TSB(\chi)\da_{S}$. Thus $G$ satisfies $\MBA$ for $\GC = O(\VC)$. The same picture holds when
we reduce $\VC$ modulo $\ell \neq 2,3$. Also, the same happens if we replace $\chi$ by the
other (faithful) real extension of $\chi\da_{S}$ to $G$. Furthermore, $\WD(\chi)$ has degree
$3060$, so $G$ cannot be irreducible on $\WD(V)$, whence $\MDA$ fails for $G$ by Proposition
\ref{v4}. Finally, the composition factor $L(3\om_{1})$ of the $\GC$-module $\SC(V)$ has degree
$1104$ if $\ell = 5$ and $1122$ if $\ell > 5$, cf. \cite{Lu2}. It follows that $G$ cannot be
irreducible on this subquotient, whence $\MCA$ fails for $G$.

2) Assume $E(G) = Sp_{4}(4)$ and $G$ satisfies $\MC$ or $\MDA$. From \cite{Atlas} we
see that any irreducible $\FF G$-module has dimension $\leq 1020$. If $d := \dim(V) \geq 22$,
then $V^{\otimes 3}$ contains $\GC$-composition factors of degree at least $1320$ by Lemma
\ref{weight}, whence $\MC$ fails. Similarly, if $d \geq 18$ then $V^{\otimes 4}$ contains
$\GC$-composition factors of degree at least $2016$ and so $\MD$ fails. Assume $d \leq 21$.
Then $d = 18$. According to \cite{Lu2}, the $\GC$-composition factor $L(\om_{1}+\om_{2})$ has
dimension $\geq 1104$ and so it is $G$-reducible. Thus $\MC$ and $\MDA$ fail for $G$.
\end{proof}

The trick of adding the center $\ZZ_{2}$ in the proof of Proposition \ref{sp44} is also used 
to construct the subgroups $(2 \times U_{3}(3)) \cdot 2$ in $Sp_{6}(\CC)$ and 
$(2 \times U_{5}(2)) \cdot 2$ in $Sp_{10}(\CC)$ recorded in Table I.
 
Next assume $S = G_{2}(q)$ with $q = 3^{f} \geq 9$. Then
$(G:N) \leq |\Out(\BS)| \cdot (S:N_{S}(Z)) = 2f(q^{6}-1)/(q-1)$, whence (\ref{bound3}) implies
$d < q^{2}(q^{2}+1) \leq \dd(\BS)$ (cf. \cite{H}), a contradiction.

Assume $S = F_{4}(q)$ with $q = 2^{f} \geq 4$. Then
$(G:N) \leq |\Out(\BS)| \cdot (S:N_{S}(Z)) = 2f(q^{12}-1)(q^{4}+1)/(q-1)$, whence (\ref{bound3})
implies $d < q^{7}(q^{3}-1)(q-1)/2 \leq \dd(\BS)$, a contradiction.

\subsection{Generic case.}

In this subsection, we consider the generic case, i.e. the quasi-simple Lie-type group
$S = E(G)$ is not of types $B_{2}$ or $F_{4}$ if $p = 2$, and type $G_{2}$ if $p = 3$. Then 
Lemma \ref{root} and (\ref{bound3}) imply that $\dim(V)$ satisfies (\ref{bound3}).
The proof of \cite[Thm. 3.1]{MMT} may now be repeated verbatim, and it shows that one of the 
following holds:

$(\star)$ $\BS = PSp_{2n}(q)$, $q = 3,5,7,9$, $V  \da_{S}$ is a Weil module of degree
$(q^{n} \pm 1)/2$;

$(\star\star)$ $\BS = U_{n}(q)$, $q = 2,3$, and $V \da_{S}$ is a Weil module of degree
$(q^{n}+q(-1)^{n})/(q+1)$ or $(q^{n}-(-1)^{n})/(q+1)$;

$(\star\star\star)$ $\BS= \Om^{\pm}_{2n}(2)$, $\tb D_{4}(3)$, $F_{4}(3)$, $E_{6}(2)$,
$\ta E_6(2)$, $\ta E_{6}(3)$, $E_{7}(2)$, or $E_{8}(2)$.

We mention right away that the proof of \cite[Thm. 3.1]{MMT} gave a specific way to deal with
quasi-simple groups $S = \Om^{\pm}_{2n}(2)$. This way also applies to the nearly simple groups
$G$ with $E(G) = S = \Om^{\pm}_{2n}(2)$, and shows that $G$ cannot satisfy $\MBA$. Now we
proceed to deal with the remaining cases individually.

\subsubsection{Symplectic groups.}\label{weil-sp}
Assume $\BS = PSp_{2n}(q)$ with $n \geq 2$, and $q = 3,5,7,9$. Notice that $V$ lifts to a
complex module $\VC$ in this case. Since $V$ is not stable under the outer diagonal
automorphism, $G$ cannot induce this automorphism on $S$. It follows that
$G = M$ if $q \neq 9$ and $(G:M) \leq 2$ if $q = 9$. This in turn implies by Lemma
\ref{g-m}(iii) that $X(V) \da_{S}$ is the sum of $1$ or $2$ irreducible
constituents of same dimension, for $X \in \{\TSB,\TWB\}$.

First we consider the case $q = 9$. Since $|\Out(S)| \leq 4 < d-1$ and $V \da_{S}$ is self-dual,
$V$ itself is self-dual. Suppose that $d = (9^{n}+1)/2$. Then the proof of \cite[Prop. 5.4]{MT1}
shows that $\SB(\VC) \da_{S} = 1_{S} + \al + \beta$, with $\al,\beta \in \Irr(S)$ and
$\al(1) = (9^{2n}-1)/16$. This implies that the $S$-module $\TSB(V)$ has a composition factor
of degree $\leq \al(1) < \dim(\TSB(V))/2$, a contradiction. Next suppose that $d = (9^{n}-1)/2$.
Then the proof of \cite[Prop. 5.4]{MT1} shows that $\WB(\VC) \da_{S} = 1_{S} + \gam + \delta$,
with $\gam,\delta \in \Irr(S)$ and $\gam(1) = (9^{n}+1)(9^{n}-9)/16$. This implies that the
$S$-module $\TWB(V)$ has a composition factor of degree $\leq \gam(1) < \dim(\TWB(V))/2$,
again a contradiction.

From now on we may assume that $q < 9$ and replace $G$ by $S$ by Lemma \ref{g-m}(iii). The
case $q = 7$ is excluded by \cite[Prop. 5.4]{MT1}. If $q = 3,5$, then $G$ satisfies
$\MBA$ at least in the case $\ell = 0$, see \cite[Prop. 5.4]{MT1}. We will show that $G$
cannot satisfy the condition $\MDA$ if $\ell \neq 2,3$. Assume the contrary. By Lemma
\ref{change1}, $S$ also satisfies $\MDA$ and we may replace $G$ by $S$.

Suppose $q = 5$. Then $V\da_{C} = A \oplus B \oplus B^{*}$, where $A = \Ker(c-1)$,
$B = \Ker(c - \eps)$ has dimension $5^{n-1}$ for $Z = \langle c \rangle$, and
$\eps \in \FF^{*}$ has order $5$. Now if $\ell \neq 2,3$ then Proposition \ref{bound2}(ii)
yields $\begin{pmatrix}(5^{n} \pm 1)/2 \\ 4 \end{pmatrix} < (S:C) = (5^{2n}-1)/2$, a
contradiction.  For $\ell = 2,3$, observe that $V^{(\ell)} \simeq V'$, where $V'$ is another
Weil module of same dimension. By \cite[Prop. 5.4]{MT1}, $(V \otimes V') \da_{S}$ is reducible.
Hence for the $\GC$-composition factor $W := L((\ell+1)\om_{1})$ of $V^{\otimes (\ell+1)}$ we
have $W \da_{S} = (V \otimes V^{(\ell)}) \da_{S} \simeq (V \otimes V')\da_{S}$ is reducible.

Suppose $q = 3$ and $n \geq 3$. Then $\VC\da_{C} = A \oplus B$, where $A = \Ker(c-1)$ has
dimension $(3^{n-1} \pm 1)/2$, $B = \Ker(c - \eps)$ has dimension $3^{n-1}$ for
$Z = \langle c \rangle$, and $\eps \in \FF^{*}$ has order $3$. Also, $V$ is not self-dual. Now
if $\ell \neq 2,3$ then Proposition \ref{bound2}(i) implies that either
$\begin{pmatrix}(3^{n} \pm 1)/2 + 3 \\ 4 \end{pmatrix} < (S:C) = (3^{2n}-1)/2$, or
$\AC(\VC)$ enters $1_{C} \ua S$ with multiplicity $\geq 2$. The former cannot hold as $n \geq 3$.
The latter is also impossible since 
$1_{C} \ua S = 1_{S} + \AC(\xi) + \AC(\eta) + \xi\bar{\eta} + \bar{\xi}\eta$ is a sum of $5$ distinct 
$S$-irreducibles, where $\VC$ is one the four Weil characters $\xi,\bar{\xi},\eta,\bar{\eta}$ of 
$S$, see \cite[p. 255]{MT1}. Assume $\ell = 2$. Then $\dim(V) = (3^{n}-1)/2$, and 
$\Vbb \simeq V^{*}$. Hence for the $\GC$-composition factor $W := L(3\om_{1})$ of $V^{\otimes 3}$ 
we have $W \da_{S} = (V \otimes \Vbb) \da_{S} \simeq (V \otimes V^{*})\da_{S}$ is reducible.

We have proved

\begin{propo}\label{sp-odd}
{\sl Assume $G \leq \GC$, $G$ satisfies $\MBA$ and $\BS \simeq PSp_{2n}(q)$ with $q$ odd.
Assume in addition that $(\ell,q) = 1$, $n \geq 2$ and $(n,q) \neq (2,3)$. Then $q = 3,5$,
$G = Z(G)S$, and $V \da_{S}$ is a Weil module of dimension $(q^{n} \pm 1)/2$. If $\ell \neq 2,3$,
or if $q = 5$ and $\ell = 3$, then $G$ cannot satisfy $\MDA$. If $\ell = 2$ then $G$ cannot
satisfy $\MCA$.
\hfill $\Box$}
\end{propo}

We also record the following remark.

\begin{lemma}\label{sp-m4}
{\sl Let $G = Sp_{2n}(q)$ with $n \geq 2$ and $q$ odd. Let $V$ be a complex Weil module for 
$G$. Then $M_{4}(G,V) \geq \lfloor (q+7)/4 \rfloor$.}
\end{lemma}

\begin{proof}
Consider the case $q \equiv 1 (\mod 4)$. According to \cite{MT1}, the dimension of any irreducible 
constituent $Y$ of $V \otimes V$ is at most $d_{1} := (q^{2n}-1)/(q-1)$; moreover, if 
$\dim(Y) < (q^{2n}-1)/(q-1)$ then $\dim(Y) \leq d_{2} := (q^{n}-1)(q^{n}+q)/2(q-1)$. Since
$V$ is self-dual (cf. Table II), $V \otimes V$ contains $1_{G}$. Thus the number $N$ of 
irreducible constituents of $V \otimes V$ is at least $1 +(\dim(V)^{2}-1)/d_{1}$, 
and so $N \geq (q+3)/4$. Assume $M_{4}(G,V) < \lfloor (q+7)/4 \rfloor$. Then $N = (q+3)/4$. But 
this is impossible as neither $1 + d_{1}(q-1)/4$ nor $1 + d_{2} + d_{1}(q-5)/4$ equals 
$\dim(V)^{2}$. The case $q \equiv 3 (\mod 4)$ is similar, with $d_{1}$ replaced by 
$(q^{2n}-1)/(q+1)$ and $d_{2}$ replaced by $(q^{n}+1)(q^{n}+q)/2(q+1)$.  
\end{proof}

To show that Weil representations cannot satisfy $\MC$, we need the following statement:

\begin{lemma}\label{uni-q}
{\sl Let $q = p^{f}$ for a prime $p$ and let $\chi$ be a complex irreducible unipotent character of 
a finite group of Lie type $G$ in characteristic $p$.

{\rm (i)} Assume $G = GL_{n}(q)$ or $GU_{n}(q)$. Then the $p$-part $\chi(1)_{p}$ of $\chi(1)$ is 
a power $q^{e}$ of $q$. Moreover, either $e = 0$ and $\chi$ is trivial, or $e = 1$ and $\chi$ is 
labeled by the partition $(n-1,1)$, or $e \geq 2$.

{\rm (ii)} Assume $p > 2$ and $G$ is of type $B_{n}$, $C_{n}$, $D_{n}$ or $\ta D_{n}$. Then 
either $\chi$ is trivial, or $\chi(1)_{p}$ is a nontrivial power of $q$.}
\end{lemma}

\begin{proof}
Use Lusztig's classification of unipotent characters of $G$ (cf. \cite{C}), and repeat the proof of 
\cite[Lemma 7.2]{MMT}.
\end{proof}

\begin{propo}\label{sp-m6}
{\sl Let $\BS = PSp_{2n}(q)$ with $n \geq 2$ and $q = 3,5$. Let $V|_{S}$ be a complex Weil module 
for $S := E(G)$. Then $G$ fails $\MC$.}
\end{propo}

\begin{proof}
Assume the contrary. Since the Weil modules of $Sp_{2n}(q)$ are not stable under outer automorphisms
of $Sp_{2n}(q)$, without loss we may assume that $G = Z(G)S$. 

1) First we consider the case $q = 3$. Observe that $\GC$-modules $(V^{*})^{\otimes 3}$ 
and $V \otimes V \otimes V^{*}$ have no common irreducible constituent, hence the same is true for 
these two modules considered over $G$. However, by \cite[Prop. 5.4]{MT1}, there is 
$X \in \{\SB,\WB\}$ such that $X(V) \simeq X(V^{*})$. It follows that $(V^{*})^{\otimes 3}|_{S}$ and
$(V \otimes V \otimes V^{*})|_{S}$ have the common summand $X(V) \otimes V^{*}$, a contradiction. 

2) From now on we may assume $q = 5$. By Lemma \ref{change2}, $V$ is self-dual, of type $\eps = \pm$ 
for $d = (5^{n}+\eps 1)/2$. Since $G$ satisfies $\MC$, $G$ is irreducible on $\TSC(V)$ (the largest
quotient of the $\GC$-module $\SC(V)$) and $\TWC(V)$ (the largest quotient of the $\GC$-module 
$\WC(V)$). However, if $2 \leq n \leq 5$ then at least one of $\TSC(V)$, $\TWC(V)$ has dimension not
divisible by $|G|$. So we may assume $n \geq 6$. If $\eps = +$ then 
$\dim(\TWC(V)) = (5^{2n}-1)(5^{n}-3)/48$. If $\eps = -$ then 
$\dim(\TSC(V)) = (5^{2n}-1)(5^{n}+3)/48$. 

It suffices to show that $S$ has no irreducible character of degree $D := (5^{2n}-1)(5^{n}-3\eps)/48$. 
Assume the contrary: $\chi(1) = D$ for some $\chi \in \Irr(S)$. Consider the dual group 
$G^{*} := SO_{2n+1}(q)$ and its natural module $U := \FF_{q}^{2n+1}$. By Lusztig's classification of
irreducible characters of $G$ (cf. \cite{C}), $\chi$ corresponds to the $G^{*}$-conjugacy class of a 
semisimple element $s \in G^{*}$, and a unipotent character $\psi$ of $C := C_{G^{*}}(s)$; moreover,
$\chi(1) = E\psi(1)$, where $E := (G^{*}:C)_{p'}$. Notice that $C$ is a subgroup of index $\leq 2$ of 
a direct product $\hat{C} := C_{O(U)}(s)$ of groups of form $GL_{m}(q^{k})$, $GU_{m}(q^{k})$, and 
$O^{\pm}_{m}(q)$. Since $\chi(1) = D$ is coprime to $5$, Lemma \ref{uni-q} implies that $\psi(1) = 1$. 
Thus $E = D$; in particular, $s \neq 1$.

Claim that if $U = U_{1} \oplus U_{2}$ is any decomposition of $U$ into $C$-invariant nonzero 
nondegenerate subspaces, then $\{\dim(U_{1}),\dim(U_{2})\} = \{1,2n\}$ or $\{2,2n-1\}$. Otherwise we 
may assume $2n-3 \geq \dim(U_{1}) = 2k+1 \geq 3$. If $2 \leq k \leq n-2$ then $E \geq q^{4n-8} > D$. If 
$k = 1$, then $E$ is divisible by $(5^{2n}-1)(5^{n-1} \pm 1)/48$, which does not divide $D$. 

It is shown in \cite[p. 2124]{TZ1} that $s$ has an eigenvalue $\lam = \pm 1$ on $U$. Observe that 
$\Ker(s-\lam)$ is a $C$-invariant proper nondegenerate subspace of $U$. The above claim now implies
that $\hat{C} = A \times O^{\pm}_{2n-a}(q)$ with $A \leq O^{\pm}_{a+1}(q)$ and $a = 0,1$, or
$\hat{C} = \ZZ_{2} \times GL^{\pm}_{m}(q^{k})$ with $n = mk$. In the former case, 
$E \leq 2(5^{2n}-1) < D$. In the latter case $E > q^{n^{2}/2} > D$, a contradiction.   
\end{proof}

\subsubsection{Unitary groups.}\label{weil-su}
Assume $\BS = U_{n}(q)$ with $n \geq 5$, and $q = 2,3$. Notice that $V$ lifts to a
complex module $\VC$ in this case.

\medskip
{\bf Case $U_{n}(3)$.}
Recall that the (reducible) Weil character
$$\omega_{n}~:~g \mapsto (-3)^{\dim_{\FF_{9}}\Ker(g-1)}$$
(where $\Ker(g-1)$ is the fixed point subspace of $g$ in the natural module $\FF_{9}^{n}$)
of $GU_{n}(3)$ is the sum of four irreducible complex Weil characters: $\al$ of degree
$(3^{n}+3(-1)^{n})/4$, and $\beta,\overline{\beta},\gam$ all of degree $(3^{n}-(-1)^{n})/4$.
All of them restrict irreducibly to $SU_{n}(3)$, and among them, $\al$ and $\gam$ are
real-valued. The irreducibility question of tensor products,
tensor squares, symmetric and alternating squares of complex Weil modules of $SU_{n}(3)$ has
been studied in \cite{MT1,LST}. The following result completely settles this question:

\begin{theor}\label{su3}
{\sl Let $W := \{\al,\beta,\overline{\beta},\gam\}$ be the set of irreducible complex
characters of $GU_{n}(3)$ with $n \geq 3$. Assume $H \in \{G := GU_{n}(3), S:= SU_{n}(3)\}$ and
$$\rho \in \{\AC(\chi),\TSB(\chi),\TWB(\chi),\chi\lam \mid \chi,\lam \in W\}.$$

{\rm (A)} Then $\rho$ is irreducible over $H$ if and only if one of the following holds.

\hspace{0.5cm}{\rm (i)} $\rho \in \{\SB(\beta),\SB(\overline{\beta})\}$ if $n$ is even, and
$\rho \in \{\WB(\beta),\WB(\overline{\beta})\}$ if $n$ is odd.

\hspace{0.5cm}{\rm (ii)} $\rho = \WB(\al)$ if $n$ is even, and $\rho = \SB(\al)$ if $n$ is odd.

\hspace{0.5cm}{\rm (iii)} $\rho = \TWB(\gam)$ if $n$ is even, and $\rho = \TSB(\gam)$ if $n$ is
odd.

\hspace{0.5cm}{\rm (iv)} $\rho = \TWB(\al)$ and $n = 3$.

{\rm (B)} Assume $n \geq 4$, and $X(\mu) \da_{S}$ is reducible  either for $(X,\mu) = (A,\beta)$,
or for $X \in \{\TSB,\TWB\}$ and $\mu \in \{\al,\gam\}$. Then $X(\mu) \da_{S}$
contains irreducible constituents of distinct degrees.}
\end{theor}

\begin{proof}
1) The case $n = 3,4$ can be checked directly using \cite{Atlas}, so we assume $n \geq 5$.
By \cite[Prop. 4.1]{MT1}, $\chi\lam$ is reducible over $H$. It is clear that $\AC(\al)$ and
$\AC(\gam)$ are reducible. We claim that $\AC(\beta)$ is also reducible.
Indeed, $\AC(\beta)$ is a subquotient of $\omega_{n}^{2}$. The proof of \cite[Prop. 3.3]{LST} shows
that $H$-composition factors of $\omega_{n}^{2}$ have degree $1$,
$e_{0} = c_{n}c_{n-1}/8$, $e_{1} = 27c_{n-1}c_{n-2}/32$, $e_{2} = 3e_{0}/4$,
$e_{3} = 3c_{n}c_{n-2}/16$, $e_{4} = e_{0}/2$, $e_{5} = 9c_{n}c_{n-3}/32$, and
$e_{6} = e_{0}/4$, where $c_{k} := (3^{k}-(-1)^{k})$ for any $k$. Now our claim follows
since all these degrees are less than $\dim(\AC(\beta))$.

2) Over $G$ we have
$$\omega_{n}^{2} = (\al + \beta + \overline{\beta} + \gam)^{2} = 4 \cdot 1_{G} + \TSB(\al) +
  \TWB(\al) + \TSB(\gam) + \TWB(\gam) +$$
$$ + \TSB(\beta) + \TWB(\beta) + \TSB(\overline{\beta}) + \TWB(\overline{\beta}) +
  2\AC(\beta) + 2\al\gam + 2\beta\gam +
  2\overline{\beta}\gam + 2\al\beta + 2\al\overline{\beta}.$$
Let $A_{1}$, $A_{2}$, $C_{1}$, $C_{2}$, $B_{1}$, $B_{2}$, $B_{3} = B_{1}$, $B_{4} = B_{2}$,
$2E$, $2D_{1}$, $2D_{2}$, $2D_{3} = 2D_{2}$, $2D_{4}$, and $2D_{5} = 2D_{4}$, denote the number of
$S$-irreducible constituents (with counting multiplicities) of the $14$ nontrivial
summands in this decomposition, in the order they are listed.
Also set $A = A_{1} + A_{2}$, $B = B_{1} + B_{2}$,
$C = C_{1} + C_{i}$, and $D = \sum^{5}_{i=1}D_{i}$. According to 1), $E \geq 2$, and
$D_{i} \geq 2$ for any $i$ so $D \geq 10$. By \cite[Prop. 3.3]{LST} and its proof,
$A_{1} + A_{2} \geq 3$, namely $A_{1} = 2$ if $n$ is even and $A_{2} = 2$ if $n$ is odd, and
$C_{2} = 1$ if $n$ is even and $C_{1} = 1$ if $n$ is odd.

It is shown in \cite{T} that $\omega_{n}^{2}$ is the sum of $40$ irreducible $H$-characters (with
counting multiplicities). Hence $2B = 40-(4+A+C+2D+2E) \leq 7$, i.e. $B \leq 3$. If $B \leq 2$
then both $\SB(\beta)$ and $\WB(\beta)$ are irreducible over $S$, whence $\AC(\beta)$ is also
irreducible, contrary to 1). Thus $B = 3$. Now $2E = 40-(4+A+C+2B+2D) \leq 5$ and $E \geq 2$,
so $E = 2$. Next, $5 \leq A+C = 40-(4+2B+2D+2E) = 26-2D \leq 6$, so $A+C = 6$ and $D = 10$.
In particular, $D_{i} = 2$ for all $i$.

Claim that $C = 3$. Assume the contrary: $C = 2$. First suppose that $n$ is odd. Then both
$\TSB(\gam)$ and $\WB(\gam)$ are irreducible of degree at least $a(a+1)/2$, where
$a := (3^{n}+1)/4$. Since $A = 4$ and $A_{2} = 2$, we get $A_{1} = 2$, and all composition
factors of $\TSB(\al)$ and $\TWB(\al)$ are of degree less than the degree of $\SB(\al)$ which
is $a(a+1)/2$. It follows that
$$(\al\gam,\al\gam)_{S} = (\gam^{2},\al^{2})_{S} =
  (1_{S}+\TSB(\gam)+\TWB(\gam),1_{S}+\TSB(\al)+\TWB(\al))_{S} = 1,$$
i.e. $D_{1} = 1$, a contradiction. Next, suppose $n$ is even. Then
$\TSB(\gam)$ is irreducible of degree $e_{2}$ and $\TWB(\gam)$ is irreducible
of degree $e_{1}$. The proof of \cite[Prop. 3.3]{LST} shows that
$A_{1} = 2$ and $\TSB(\al)$ is the sum of two characters of degree $e_{4}$ and
$e_{5}$. Hence $A_{2} = 2$ and $\TWB(\al)$ is reducible, of the same
degree as of $\TSB(\gam)$. On the other hand, $D_{1} = 2$ implies
$$2 = (\al\gam,\al\gam)_{S} = (\gam^{2},\al^{2})_{S} =
  (1_{S}+\TSB(\gam)+\TWB(\gam),1_{S}+\TSB(\al)+\TWB(\al))_{S}.$$
The above degree consideration shows that $\TWB(\gam)$ has to be a constituent of
$\TWB(\al)$. It follows that $\TWB(\al)$, and so $\omega_{n}^{2}$ as well, has a nontrivial
$S$-irreducible constituent of degree at most
$\TWB(\al)(1)-\TWB(\gam)(1) = (3^{n}+3)/4$. But this cannot happen, see \cite[p. 196]{LST}.

Thus $C = 3 = A$, and moreover $(A_{1},A_{2}) = (C_{1},C_{2}) = (1,2)$ if $n$ is odd, and
$(2,1)$ if $n$ is even. It remains to determine the $B_{i}$'s, which amounts to show that
$B_{1} > 1$ if $n$ is odd and $B_{2} > 1$ if $n$ is even. Assume $n$ is odd but $B_{1} = 1$,
or $n$ is even but $B_{2} = 1$. Then $\omega_{n}^{2}$ has an $S$-irreducible constituent of degree
$(3^{n}-(-1)^{n})(3^{n}-5(-1)^{n})/32$, which is none of the $e_{i}$'s, a contradiction.

3) To prove (B), we assume that $X(\mu)$ is as in (B) but all $S$-composition factors of
$X(\mu) \da_{S}$ are of same degree say $d'$. By the results of 2),
$d' = \dim(X(\mu))/2$. If $(X,\mu) = (A,\beta)$, then $d' = (d^{2}-1)/2$. If $\mu = \al$, then
$X = \TSB$ if $n$ is even and $X = \TWB$ if $n$ is odd, and
$d' = (3^{n}-(-1)^{n})(3^{n}+11(-1)^{n})/64$. If $\mu = \gam$, then $X = \TSB$ if $n$ is even
and $X = \TWB$ if $n$ is odd, and $d' = 3e_{0}/8$. In all cases, none of the $S$-composition
factors of $\omega_{n}^{2}$ can have this degree $d'$, a contradiction.
\end{proof}

Now suppose that $G < \GC$ satisfies $\MBA$. Then $G$ is irreducible on $X(V)$ for some
$X \in \{\AC,\TSB,\TWB\}$. Hence the $G$-module $X(\VC) (\mod \ell)$ has an irreducible
constituent (namely $X(V)$) of degree $\geq \dim(X(\VC)) -2$, so the same is true for
$X(\VC)$. Let $U$ be any other irreducible constituent of the $G$-module $X(\VC)$. Then
$\dim(U) \leq 2$. Since $\dd(\BS) > 2$, it follows that $S$ acts trivially on $U$.

First assume that $\GC = GL(V)$ and let $X = \AC$. Since $\al$ and $\gam$ are self-dual and
$|\Out(\BS)| < d-1$, by Lemma \ref{rare} $V \da_{S}$ is not self-dual, whence we may assume
that $\VC \da_{S}$ affords the character $\beta$. As $S$ acts irreducibly on $\VC$, $1_{S}$
enters $\VC \otimes \VC^{*}$ with multiplicity $1$, and moreover it cannot enter $\AC(\VC)$. Thus
the $G$-module $\AC(\VC)$ is irreducible. By Theorem \ref{su3}(A), $S$ is reducible on
$\AC(\VC)$. Now by Theorem \ref{su3}(B), $\AC(\VC) \da_{S}$ has irreducible constituents of
distinct degrees, which contradicts Clifford's theorem.

Next we consider the case $\ell \neq 2,3$ and $\GC = Sp(V)$ or $O(V)$. Set $Y = \SB$ if
$X = \TSB$, and $Y = \WB$ if $X = \TWB$. Since $\beta (\mod \ell)$ is not self-dual, we may
assume that $\VC \da_{S}$ affords the character $\al$ or $\gam$. As $S$ acts irreducibly on
$\VC$ and $\ell \neq 2$, $\VC$ has the same types as a $G$-module or as an $S$-module.
Hence $1_{S}$ enters $Y(\VC)$ with multiplicity $\leq 1$, and moreover it cannot enter
$X(\VC)$. Thus the $G$-module $X(\VC)$ is irreducible. Now choose $X = \TSB$ if $n$ is even
and $X = \TWB$ if $n$ is odd. By Theorem \ref{su3}(A), $S$ is reducible on $X(\VC)$. By
Theorem \ref{su3}(B), $X(\VC) \da_{S}$ has irreducible constituents of distinct degrees,
contrary to Clifford's theorem.

Finally we consider the case $\ell = 2$ and $\GC = Sp(V)$ or $O(V)$. Then $X = \TWB$ and
$Z(G) = 1$. Since $\ell = 2$, we may assume that $\VC \da_{S}$ affords the character $\al$
if $n$ is odd and $\gam$ if $n$ is even. We claim that $\TWB(V)$ is irreducible if and only
if $n = 3$. The case $n = 3,4$ can be checked directly, so we suppose that $n \geq 5$ but
$\TWB(V)$ is irreducible. Let $K$ be the subgroup consisting of the elements of $G$ that
induce only inner-diagonal automorphisms on $S$. Then $(G:K) \leq 2$ and $K$ can be embedded
in $PGU_{n}(3)$. First assume $n$ is even. Then we may assume that $\VC \da_{S}$ affords the
character $\gam$. The proof of \cite[Prop. 3.3]{LST} shows that
$\SB(\al)$ is the sum of three irreducible $K$-characters of degree $1$,
$e_{4} = (3^{n}-1)(3^{n-1}+1)/16$, and $e_{5} = (3^{n}-1)(3^{n-1}+9)/32$. The same is true
for $G$, since $K \lhd G$ and $\al$ is $G$-stable. Since $\al (\mod 2) = \gam (\mod 2) + 1$,
we have $\SB(\al) (\mod 2) = \SB(\gam) (\mod 2) + \gam (\mod 2) + 1$. Also,
$\SB(\gam (\mod 2)) = \SB(V) = \WB(V) + \Vbb = \TWB(V) + \Vbb + 2$ since
$4|d$. Thus the irreducible constituents of the $G$-character $\SB(\al) (\mod 2)$ are of
degree $1$, $d = (3^{n}-1)/4$, and $\dim(\TWB(V)) = (3^{n}-1)(3^{n}-5)/32$. This yields a
contradiction as $\dim(\TWB(V))$ is larger than $e_{4}$ and $e_{5}$. Now we assume that
$n$ is odd. Then we may assume that $\VC \da_{S}$ affords the character $\al$. The proof of
\cite[Prop. 3.3]{LST} shows that $\WB(\al)$ is the sum of three irreducible $K$-characters
of degree $1$, $e_{4} = (3^{n}+1)(3^{n-1}-1)/16$, and $e_{5} = (3^{n}+1)(3^{n-1}-9)/32$.
The same is true for $G$, since $K \lhd G$ and $\al$ is $G$-stable. Since
$\dim(V) =  \al(1)$ and $\dim(\TWB(V)) = \dim(\WB(\al)) - 2$, it follows that
$\TWB(V)$ cannot be irreducible, a contradiction. We record this in the following statement,
which complements \cite[Thm. 3.1(ii)]{MMT}:

\begin{propo}\label{su3-2}
{\sl Let $G$ be a finite group with normal subgroup $S = U_{n}(3)$ with $n \geq 4$ and
a (faithful) irreducible $\overline{\FF}_{2}G$-module $V$ such that $V \da_{S}$ is a
Weil module. Then $\TWB(V)$ is not irreducible.
\hfill $\Box$}
\end{propo}

\medskip
{\bf Case $U_{n}(2)$, $n \geq 5$.}
The (reducible) Weil character
$$\omega_{n}~:~g \mapsto (-2)^{\dim_{\FF_{4}}\Ker(g-1)}$$
of $GU_{n}(2)$ is the sum of three irreducible complex Weil characters: $\al = \zeta^{0}_{n,2}$
of degree $(2^{n}+2(-1)^{n})/3$, and
$\beta = \zeta^{1}_{n,2},\overline{\beta} = \zeta^{2}_{n,2}$ of degree $(2^{n}-(-1)^{n})/3$.
All of them restrict irreducibly to $SU_{n}(2)$. According to \cite[Prop. 3.3]{LST},
all these characters satisfy $M_{4}(\GCC)$. However, we will show that $G$ cannot satisfy
$\MDA$. Assume the contrary.

First we consider the case $\ell \neq 2,3$. Assume $\GC = GL(V)$. Let $Z \simeq \ZZ_{2}$
be a long-root subgroup of $S$, $C := C_{G}(Z) = N_{G}(Z)$. Since
$Z \not\leq Z(G)$, we have $\VC \da_{C} = A \oplus B$ with $A = C_{\VC}(Z)$ and
$B = [\VC,Z]$ both being nonzero. Thus the assumptions of Proposition \ref{bound2}(i) are
fulfilled. By Lemma \ref{root}, $(G:C) = (S:C \cap S)$, hence
$(1_{C} \ua G) \da_{S} = 1_{C \cap S} \ua S$, and moreover the last character is a rank $3$
permutation character (indeed, it is $1_{S} + \TWB(\al) + \AC(\beta)$ if $n$ is odd and
$1_{S} + \TSB(\al) + \AC(\beta)$ if $n$ is even). It follows that $\AC(\VC)$ cannot enter
$1_{C} \ua G$ with multiplicity $\geq 2$. So
$\begin{pmatrix}d+3 \\4 \end{pmatrix} \leq (G:C)$, which is a contradiction since $n \geq 5$,
$d \geq (2^{n}-2)/3 \geq 10$ and
$(G:C) = (S:C \cap S) = (2^{n}-(-1)^{n})(2^{n-1}-(-1)^{n-1})/3$. Next assume that
$\GC = Sp(V)$, resp. $O(V)$. Since $\ell \neq 3$, $\beta (\mod \ell)$ is not self-dual, whence
we may assume that $\VC \da_{S}$ affords the character $\al$. Again
$d \geq (2^{n}-2)/3 \geq 10$. By Proposition \ref{v4}, the condition $\MDA$ implies that
$\SD(\VC)$, resp. $\WD(\VC)$, is irreducible. Let $D := SU_{n-1}(2)$ and $C := Z(G)D$. Then
$\al_{D} = \zeta^{1}_{n-1,2} + \zeta^{2}_{n-1,2}$, so $\VC \da_{C} = B \oplus B^{*}$ with
$\dim(B) = (2^{n-1}-(-1)^{n-1})/3 \geq 5$. The proof of Proposition \ref{bound2}(ii) now
implies that $\begin{pmatrix}d \\4 \end{pmatrix} < (G:C)$, which is a contradiction when
$n \geq 7$ since $(G:C) \leq |\Out(\BS)| \cdot (S: S \cap C) = 6(2^{n}-(-1)^{n})2^{n-1}$.
When $n = 5,6$, we can check directly that $\SD(\VC)$, resp. $\WD(\VC)$, is reducible.

Next we consider the case $\ell = 3$. Assume $n$ is odd. Then
$\beta (\mod 3) = \al (\mod 3) + 1$, so we may assume that $\VC \da_{S}$ affords the character
$\al$. Since $\QQ(\al) = \QQ$, $V \da_{S} \simeq \Vcc \da_{S}$ and moreover
$V \da_{S}$ is self-dual. The condition $\MDA$ implies that $G$ is irreducible on
$W = L(4\om_{1}) \simeq V \otimes \Vcc$. However,
$W \da_{S} \simeq (V \otimes \Vcc) \da_{S} \simeq (V \otimes V)\da_{S}$ contains $1_{S}$,
and so $W$ is reducible, a contradiction. Next assume that $n$ is even. Then
$\beta (\mod 3) = \al (\mod 3) - 1 = \overline{\beta} (\mod 3)$, so we may assume that
$V \da_{S}$ affords the character $\al (\mod 3)-1$. Since $\QQ(\al) = \QQ$,
$V \da_{S} \simeq \Vcc \da_{S}$ and moreover $V \da_{S}$ is self-dual. Arguing as in the
case $n$ is odd, we conclude that $G$ is reducible on $L(4\om_{1})$.

We have proved

\begin{propo}\label{su-m8}
{\sl Assume $G \leq \GC$, $G$ satisfies $\MBA$ and $\BS \simeq U_{n}(q)$ with $n \geq 3$.
Assume in addition that $(\ell,q) = 1$ and $(n,q) \neq (3,2)$, $(3,3)$, $(4,2)$, $(4,3)$,
$(6,2)$. Then $q = 2$ and $V \da_{S}$ is a Weil module of dimension $(2^{n} + 2(-1)^{n})/3$
or $(2^{n}-(-1)^{n})/3$. Moreover, $G$ cannot satisfy $\MDA$.
\hfill $\Box$}
\end{propo}

\begin{propo}\label{su-m6}
{\sl Let $\BS = U_{n}(2)$ with $n \geq 4$ and let $V|_{S}$ be a complex Weil module for 
$S := E(G)$. Then $G$ fails $\MC$, except for the case $G = (\ZZ_{2} \times U_{5}(2)) \cdot 2$ inside
$\GC := Sp_{10}(\CC)$ in which case $G$ fails $\MD$.}
\end{propo}

\begin{proof}
1) Assume the contrary and let $\eps := -1$. Since $G$ satisfies $\MC$, $G$ is irreducible on 
$\TSC(V)$ and $\TWC(V)$. However, if $n = 4$, $6$, or $7$, then at least one of $\TSC(V)$, 
$\TWC(V)$ has dimension not divisible by $|G|$. If $n = 5$ then direct computation shows 
that $G$ satisfies $\MC$ precisely when $G = (\ZZ_{2} \times U_{5}(2)) \cdot 2$ inside
$\GC := Sp_{10}(\CC)$ (but no proper subgroup of $G$ can satisfy $\MC$), and $\MD$ fails. 

From now on we assume $n \geq 8$. If $2|n$, then we consider 
$D := \dim(\TSC(V))$ which is  $(2^{n}-1)(2^{n}+2)(2^{n}+5)/162$ for $d = (2^{n}-1)/3$ and 
$(2^{n}-1)(2^{n}+2)(2^{n}+14)/162$ for $d = (2^{n}+2)/3$. If $n$ is odd, then we consider 
$D := \dim(\TWC(V))$ which is  $(2^{n}+1)(2^{n}-2)(2^{n}-5)/162$ for $d = (2^{n}+1)/3$ and 
$(2^{n}+1)(2^{n}-2)(2^{n}-14)/162$ for $d = (2^{n}-2)/3$. These dimension formulae follow
from Table II and Lemma \ref{change2}. Thus 
$$D := (2^{n}-\eps^{n})(2^{n-1}+\eps^{n})(2^{n}+\gam\eps^{n})/81$$ 
with $\gam \in \{5,14\}$. It suffices to show that $G$ has no irreducible character of degree 
$D$.  

Recall that for $X \in \{\SC,\WC\}$, $\tilde{X}(V)$ is either $X(V)$ or $X(V)/V$. Also, 
the Weil representation $V$ of $SU_{n}(2)$ extends to $H := GU_{n}(2)$. Hence, by 
\cite[Lemma 2.2]{MT1}, the $H$-module $X(V)$ has a quotient of dimension equal to $D$. 
So without loss we may assume that $G \geq H$. Since $(G:Z(G)H) \leq 2$, it suffices to show 
that $H$ has no irreducible character of degree $D$ or $D/2$. Assume the contrary: 
$\chi(1) \in \{D,D/2\}$ for some $\chi \in \Irr(H)$. We can identify the dual group $H^{*}$ with 
$H$ and consider its natural module $U := \FF_{4}^{n}$. By Lusztig's classification of 
irreducible characters of $G$, $\chi$ corresponds to the $H$-conjugacy class of a semisimple 
element $s \in H$, and a unipotent character $\psi$ of $C := C_{H}(s)$; moreover, 
$\chi(1) = E\psi(1)$, where $E := (H:C)_{2'}$. Notice that $C$ is a direct product of groups of 
form $GU_{m}(2^{k})$ (with $k$ odd) and $GL_{m}(4^{k})$, and that $4 \not{|} D$. 

2) Claim that if $U = W_{1} \oplus W_{2}$ is any decomposition of $U$ into $C$-invariant nonzero 
nondegenerate subspaces, then $\{\dim(W_{1}),\dim(W_{2})\} = \{1,n-1\}$ or $\{2,n-2\}$. 
Otherwise we may assume $n/2 \geq \dim(W_{1}) = k \geq 3$. If $k \geq 4$ and 
$n \geq 9$ then $E > D$. If $k = 4$ and $n = 8$ then $D/E$ is not an integer. If $k = 3$, then 
$E$ is divisible by $(2^{n}-\eps^{n})(2^{n-1}+\eps^{n})(2^{n-2}-\eps^{n})/27$, which does 
not divide $D$. In any of these cases we get a contradiction.

Now assume that $U = \oplus^{t}_{i=1}V_{i}$ is any decomposition of $U$ into $C$-invariant 
nondegenerate subspaces, with $1 \leq \dim(V_{1}) \leq \ldots \leq \dim(V_{t})$, and with 
$t$ as large as possible. The above claim implies that $t \leq 3$. Moreover, if $t = 3$ 
then $\dim(V_{1}) = \dim(V_{2}) = 1$, and if $t = 2$ then $\dim(V_{1}) = 1$ or $2$.

3) With the notation as in 2), here we consider the case $t = 3$. Then 
$C = GU_{1}(2) \times GU_{1}(2) \times A$, with $A = C_{GU(V_{3})}(s|_{V_{3}})$. In particular,
$E = F \cdot (2^{n}-\eps^{n})(2^{n-1}+\eps^{n})/9$, where $F := (GU_{n-2}(2):A)_{2'}$. Using 
the results of \cite{TZ1}, one can show that either $F = 1$, or $F = (2^{n-2}-\eps^{n})/3$, or 
$F \geq (2^{n-2}-\eps^{n})(2^{n-3}-4)/9$. Since $E|D$, the last two possibilities cannot 
occur. Thus $F = 1$ and so $A = GU_{n-2}(2)$. Since $4 \not{|} D$, by Lemma \ref{uni-q} either 
$\psi(1) = 1$ or $\psi(1) = (2^{n-2}+2\eps^{n})/3$. In none of these cases can $\chi(1)$ be 
equal to $D$ or $D/2$, a contradiction.    

Arguing similarly, we see that the case $t = \dim(V_{1}) = 2$ cannot happen.

Assume $t = 2$ and $\dim(V_{1}) = 1$. By the choice of $t$, $V_{2}$ cannot be decomposed
into a direct sum of $C$-invariant proper nondegenerate subspaces. It follows that 
$C = GU_{1}(2) \times A$, and either $A = GU_{n-1}(2)$, or $A = GU_{m}(2^{k})$ with 
$mk = n-1$ and $k \geq 3$ odd, or $A = GL_{m}(4^{k})$ with $2mk = n-1$. In the last two cases one
can check that $E > D$ (as $n \geq 8$). So $A = GU_{n-1}(2)$. Since $4 \not{|} D$, by Lemma 
\ref{uni-q} either $\psi(1) = 1$ or $\psi(1) = (2^{n-1}-2\eps^{n})/3$. In none of these cases
can $\chi(1)$ be equal to $D$ or $D/2$, a contradiction.         

4) We have shown that $t = 1$. By the choice of $t$, $V$ cannot be decomposed
into a direct sum of $C$-invariant proper nondegenerate subspaces. It follows that 
either $C = GU_{n}(2)$, or $C = GU_{m}(2^{k})$ with $mk = n$ and $k \geq 3$ odd, or 
$A = GL_{m}(4^{k})$ with $2mk = n$. In the last two cases one can check that $E > D$ (as 
$n \geq 8$). So $C = GU_{n}(2)$ and $E = 1$. Since $4 \not{|} D$, by Lemma 
\ref{uni-q} either $\psi(1) = 1$ or $\psi(1) = (2^{n}+2\eps^{n})/3$. In both of these cases
$\chi(1) < D/2$, again a contradiction. 
\end{proof}

\subsection{Some exceptional groups.}
In this subsection we handle the exceptional groups listed in $(\star\star\star)$. Our results
complement \cite[Thm. 3.1]{MMT}.

\begin{propo}\label{e678}
{\sl Assume $G \rhd S = E_{m}(2)$ with $m = 6,7,8$ and $\ell \neq 2$. Then $G$ has no faithful
irreducible $\ell$-modular representation $V$ such that $X(V)$ is irreducible, if $X = \AC$, resp.
$\TSB$, $\TWB$, for $V$ of type $\circ$, resp. $+$, $-$.}
\end{propo}

\begin{proof}
The proof is the same for $m = 6$, $7$, or $8$, so we give the details for $m = 8$.
Assume the contrary. Then we are in the so-called {\it good case} considered in \cite{MMT}.
Let $Z$ be a long-root subgroup of $S$ and $N := N_{G}(Z)$. By Lemma \ref{root},
$(G:N) = (S:N \cap S)$. The proof of \cite[Thm. 3.1]{MMT} yields
\setcounter{equation}{0}
\begin{equation}\label{e8-1}
  \dim(V) \leq 765,625,740.
\end{equation}
On the other hand,
\begin{equation}\label{e8-2}
  \dim(V) \geq 545,925,247
\end{equation}
by \cite{Ho}. Since $\sqrt{(G:N)+2} < 541,379,153 < \dim(V)$, the proof of
\cite[Thm. 3.1]{MMT} shows that $V$ is self-dual.

Recall that $N = QL$ with $Q = 2^{1+56}$, $Z = Z(Q) = \langle z \rangle$, and $L = E_{7}(2)$. We
write $V \da_{N} = \oplus^{2}_{i=0}V_{i}$ with $V_{0} := C_{V}(Q)$, $V_{1} := [C_{V}(Z),Q]$, and
$V_{2} := [V,Z]$. Since $Z \not\leq Z(G)$, $V_{2} \neq 0$, and so
$\dim(V_{2}) = 2^{28}a$ for some integer $a \geq 1$.

Assume $V_{i} \neq 0$ for all $i$. Since $V$ is self-dual, all $V_{i}$ are self-dual.
Furthermore, observe that $\Hom_{Z}(V_{2},V_{i}) = 0$ for $i = 0,1$ as $z = -1$ on $V_{2}$
and $z = 1$ on $V_{i}$. Similarly, $\Hom_{Q}(V_{0},V_{1}) = 0$. Thus
$\Hom_{N}(V_{i},V_{j}^{*}) = 0$ whenever $i \neq j$. The proof of
\cite[Prop. 2.3]{MMT} applied to the decomposition $V \da_{N} = \oplus^{2}_{i=0}V_{i}$ shows that
$X(V)$ enters $1_{N} \ua G$ with multiplicity $\geq 2$. Thus
$2(d(d-1)/2 -2) \leq (G:N)$, and so $d < \sqrt{(G:N)} + 1 < 541,379,154$, contrary to
(\ref{e8-2}). On the other hand, $V_{2} \neq 0$, and $V_{1} \neq 0$
according to (the proof of) \cite[Lemma 1]{Ho}. So $V_{0} = 0$.

As mentioned in \cite{Ho}, $L$ acts on $\Irr(Q/Z) \setminus \{1_{Q/Z}\}$ with five orbits.
Among them, four have length larger than $2^{44} > \dim(V)$. So $V_{1}$ has to afford only
the smallest orbit, call it $\OC$, of length $(2^{5}+1)(2^{9}+1)(2^{14}-1)$.
By (\ref{e8-1}), $(\dim(V)-|\OC|)/2^{28} < 1.82$, so $a = 1$ and $\dim(V_{2}) = 2^{28}$.
Next, $(\dim(V)-2^{28})/|\OC| < 1.79$, so $V_{1}$ cannot afford $\OC$ twice, i.e.
$\dim(V_{1}) = |\OC|$. Consequently, $\dim(V) = 2^{28} + |\OC| = 545,783,263$, contradicting
(\ref{e8-2}).
\end{proof}

\begin{propo}\label{2e6}
{\sl Assume $G \rhd S = \ta E_{6}(3)$ or $\tb D_{4}(3)$, and $\ell \neq 3$.

{\rm (i)} Assume in addition that $\ell \neq 2$ if $S = \tb D_{4}(3)$. Then
$G$ has no faithful irreducible $\ell$-modular representation $V$ such that $X(V)$ is
irreducible, if $X = \AC$, resp. $\TSB$, $\TWB$, for $V$ of type $\circ$, resp. $+$, $-$.

{\rm (ii)} Assume $S = \tb D_{4}(3)$. Then $G$ cannot satisfy $\MBA$.}
\end{propo}

\begin{proof}
Assume the contrary. Then we are in the {\it good case} considered in \cite{MMT}.
Let $Z$ be a long-root subgroup of $S$ and $N := N_{G}(Z)$. We write
$V \da_{N} = \oplus^{2}_{i=0}V_{i}$ with $V_{0} := C_{V}(Q)$, $V_{1} := [C_{V}(Z),Q]$, and
$V_{2} := [V,Z]$. By Lemma \ref{root}, $(G:N) = (S:N \cap S)$.

1) First we consider the case $S = \ta E_{6}(3)$. The proof of \cite[Thm. 3.1]{MMT} yields
\begin{equation}\label{e6-1}
  \dim(V) \leq 175,030.
\end{equation}
On the other hand,
\begin{equation}\label{e6-2}
  \dim(V) \geq 172,936
\end{equation}
by \cite[Thm. 4.2]{MMT}. Since $\sqrt{(G:N)+2} < 123,765 < \dim(V)$, the proof of
\cite[Thm. 3.1]{MMT} shows that $V$ is self-dual.

Recall that $N = QL$ with $Q = 3^{1+20}$, $Z = Z(Q)$, and $L = SU_{6}(3) \cdot \ZZ_{2}$. Since $V$ is self-dual, all $V_{i}$ are self-dual. Since $Z \not\leq Z(G)$,
$V_{2} \neq 0$, and $\dim(V_{2}) = 2a \cdot 3^{10}$ for some integer $a \geq 1$.

Assume $V_{i} \neq 0$ for all $i$. Notice that $\Hom_{Z}(V_{2},V_{i}) = 0$ for $i = 0,1$, and
$\Hom_{Q}(V_{0},V_{1}) = 0$. Thus $\Hom_{N}(V_{i},V_{j}^{*}) = 0$ whenever $i \neq j$. The proof
of \cite[Prop. 2.3]{MMT} applied to the decomposition $V \da_{N} = \oplus^{2}_{i=0}V_{i}$ shows
that $X(V)$ enters $1_{N} \ua G$ with multiplicity $\geq 2$. Thus
$2(d(d-1)/2 -2) \leq (G:N)$, and so $d < \sqrt{(G:N)} + 1 < 123,766$, contrary to
(\ref{e6-2}). On the other hand, $V_{2} \neq 0$, and $V_{1} \neq 0$
according to (the proof of) \cite[Lemma 1]{Ho}. So $V_{0} = 0$.

As mentioned in the proof of \cite[Thm. 4.2]{MMT}, $L$ acts on $\Irr(Q/Z) \setminus \{1_{Q/Z}\}$
with five orbits. Among them, four have length larger than $\dim(V)$. So $V_{1}$ has to afford
only the smallest orbit, call it $\OC$, of length $(3^{2}-1)(3^{3}+1)(3^{5}+1)$.
By (\ref{e6-1}), $(\dim(V)-|\OC|)/(2 \cdot 3^{10}) < 1.05$, so $a = 1$ and
$\dim(V_{2}) = 2 \cdot 3^{10}$. Next, $(\dim(V) - 2 \cdot 3^{10})/|\OC| < 1.05$, so $V_{1}$
cannot afford $\OC$ twice, i.e. $\dim(V_{1}) = |\OC|$. Consequently,
$\dim(V) = 2 \cdot 3^{10} + |\OC| = 172,754$, contradicting (\ref{e6-2}).

2) Now we consider the case $S = \tb D_{4}(3)$. Arguing as above and using \cite[Thm. 4.1]{MMT}
instead of \cite[Thm. 4.2]{MMT}, we get $V \simeq V^{*}$, $V_{0} = 0$, $\dim(V_{2}) = 162$,
$\dim(V_{1}) = 56$, whence $\dim(V) = 218$. According to \cite{Lu1}, the last equality cannot
hold if $\ell \neq 2,3$. Assume $\ell = 2$. By \cite{Him}, $V|_{S}$ is contained in the reduction modulo
$2$ of the unique irreducible complex representation of degree $219$ of $S$. For this representation
$V|_{S}$, Hiss (private communication) has shown that $(V \otimes V)|_{S}$ contains an irreducible
constituent of degree $3942$. However, $V \otimes V$ has $\GC$-composition factors $\TWB(V)$ 
(of dimension $23652 = 6 \cdot 3942)$) and $\Vbb$ (of dimension $218$), and $(G : Z(G)S) \leq 3$. So
$G$ cannot satisfy $\MBA$ by Lemma \ref{g-m}.
\end{proof}

\begin{lemma}\label{f43}
{\sl Assume $G \rhd S = F_{4}(3)$ and $\ell \neq 3$. Then $G$ has no faithful
irreducible $\ell$-modular representation $V$ such that $X(V)$ is irreducible for some
$X \in \{\AC,\TSB,\TWB\}$.}
\end{lemma}

\begin{proof}
Assume the contrary. Then we are in the {\it good case} considered in \cite{MMT}. Let $Z$ be a
long-root subgroup of $S$ and $N := N_{G}(Z)$. By Lemma \ref{root}, $(G:N) = (S:N \cap S)$.
The proof of \cite[Thm 3.1]{MMT} yields $\dim(V) \leq 6601$. On the other hand,
$\dim(V) \geq 3^{8}+3^{4}-2 = 6640$ by \cite{MT2}, a contradiction.
\end{proof}

\begin{lemma}\label{2e62}
{\sl Assume $\BS = F_{4}(2)$, resp. $\ta E_{6}(2)$, and $\ell \neq 2$. If $V$ is any irreducible 
$\ell$-modular $\FF G$-representation $V$ such that $X(V)$ is irreducible for some 
$X \in \{\AC,\TSB,\TWB\}$ then $\dim(V) \leq 528$, resp. $2817$. In particular, $G$ cannot satisfy 
$\MBA$ for any faithful irreducible complex representation $V$, unless $S = 2F_{4}(2)$ and 
$\dim(V) = 52$.}
\end{lemma}

\begin{proof}
Assume $X(V)$ is irreducible. Then we are in the {\it good case} considered in \cite{MMT}. Let $Z$ 
be a long-root subgroup of $S$ and $N := N_{G}(Z)$. By Lemma \ref{root}, $(G:N) = (S:N \cap S)$ if 
$\BS = \ta E_{6}(2)$, and $(G:N) \leq 2(S:N \cap S)$ if $\BS = F_{4}(2)$.  
The proof of \cite[Thm 3.1]{MMT} now yields $\dim(V) \leq 2817$, resp. $528$. Now assume $\ell = 0$
and $\BS = \ta E_{6}(2)$. According to \cite{Lu1}, $S = \BS$ or $2\BS$, and $\dim(V) = 1938$ or 
$2432$. Using \cite{Atlas} one can see that $\MBA$ fails. The case $\BS = F_{4}(2)$ is similar.
\end{proof}

Now we can state the main result concerning the finite groups of Lie type:

\begin{theor}\label{Lie}
{\sl Let $G$ be a finite subgroup of $\GC$ that satisfies $\MBA$.  Assume $S = E(G)$ is a
covering group of a finite simple group $\BS$ of Lie type in characteristic $\neq \ell$.
Then one of the following holds.

{\rm (i)} $\BS = PSp_{2n}(q)$, $q = 3,5$, $G = Z(G)S$, and $V \da_{S}$ is a Weil module of
dimension $(q^{n} \pm 1)/2$. If $\ell \neq 2,3$, or if $q = 5$ and $\ell = 3$, then $G$ cannot
satisfy $\MDA$. If $\ell = 2$ then $G$ cannot satisfy $\MCA$.

{\rm (ii)} $\BS = U_{n}(2)$ and $V \da_{S}$ is a Weil module of dimension
$(2^{n} + 2(-1)^{n})/3$ or $(2^{n}-(-1)^{n})/3$. Moreover, $G$ cannot satisfy $\MDA$.

{\rm (iii)} $\BS \in \EC_{3}$ and $V$ is as listed in Tables {\rm IV} and {\rm V} below.

{\rm (iv)} $\BS = F_{4}(2)$ or $\ta E_{6}(2)$. Moreover, if $\ell = 0$ then 
$S = 2F_{4}(2)$ and $\dim(V) = 52$.
\hfill $\Box$}
\end{theor}

\section{Nearly Simple Groups. III}\label{nearly3}
In this section we treat the {\bf small} nearly simple groups, that is, the ones belonging to
$\cup^{3}_{i=1}\EC_{i}$. The outline of the arguments is as follows. For the large groups,
modular character tables of which are not completed yet, we use Lemma \ref{weight}, upper bounds on
complex character degrees of $G$ (such as $\sqrt{|G|}$), and lower bounds for modular
degrees of $G$ (such as given in \cite{LS} and \cite{Jans}) to show that
$\MK$ fails for $G$ (say for $k \geq 6$), see the example of $F_{1}$ below. In several cases
(see the example of $Co_{1}$ below), this can be done only for $\FF G$-modules of dimension
larger than certain bound. In this situation, we can use the results of \cite{HM} to
identify the modules $V$ below the bound and handle these modules individually. We also use
\cite{Lu2} to find the dimensions of certain $\GC$-composition factors in $V^{\otimes k}$,
if $\rank(\GC)$ is not large. For groups of Lie type, sometimes we can also use 
Propositions \ref{longroot} and \ref{bound1}, and argue somewhat as in \S\ref{nearly2}.
Most of the quasi-simple groups arising here have a lot of different extensions of $S := E(G)$
and sometimes isoclinic groups behave completely differently (say in the case of
$S = 2 \cdot F_{4}(2)$). In many cases, Lemmas \ref{change2} and \ref{rare} allow us to cut off a
lot of extensions.

A further useful observation is the following lemma:

\begin{lemma}\label{m6-8}
{\sl Assume $V = \FF^{d}$, $d \geq 6$, and $V \simeq V^{*} \simeq V^{(\ell)}$ as $\FF G$-modules.
Then $\MCA$, $\MDA$, and $\ME$ fail for the $\FF G$-module $V$ if $\ell = 2$. Furthermore,
$\MDA$ and $\ME$ fail if $\ell = 3$.}
\end{lemma}

\begin{proof}
Assume $\ell = 2$. Then the $\GC$-composition factor $L(3\om_{1})$ of $V^{\otimes 3}$
is $V \otimes \Vbb$, so it restricts to $G$ as $V \otimes V^{*}$. In particular, it contains
$1_{G}$, whence $\MCA$ fails. A similar argument applied to $L(5\om_{1})$ shows that $\ME$ fails.
Next, the $\GC$-composition factor $L(2\om_{1}+\om_{2})$ of $V^{\otimes 4}$
is $\Vbb \otimes \TWB(V)$, and it restricts to $G$ as $V \otimes \TWB(V)$. For the corresponding
complex module $\VC$ of $Sp(\VC)$ we have
$\VC \otimes \WB(\VC) \simeq \VC \oplus \wedge^{3}(\VC) \oplus L_{\CC}(\om_{1}+\om_{2})$. It follows
that all $G$-composition factors of $V \otimes \TWB(V)$ have dimension $\leq d(d^{2}-4)/3$, which is
less than $\dim(V \otimes \TWB(V))$. Thus $L(2\om_{1}+\om_{2})$ is reducible over $G$, whence
$\MDA$ fails. Now we assume that $\ell = 3$. The above arguments applied to $L(4\om_{1})$ and
$L(3\om_{1}+\om_{2})$ show that $\MDA$ and $\ME$ fails.
\end{proof}

\medskip
In what follows we give detailed arguments for some of the (nontrivial) cases. Denote
$S = E(G)$ and $\BS = S/Z(S)$. Following Lemma \ref{weight}, we set $e := \lfloor d/2 \rfloor$
and $d_{k} := 2^{k} \begin{pmatrix} e\\k \end{pmatrix}$. Also, let $D$ be the largest complex
irreducible character degree of $G$.

\smallskip
{\bf Case:} $\BS = M = F_{1}$.
According to \cite{Jans}, $d \geq 196882$. Assume $k = 6,7$. Then
$d_{k} \geq 8 \cdot 10^{28} > \sqrt{|\BS|} \geq D$. By Lemma \ref{weight}, $\MK$ fails for $G$.

\smallskip
{\bf Case:} $\BS = Co_{1}$.
Assume $d \geq 170$ and $k = 5,6$. Then $d_{k} > (1.04) \cdot 10^{9}$, meanwhile
$D < (1.03) \cdot 10^{9}$, so $\MK$ fails. Now we assume $d \leq 169$. By \cite{HM}, in this case
$d = 24$, $S = 2\BS$, and $V = \FF_{\ell}^{24}$. As an $S$-module, $V = V^{*} = V^{(\ell)}$ and
$V$ lifts to a complex module $\VC$. Also observe that $G = Z(G)S$ as $\Out(\BS) = 1$. Now Lemmas
\ref{m6-8} and \ref{change1} imply that $\MDA$ and $\ME$ fail when $\ell = 2,3$. Assume $\ell > 3$.
Since $V\da_{S}$ is of type $+$, $\GC = GL(V)$ or $O(V)$. By Lemma \ref{weight},
$\wedge^{k}(V)$ is $\GC$-irreducible for $k = 4,5$. But $\wedge^{k}(\VC)$ is $S$-reducible for
$k = 4,5$, so $\MDA$ and $\ME$ fail here. We conclude that $\ME$ and $\MF$ always fail for $G$. Notice
that $\MC$ holds when $d = 24$ and $\ell = 0$.

\smallskip
{\bf Case:} $\BS = J_{2}$.
Here $D = 448$. By Lemma \ref{weight}, $\MC$ and $\MDA$ fail if $d \geq 18$. If $d \geq 31$, then
$d(d-1)/2 -2 > D$ and so $\MBA$ fails. It remains to check through the cases where $d \leq 30$.
Direct check using \cite{Atlas} and \cite{JLPW} shows that $\MBA$ and $\MC$ fail, except possibly
for $d = 6$. Let $d = 6$. Then $S = 2\BS$ and $V$ lifts to a complex module $\VC$. Also,
$V\da_{S}$ is of type $-$. By Lemma \ref{rare} we may assume that $\GC = Sp(V)$. Observe $\MF$ fails
for $\ell > 5$ as $S^{6}(V)$ has dimension $462 > D$. Assume $\ell = 0$ (or $\ell > 7$). Then
$M_{10}(S,\VC) = 909$. A. Cohen and N. Wallach have kindly performed for us a computation using the 
package LIE to decompose the $\GC$-module $\VC^{\otimes m}$ with $m = 5,6$. It turns out in particular 
that $M_{10}(\GC,\VC) = 909$ and $M_{12}(\GC,\VC) = 9449$,
whence $\ME$ holds for $\ell = 0$. The same is true for any positive characteristic $\ell > 7$, as
this $\ell$ is coprime to $|G|$. Assume $\ell = 7$. Then $L(2\om_{1}+\om_{3})$ (of dimension $158$)
is $G$-reducible, so $\ME$ fails. One checks that $\MC$ and $\MD$ hold here. Assume $\ell = 3$.
Then $\MC$ and $\MD$ hold, but $\ME$ and $\MF$ fail, as $G$ is reducible on
$L(\om_{1}+2\om_{2})$ (of degree $286$) and $L(\om_{1}+\om_{2}+\om_{3})$ (of degree $358$).
Assume $\ell = 5$. Then $\MB$ and $\MC$ hold, but $\MDA$ and $\ME$ fail, as $G$ is reducible on
$L(4\om_{1})$ (of degree $126$) and $L(\om_{2}+\om_{3})$ (of degree $62$). Finally, assume
$\ell = 2$. Then $\MB$ holds, but $\MC$ and $\MDA$ fail, as $G$ is reducible on
$L(\om_{3})$ (of degree $8$) and $L(\om_{1}+\om_{3})$ (of degree $48$).

\smallskip
{\bf Case:} $\BS = \ta E_{6}(2)$.
According to \cite{LS}, $d \geq 1536$, so $d_{5} \geq 7 \cdot 10^{13}$, meanwhile
$D < 3 \cdot 10^{12}$, so $\ME$ fails. We claim that $\MD$ fails as well. Assume the contrary:
$\MD$ holds for some $V$. Since $\MD$ implies $\MB$, $G$ is irreducible and primitive on $V$
by Proposition \ref{irred1}. It is known that $\Out(\BS) = \SSS_{3}$. Observe
that the largest complex degree of $\BS$ and $\BS \cdot 2$ is at most $(1.4) \cdot 10^{11}$.
Next, the largest complex degree of $\BS \cdot 3$ is at most $(1.1) \cdot 10^{11}$. It follows
that the largest complex degree $D'$ of any subgroup of $\Aut(\BS)$ is at most
$(2.2) \cdot 10^{11}$. First we consider the case $S = 6\BS$ or $3\BS$. Since $Z(S) \geq \ZZ_{3}$
acts scalarly and faithfully on $V$, $\GC = GL(V)$. By Lemma \ref{nolan2},
$L(\om_{2} + \om_{d-2}) \hra V^{\otimes 2} \otimes (V^{*})^{\otimes 2}$. The orbit length of
$\om_{2}+\om_{d-2}$ under the Weyl group of $\GC$ is $d(d-1)(d-2)(d-3)/4 > 10^{12} > D'$. Hence
$G$ is reducible on $L(\om_{2} + \om_{d-2})$ (as $Z(G)$ acts trivially on
$V^{\otimes 2} \otimes (V^{*})^{\otimes 2}$), a contradiction. Next assume that
$S = 2\BS$ or $\BS$. Then $Z(S) \leq \ZZ_{2}$ acts trivially on $V^{\otimes 4}$. Now
$d_{4} > (2.3) \cdot 10^{11} > D'$, whence $G$ is reducible on $L(\om_{4})$, again a contradiction.

\smallskip
{\bf Case:} $\BS = F_{4}(2)$.
Here $D < (2.93) \cdot 10^{7}$. If $d \geq 186$, then $d_{4}, d_{5} > (2.94) \cdot 10^{7} > D$
and so $\MDA$ and $\ME$ fail. Assume $d \leq 185$. By \cite{HM}, $d = 52$, $S = 2\BS$,
$V$ lifts to a complex module $\VC$, and $V\da_{S}$ is of type $+$. By Lemmas \ref{rare} and
\ref{change2}, we may assume that $\GC = O(V)$. Here either $G = S$ or $G = S \cdot 2$ (in the
latter case $G$ is isoclinic to the extension $2\BS \cdot 2$ indicated in \cite{Atlas}). One can
check that $\MC$ holds for $S \cdot 2$ if $\ell = 0$, and $\MC$ fails for $S$ for all
$\ell \neq 2$. We show that $\MDA$ and $\ME$ fail. It suffices to prove it for $G = S \cdot 2$.
Observe that $V = V^{*} = V^{(\ell)}$. Hence $\MDA$ and $\ME$ fail for $\ell = 2,3$. Assume
$\ell \geq 5$. By Premet's theorem \cite{Pre}, $L(4\om_{1})$ affords the weights $4\om_{1}$,
$\om_{1}+\om_{3}$, $2\om_{1}+\om_{2}$, whose orbit length under the Weyl group of $\GC = O(V)$ is
$239200$, $62400$, $2600$. Assuming $\MDA$ holds, we see that $\alpha := S^{4}(\VC)/S^{2}(\VC)$
contains a $G$-irreducible constituent say $\beta$ with
$304200 \leq \beta(1) \leq 339677 = \alpha(1)$. Inspecting the character table of $S$ we see that
$\beta(1) = 322218$ or $324870$. Thus $\alpha-\beta$ is a character of $S$ of degree $17459$, resp.
or $14807$. One can also see that the value of $\alpha-\beta$ at an element of class $2A$ of
$S$ is $11315$, resp. $9815$. It is easy to see that $S$ has no such character, a contradiction.
A similar argument, applied to $L(5\om_{1})$ if $\ell > 5$ and $L(3\om_{1}+\om_{2})$ if
$\ell = 5$ (both of them are inside $S^{5}(V)$), shows that $\ME$ fail.

\smallskip
{\bf Case:} $\BS = \Omega^{+}_{8}(2)$.
If $d \geq 54$ then $d_{3}, d_{4} \geq 23400 > D$ and $\MC$, $\MDA$ fail.
Assume $V$ satisfies $\MBA$. By Proposition \ref{irred1}, $S$ is irreducible on $V$, so
$S = \BS$ or $2\BS$. Notice that $\BS$ has an involution $\bar{t}$ of class $2A$ which lifts to
an involution $t \in S$, and $C := C_{G}(t)$ has index $\leq 3150$ in $G$. Since $\ell \neq 2$ and
$t$ affords both eigenvalues $1$ and $-1$ on $V$, we can apply Proposition \ref{bound1} to $C$.
It follows in particular that $d \leq 79$. Checking through the $\FF G$-representations of
dimension $\leq 79$ we see that both $\MBA$ and $\MC$ fail unless $d = 8$. Assume $d = 8$. In
this case $V$ lifts to a complex module $\VC$, $S = 2\BS$, and as an $S$-module,
$V = V^{*} = V^{(\ell)}$. One can check that $\MB$ and $\MC$ hold for this $V$. However,
$\MDA$ and $\ME$ fail. Indeed, $G$ is reducible on $L(4\om_{1})$ if $\ell > 3$ and on
$L(3\om_{1}+\om_{2})$ if $\ell > 5$. Assume $\ell = 3$. Then
$L(4\om_{1})\da_{S} = (V \otimes \Vcc)\da_{S} = (V \otimes V^{*})\da_{S}$ contains $1_{S}$, so
$L(4\om_{1})$ is reducible over $G$ as well. Furthermore, $L(3\om_{1}+\om_{2})$ (of degree $224$)
restricted to $S$ contains an irreducible constituent of degree $8$, whence it is also reducible
over $G$. Assume $\ell = 5$. Here $\GC = GL(V)$ or $O(V)$. Notice that $S^{3}(\VC) = 112a + 8a$ and
$S^{5}(\VC) = 560a + 2 \cdot 112a + 8a$ over $S$ (where $560a$ is some irreducible character of
$S$ or degree $560$, similarly for $112a$ and $8a$). Moreover, $560a$ and $8a$ are irreducible
modulo $5$, and $112a (\mod 5) = 104a + 8a$. If $\GC = GL(V)$ then
$S^{5}(V) = L(5\om_{1}) + L(3\om_{1}+\om_{2})$ as a $\GC$-module, and so
$L(3\om_{1}+\om_{2})\da_{S} = 560a + 2 \cdot 104a + 2 \cdot 8a$ and therefore $L(3\om_{1}+\om_{2})$
is reducible over $G$.  If $\GC = O(V)$ then $S^{5}(V) = S^{3}(V) + L(5\om_{1}) + L(3\om_{1}+\om_{2})$
as a $\GC$-module, and so $L(3\om_{1}+\om_{2})\da_{S} = 560a + 104a$ and therefore
$L(3\om_{1}+\om_{2})$ is reducible over $G$.

\smallskip
{\bf Case:} $\BS = U_{4}(3)$.
If $d \geq 26$ then $d_{4} \geq 11440 > D$ and so $\MDA$ fail. Assume $d \geq 32$, then
$d_{3} \geq 4480$. By Lemma \ref{weight}, $V^{\otimes 3}$ contains a $\GC$-composition factor $W$ of
dimension $\geq d_{3}$. If $\ell = 2$, then some $S$-constituent $U$ of $W$ is acted on by
$3_{1}\BS$ or $3_{2}\BS$, whence $G$-composition factors of $U\ua^{G}$ have dimension
$\leq 3780$ and so $G$ is reducible on $W$. If $\ell > 3$ then we can take $W = X(V)$ with
$X \in \{\wedge^{3},S^{3}\}$. Now $O_{3}(\Mult(\BS))$ acts trivially on some $S$-constituent $U$ of
$X(V)$, whence $W$ is acted on by $4\BS$ and so $G$-composition factors of $U\ua^{G}$
have dimension $\leq 2560$, and again $G$ is reducible on $W$. Thus $\MC$ fails if $d \geq 32$.
Assume $V$ satisfies $\MBA$. By Proposition \ref{irred1}, $S$ is irreducible on $V$, so
$S = \BS$ or $2\BS$. Notice that $\BS$ has an element $\bar{t}$ of class $3A$ which lifts to
an element $t \in S$ of order $3$, and $N := N_{G}(\langle t \rangle)$ has index $280$ in $G$.
If $3$ divides $|Z(S)|$ then (since $Z(S)$ acts scalarly on $V$) $\GC = GL(V)$, so Proposition
\ref{bound1}(i) applies to $N$. If $(|Z(S)|,3) = 1$, then $S$ is a quotient of
$SU_{4}(3)$, in which case by Proposition \ref{longroot}, Proposition \ref{bound1}(i), (ii)
applies to $N$. It follows in particular that $\MBA$ fails if $d \geq 30$.
Checking through the $\FF G$-representations of dimension $\leq 31$ we see that both $\MBA$ and $\MC$
fail unless $d = 6$. Assume $d = 6$. Then $\MC$ holds if $\ell \neq 2,3$, and $\MB$ holds if
$\ell \neq 3$. However, $G$ is reducible on $L(4\om_{1})$ and $L(3\om_{1}+2\om_{5})$ if
$\ell > 3$, and $G$ is reducible on $L(3\om_{1})$ and $L(2\om_{1}+\om_{2})$ if $\ell = 2$.
Thus $\MDA$ and $\ME$ fail.

\medskip
The following tables list results for small nearly simple groups. In the column for
$\MK$, $+$ means $\MK$ holds and $-$ means it does not. In the column for $\MB$, if
$\MB$ holds then we also indicate the values for the corresponding triple
$(d,S,\ell)$. In this case, we find the largest $k$ such that $\MK$ holds for the
$\FF G$-module $V$ of this dimension $d$ and indicate a smallest of the corresponding groups
$G$ in place of $S$. This $G$ may be larger than $S$. For instance, in the third row for
$\AAA_{6}$ in Table III, $\MC$ holds for $G = \HA_{6} \cdot 2_{1}$ but not for $\HA_{6}$.
The notation for outer automorphisms is taken from \cite{Atlas}.

\medskip
Table III lists all $G$-modules $V$, where $S/Z(S) = \AAA_{n}$, $V$ satisfies $\MBA$ if
$n \geq 9$ and $V$ satisfies either $\MBA$ or $\MC$ if $5 \leq n \leq 8$. We also assume
that $d := \dim(V) \geq 4$ (which excludes a few representations of degree $2$, $3$, of 
$2\AAA_{5}$ and $3\AAA_{6}$. See also Theorem \ref{small} in this regard). 
``Natural'' means that $S = \AAA_{n}$, either $d = n-2$ and $\ell|n$
or $d = n-1$ and $\ell = 2$,  and $V \da_{S}$ is the heart of the natural permutation module
(that is, labeled by $(n-1,1)$).

\vskip10pt
\centerline
{{\sc Table III.} Groups $G$ with $S/Z(S) = \AAA_{n}$, $n \geq 5$, in dimension $\geq 4$}
\begin{center}
\begin{tabular}{|c||c|c|c|c|} \hline \skipa
     $S/Z(S)$ & $\MB$ & $\MC$ & $\MD$ & $\ME$ \\ \skipa \hline \hline
     $\AAA_{n}$ & natural & $-$ & $-$ & $-$ \\ \hline \skipa
     $\AAA_{6}$ & $(4,\AAA_{6} \cdot 2_{1},2)$ & $-$ & $+$ & $-$ \\
     $\AAA_{6}$ & $(4,\AAA_{6},3)$ & $-$ & $-$ & $-$ \\
     $\AAA_{6}$ & $(4,\HA_{6} \cdot 2_{1},0)$ & $+$ & $-$ & $-$ \\
     $\AAA_{6}$ & $(4,\HA_{6},5)$ & $-$ & $-$ & $-$ \\ \hline
     $\AAA_{7}$ & $(6,6\AAA_{7},0)$ & $-$ & $-$ & $-$ \\
     $\AAA_{7}$ & $(4,\HA_{7},0 \mbox{ or } 7)$ & $+$ & $-$ & $-$ \\
     $\AAA_{7}$ & $(4,\HA_{7},3 \mbox{ or } 5)$ & $-$ & $-$ & $-$ \\
     $\AAA_{7}$ & $(4,\AAA_{7},2)$ & $-$ & $-$ & $-$ \\ \hline
     $\AAA_{8}$ & $(4,\AAA_{8},2)$ & $-$ & $-$ & $-$ \\  \hline
     $\AAA_{9}$ & $(8,\AAA_{9},2)$ & $-$ & $-$ & $-$ \\
     $\AAA_{9}$ & $(8,\HA_{9},0)$ & $+$ & $-$ & $-$ \\
     $\AAA_{9}$ & $(8,\HA_{9},5 \mbox{ or } 7)$ & $-$ & $-$ & $-$ \\ \hline
     $\AAA_{10}$ & $(8,\HA_{10},5)$ & $-$ & $-$ & $-$ \\ \hline
\end{tabular}
\end{center}

\newpage
\vskip10pt
{\small
\centerline
{{\sc Table} IV. Groups $G$ with $S/Z(S)$ being a sporadic finite simple group:
Modular case}
\vspace{0.1cm}
\begin{center}
\begin{tabular}{|c||c|c|c|c|c|c|} \hline
     $S/Z(S)$ & $\MB$ & $\MC$ & $\MD$ & $\ME$ & $\MF$ & $\MG$ \\ \skipa \hline \hline
     $M_{11}$ & $(10,M_{11},2)$ & $-$ & $-$ & $-$ & $-$ & $-$ \\ \hline
     $M_{12}$ & $(10,M_{12}, 2 \mbox{ or }3)$ & $-$ & $-$ & $-$ & $-$ & $-$ \\
     $M_{12}$ & $(10,2M_{12}, \neq 2,3)$ & $-$ & $-$ & $-$ & $-$ & $-$ \\ \hline
     $J_{1}$ & $-$ & $-$ & $-$ & $-$ & $-$ & $-$ \\ \hline
     $M_{22}$ & $(6,3M_{22}, 2)$ & $-$ & $-$ & $-$ & $-$ & $-$ \\
     $M_{22}$ & $(10,2M_{22}, 7)$ & $+$ & $-$ & $-$ & $-$ & $-$ \\
     $M_{22}$ & $(10,2M_{22}, \neq 2,7)$ & $-$ & $-$ & $-$ & $-$ & $-$ \\ \hline
     $J_{2}$ & $(6,2J_{2}, 0)$ & $+$ & $+$ & $+$ & $-$ & $-$ \\
     $J_{2}$ & $(6,2J_{2}, 3 \mbox{ or }7)$ & $+$ & $+$ & $-$ & $-$ & $-$ \\
     $J_{2}$ & $(6,2J_{2}, 5)$ & $+$ & $-$ & $-$ & $-$ & $-$ \\
     $J_{2}$ & $(6,J_{2}, 2)$ & $-$ & $-$ & $-$ & $-$ & $-$ \\ \hline
     $M_{23}$ & $(21,M_{23}, 23)$ & $-$ & $-$ & $-$ & $-$ & $-$ \\
     $M_{23}$ & $(45,M_{23}, 0 \mbox{ or 7})$ & $-$ & $-$ & $-$ & $-$ & $-$ \\ \hline
     $HS$ & $-$ & $-$ & $-$ & $-$ & $-$ & $-$ \\ \hline
     $J_{3}$ & $(18,3J_{3}, 0 \mbox{ or }5)$ & $+$ & $-$ & $-$ & $-$ & $-$ \\
     $J_{3}$ & $(18,3J_{3}, 17 \mbox{ or } 19)$ & $-$ & $-$ & $-$ & $-$ & $-$ \\ \hline
     $M_{24}$ & $(22,M_{24},3)$ & $-$ & $-$ & $-$ & $-$ & $-$ \\
     $M_{24}$ & $(45,M_{24}, 0, 7, \mbox{ or } 23)$ & $-$ & $-$ & $-$ & $-$ & $-$ \\ \hline
     $McL$ & $(21,McL, 5)$ & $-$ & $-$ & $-$ & $-$ & $-$ \\
     $McL$ & $(22,McL, \neq 3,5)$ & $-$ & $-$ & $-$ & $-$ & $-$ \\ \hline
     $He$ & $-$ & $-$ & $-$ & $-$ & $-$ & $-$ \\ \hline
     $Ru$ & $(28,2Ru, \neq 5)$ & $-$ & $-$ & $-$ & $-$ & $-$ \\ \hline
     $Suz$ & $(12,6Suz, \neq 2,3)$ & $+$ & $-$ & $-$ & $-$ & $-$ \\
     $Suz$ & $(12,2Suz, 3)$ & $+$ & $-$ & $-$ & $-$ & $-$ \\
     $Suz$ & $(12,Suz, 2)$ & $-$ & $-$ & $-$ & $-$ & $-$ \\ \hline
     $O'N$ & $(342,3O'N, 0 \mbox{ or } 19)$ & $-$ & $-$ & $-$ & $-$ & $-$ \\
       \hline
     $Co_{3}$ & $(22,Co_{3}, 2)$ & $-$ & $-$ & $-$ & $-$ & $-$ \\
     $Co_{3}$ & $(23,Co_{3}, \neq 2,3)$ & $-$ & $-$ & $-$ & $-$ & $-$ \\ \hline
     $Co_{2}$ & $(22,Co_{2}, 2)$ & $-$ & $-$ & $-$ & $-$ & $-$ \\
     $Co_{2}$ & $(23,Co_{2},\neq 2)$  & $+$ & $-$ & $-$ & $-$ & $-$ \\ \hline
     $Fi_{22}$$^{(\sharp)}$ &  &  & $-$ & $-$ & $-$ & $-$ \\ \hline
     $HN$$^{(\sharp)}$ &  &  & $-$ & $-$ & $-$ & $-$ \\ \hline
     $Ly$$^{(\sharp)}$ &  &  &  & $-$ & $-$ & $-$ \\ \hline
     $Th$$^{(\sharp)}$ &  &  &  & $-$ & $-$ & $-$ \\ \hline
     $Fi_{23}$$^{(\sharp)}$ &  &  &  & $-$ & $-$ & $-$ \\ \hline
     $Co_{1}$$^{(\sharp)}$ &  &  &  & $-$ & $-$ & $-$ \\ \hline
     $J_{4}$$^{(\sharp)}$ &  &  &  & $-$ & $-$ & $-$ \\ \hline
     $Fi'_{24}$$^{(\sharp)}$ &  &  &  & $-$ & $-$ & $-$ \\ \hline
     $BM$$^{(\sharp)}$ &  &  &  &  & $-$ & $-$ \\ \hline
     $M$$^{(\sharp)}$ &  &  &  &  & $-$ & $-$ \\ \hline
\end{tabular}
\end{center}}

\newpage
\vskip10pt
\centerline
{{\sc Table} V. Groups $G$ with $S/Z(S)$ being a ``small'' group of Lie type:}
\centerline{Cross-characteristic case in dimension $> 4$}
\begin{center}
\begin{tabular}{|c||c|c|c|c|} \hline
     $S/Z(S)$ & $\MB$ & $\MC$ & $\MD$ & $\ME$ \\ \skipa \hline \hline
     $L_{2}(7) \simeq L_{3}(2)$ & $-$ & $-$ & $-$ & $-$ \\ \hline
     $L_{2}(11)$ & $-$ & $-$ & $-$ & $-$ \\ \hline
     $L_{2}(13)$ & $(6,L_{2}(13),2)$ & $-$ & $-$ & $-$ \\ \hline
     $L_{3}(4)$ & $(6,6 \cdot L_{3}(4) \cdot 2_{1}, 0)$ & $+$ & $-$ & $-$ \\
     $L_{3}(4)$ & $(6,2 \cdot L_{3}(4) \cdot 2_{1}, 3)$ & $+$ & $-$ & $-$ \\
     $L_{3}(4)$ & $(6,6 \cdot L_{3}(4), 5 \mbox{ or }7)$ & $-$ & $-$ & $-$ \\
     $L_{3}(4)$ & $(8,4_{1} \cdot L_{3}(4), 0 \mbox{ or }5)$ & $-$ & $-$ & $-$ \\ \hline
     $U_{3}(3)$ & $(6,U_{3}(3) \cdot 2, 0)$ & $+$ & $-$ & $-$ \\
     $U_{3}(3)$ & $(6,U_{3}(3), 2 \mbox{ or }7)$ & $-$ & $-$ & $-$ \\ \hline
     $U_{4}(2) \simeq PSp_{4}(3)$ & $(5,U_{4}(2),3)$ & $+$ & $-$ & $-$ \\
     $U_{4}(2) \simeq PSp_{4}(3)$ & $(5,U_{4}(2),0 \mbox{ or }5)$ & $-$ & $-$ & $-$ \\
     $U_{4}(2) \simeq PSp_{4}(3)$ & $(6,U_{4}(2),\neq 3)$ & $-$ & $-$ & $-$ \\ \hline
     $U_{4}(3)$ & $(6,6_{1} \cdot U_{4}(3), \neq 2,3)$ & $+$ & $-$ & $-$ \\
     $U_{4}(3)$ & $(6,3_{1} \cdot U_{4}(3), 2)$ & $-$ & $-$ & $-$ \\ \hline
     $U_{5}(2)$ & $(10,(2 \times U_{5}(2)) \cdot 2, \neq 2,3)$ & $+$ & $-$ & $-$ \\
     $U_{5}(2)$ & $(10,U_{5}(2),3)$ & $-$ & $-$ & $-$ \\
     $U_{5}(2)$ & $(11,U_{5}(2),\neq 2,3)$ & $-$ & $-$ & $-$ \\ \hline
     $U_{6}(2)$ & $(21,3 \cdot U_{6}(2), \neq 2,3)$ & $-$ & $-$ & $-$ \\
     $U_{6}(2)$ & $(21,U_{6}(2),3)$ & $-$ & $-$ & $-$ \\
     $U_{6}(2)$ & $(22,U_{6}(2), \neq 2,3)$ & $-$ & $-$ & $-$ \\ \hline
     $Sp_{4}(4)$ & $(18,(2 \times Sp_{4}(4)) \cdot 4, \neq 2,3)$ & $-$ & $-$ & $-$ \\ \hline
     $Sp_{6}(2)$ & $(7,Sp_{6}(2), \neq 2)$ & $-$ & $-$ & $-$ \\ \hline
     $\Omega^{+}_{8}(2)$ & $(8,2 \cdot \Omega^{+}_{8}(2),\neq 2)$ & $+$ & $-$ & $-$ \\ \hline
     $\Omega_{7}(3)$ & $-$ & $-$ & $-$ & $-$ \\ \hline
     $\ta B_{2}(8)$ & $(14,\ta B_{2}(8) \cdot 3, 0)$ & $-$ & $-$ & $-$ \\ \hline
     $\tb D_{4}(2)$ & $-$ & $-$ & $-$ & $-$ \\ \hline
     $G_{2}(3)$ & $(14,G_{2}(3), \neq 3)$ & $-$ & $-$ & $-$ \\ \hline
     $G_{2}(4)$ & $(12,2 \cdot G_{2}(4) \cdot 2, \neq 2,5)$ & $+$ & $-$ & $-$ \\
     $G_{2}(4)$ & $(12,2 \cdot G_{2}(4),5)$ & $-$ & $-$ & $-$ \\ \hline
     $\ta F_{4}(2)'$ & $(26,\ta F_{4}(2)', 0 \mbox{ or }13)$ & $-$ & $-$ & $-$ \\ \hline
     $F_{4}(2)$$^{(\sharp)}$ & $(52,2 \cdot F_{4}(2) \cdot 2, \neq 2,3)$ & $+$ & $-$ & $-$ \\
     $F_{4}(2)$$^{(\sharp)}$ & $(52,2 \cdot F_{4}(2), 3)$ & $-$ &  $-$ & $-$ \\ \hline
     $\ta E_{6}(2)$$^{(\sharp)}$ &  &  &  $-$ & $-$ \\ \hline
\end{tabular}
\end{center}

\medskip
Table IV lists results in the modular case of sporadic groups. 
For $10$ large sporadic groups ($Fi_{22}$, $HN$, $Ly$, $Th$,
$Fi_{23}$, $Co_{1}$, $J_{4}$, $Fi'_{24}$, $BM$, $M$), we only prove that for some $k \leq 6$,
$\MK$ can never hold for $\FF G$-modules in the modular case; these cases are marked with 
$^{(\sharp)}$. For the $16$ smaller
sporadic groups, we list all modules which satisfy either $\MBA$ or $\MC$. (Notice that
{\it complex} modules for sporadic groups that satisfy $\MB$ have been classified in Theorem 
\ref{complex}.)

Table V list all modules of dimension $> 4$ which satisfy either $\MBA$ or $\MC$ for
small groups of Lie type (except for the cases of $F_{4}(2)$ and $\ta E_{6}(2)$, where we only 
show that $\MD$ and $\ME$ both fail; these two cases are distinguished by $^{(\sharp)}$).
Notice that the cases $S/Z(S) \simeq \AAA_{n}$ with
$n = 5,6,8$ are already listed in Table III. 

In Table VI we list all groups $G$ (with no composition factor belonging to $\EC_{1} \cup \EC_{2}$) 
that satisfy $\MBA$ but not $\MB$. In all these cases,
$\GC = GL(V)$, and it turns out that $\MC$ and $\MDA$ fail.

\vskip10pt
\centerline
{{\sc Table} VI. Examples where $\MBA$, but not $\MB$, holds}
\begin{center}
\begin{tabular}{|c||c|c|c|} \hline
     $S/Z(S)$ & $(d,G,\ell)$ & $\MC$ & $\MDA$ \\ \skipa \hline \hline
     $\AAA_{7}$ & $(6,6 \cdot \AAA_{7}, 5 \mbox{ or } 7)$ & $-$ & $-$\\ \hline
     $L_{3}(4)$ & $(8,4_{1} \cdot L_{3}(4), 3 \mbox{ or } 7)$ & $-$ & $-$\\ \hline
     $\ta B_{2}(8)$ & $(14,\ta B_{2}(8) \cdot 3, 5 \mbox{ or } 13)$ & $-$ & $-$\\ \hline
     $\ta F_{4}(2)'$ & $(26,\ta F_{4}(2)', 3 \mbox{ or } 5)$ & $-$ & $-$\\ \hline
     $M_{11}$ & $(5,M_{11}, 3)$ & $-$ & $-$\\ \hline
     $M_{12}$ & $(6,2 \cdot M_{12}, 3)$ & $-$ & $-$\\
     $M_{12}$ & $(10,2 \cdot M_{12}, 3)$ & $-$ & $-$\\ \hline
     $M_{22}$ & $(10,M_{22}, 2)$ & $-$ & $-$\\ \hline
     $M_{23}$ & $(11,M_{23}, 2)$ & $-$ & $-$\\
     $M_{23}$ & $(45,M_{23}, 11 \mbox{ or } 23)$ & $-$ & $-$\\ \hline
     $J_{3}$ & $(9,3 \cdot J_{3}, 2)$ & $-$ & $-$\\
     $J_{3}$ & $(18,3 \cdot J_{3}, 2)$ & $-$ & $-$\\ \hline
     $M_{24}$ & $(11,M_{24}, 2)$ & $-$ & $-$\\
     $M_{24}$ & $(45,M_{24}, 11)$ & $-$ & $-$\\ \hline
     $Ru$ & $(28,2Ru, 5)$ & $-$ & $-$\\ \hline
     $O'N$ & $(342,3O'N, 5,7,11, \mbox{ or }31)$ & $-$ & $-$\\ \hline
\end{tabular}
\end{center}

\section{Proofs of Main Theorems}

\begin{lemma}\label{sp-o}
{\sl Let $\GC = Sp(V)$ with $V = \FF^{2n}$, $\ell = 2$. Assume
$\GC \geq G \geq \Om^{\pm}_{2n}(q)$ for some even $q$. Then $G$ is irreducible on every
$\GC$-composition factor of $V^{\otimes k}$ with $k \leq \min\{q,n-1\}$. Moreover $G$
satisfies $M_{2k}(\GC)$ for any $k \leq \min\{q-1,n-1\}$.}
\end{lemma}

\begin{proof}
Let $k \leq \min\{q,n-1\}$. First we show that $H := \Om^{\pm}_{2n}(q)$ is irreducible
on any $\GC$-composition factor $W = L(\om)$ of $V^{\otimes k}$.

1) Clearly, $\om = \sum^{n}_{i=1}z_{i}\om_{i}$ is a dominant weight occuring in
$V^{\otimes k}$, so $\om = k\om_{1} - \sum^{n}_{i=1}x_{i}\al_{i}$ for some nonnegative
integers $z_{i},x_{i}$. Here $\al_{1}, \ldots ,\al_{n}$ are the simple roots, and
$\om_{1}, \ldots ,\om_{n}$ are the corresponding fundamental weights. In particular,
$\al_{1} = 2\om_{1} - \om_{2}$, $\al_{i} = -\om_{i-1}+2\om_{i}-\om_{i+1}$ for
$2 \leq i \leq n-1$, and $\al_{n} = -2\om_{n-1}+2\om_{n}$. One checks that
\setcounter{equation}{0}
\begin{equation}\label{sp-o1}
  x_{i-1} - x_{i} = x_{n} + \sum^{n}_{j=i}z_{j}
\end{equation}
for $2 \leq i \leq n$, and
\begin{equation}\label{sp-o2}
  2x_{1} - x_{2} = k-z_{1}.
\end{equation}

We claim that $tz_{t} \leq k$ for $1 \leq t \leq n$. Indeed, (\ref{sp-o1}) implies that
$x_{t-1} \geq z_{t}$, $x_{t-2} \geq 2z_{t}$, $\ldots$, and $x_{1} \geq (t-1)z_{t}$. Now
(\ref{sp-o2}) yields $k \geq x_{1} + (x_{1}-x_{2}) \geq (t-1)z_{t} + z_{t} = tz_{t}$,
as required. In particular, $t \leq k$ if $z_{t} > 0$.

We also claim that either $\om = q\om_{1}$ and $k = q$, or $z_{i} \leq q-1$ for all $i$. Indeed,
if $i \geq 2$ then $2z_{i} \leq iz_{i} \leq k \leq q$, whence $z_{i} < q$. Assume $z_{1} \geq q$.
Since $z_{1} \leq k \leq q$, we must have $k = q = z_{1}$. In this case, (\ref{sp-o2})
implies $0 = x_{1} = x_{1}-x_{2}$, whence $x_{i-1}-x_{i} = 0$ for all $i \geq 2$ by
(\ref{sp-o1}), i.e. $x_{i} = 0$ for all $i$ and so $\om = q\om_{1}$.

2) Now we decompose $\om = \sum^{\infty}_{i=0}2^{i}\lam_{i}$, where all $\lam_{i}$ are restricted
weights. Since $k < n$, the results of 1) imply that $\lam_{i}$ does not involve $\om_{n}$.
Thus the simple $\GC$-module $L(\lam_{i})$ is a restricted irreducible
representation of $\HC := \Omega(V)$ by the results of \cite{S2}. Notice that
$\om_{i} = \pi_{i}$ for $i \leq n-2$ and $\om_{n-1} = \pi_{n-1} + \pi_{n}$, where
$\pi_{1}, \ldots ,\pi_{n}$ are the fundamental weights of $\HC$.

If $k = q$ and $\om = q\om_{1}$, then $L(\om) = V^{(q)}$ is irreducible over $H$. Otherwise
$z_{i} \leq q-1$ for all $i$ by the results of 1). In this case,
$\om = \sum^{f-1}_{i=0}2^{i}\lam_{i}$ if $q = 2^{f}$, and so
$L(\om) \da_{H} \simeq (\otimes^{f-1}_{i=0}L(\lam_{i})^{(2^{i})} \da_{\HC})\da_{H}$ is
irreducible over $H$ by \cite[Thm. 5.4.1]{KL}.

3) Finally, if we take $k \leq \min\{q-1,n-1\}$, then for any dominant weight $\om$ occuring
in $V^{\otimes k}$ we have $\om = \sum^{n}_{i=1}z_{i}\om_{i}$ with $0 \leq z_{i} < q$ and
$z_{n} = 0$. In terms of fundamental weights of $\HC$, we have
$\om = \sum^{n-2}_{i=1}z_{i}\pi_{i} + z_{n-1}(\pi_{n-1} + \pi_{n})$. So if $L(\om)$ and
$L(\om')$ are two distinct $\GC$-composition factors of $V^{\otimes k}$, then they are also
nonisomorphic over $\HC$ and over $H$, meaning $H$ satisfies $M_{2k}(\GC)$.
\end{proof}

Assume $\ell = 2$. Lemma \ref{sp-o} shows that only irreducible modules $L(\om)$ of $Sp(V)$
with $\om$ involving the last fundamental weight $\om_{n}$ can detect the subgroups of $O(V)$
inside $Sp(V)$. Notice that the spin module $L(\om_{n})$ is reducible over $\Omega(V)$ but irreducible
over $O(V)$. The next lemma gives examples of such $L(\om)$. Another way to detect subgroups of
$O(V)$ inside $Sp(V)$ is to view $Sp_{2n}(\FF)$ as $SO_{2n+1}(\FF)$ and use the indecomposable
module $\FF^{2n+1}$ (observe that this module is semisimple over $O(V)$); see Lemma \ref{indecomp}
below.

\begin{lemma}\label{sp-o3}
{\sl Let $\ell = 2$, $V = \FF^{d}$ with $d = 2n \geq 6$. Then $O(V)$ is reducible on both
$Sp(V)$-modules $L(\om_{1}+\om_{n})$ and $L(\sum^{n}_{i=1}\om_{i})$.}
\end{lemma}

\begin{proof}
Let $\pi_{1}, \ldots ,\pi_{n}$ be the fundamental weights of the subgroup $\HC = \Omega(V)$ of
$Sp(V)$. First we consider $U := L(\om_{1} + \om_{n}) \simeq L(\om_{1}) \otimes L(\om_{n})$, of
dimension $2^{n+1}n$. Then $L(\om_{1})|_{\HC} = L(\pi_{1})$ and
$L(\om_{n})|_{\HC} = L(\pi_{n-1}) \oplus L(\pi_{n})$, whence $U|_{\HC}$ contains a subquotient
$L(\pi_{1}+\pi_{n})$. Direct calculation using Weyl's character formula (cf. \cite[p. 410]{FH})
shows that the dimension of the Weyl module of $\HC$ with highest weight $\pi_{1}+\pi_{n}$ is
$2^{n-1}(2n-1)$. In particular, $\dim(L(\pi_{1}+\pi_{n})) < \dim(U)/2$, so $U$ is reducible on
$O(V) = \HC.2$.

Next we consider
$W := L(\sum^{n}_{i=1}\om_{i}) \simeq L(\sum^{n-1}_{i=1}\om_{i}) \otimes L(\om_{n})$. Then
$L(\sum^{n-1}_{i=1}\om_{i})|_{\HC} = L(\sum^{n}_{i=1}\pi_{i})$ has dimension $2^{n(n-1)}$.
As above $L(\om_{n})|_{\HC} = L(\pi_{n-1}) \oplus L(\pi_{n})$, so $W|_{\HC}$ contains a subquotient
$A := L(\pi_{n}+\sum^{n}_{i=1}\pi_{i}) \simeq L(\gamma) \otimes L(\pi_{n-1})^{(2)}$ with
$\gamma := \sum^{n-2}_{i=1}\pi_{i} + \pi_{n}$. Let $f(n)$ be the dimension of the Weyl module
with highest weight $\gamma$. Again using Weyl's character formula we see that $f(3) = 20$ and
$f(n+1)/f(n) = 10\prod^{n-1}_{k=2}\frac{16k^{2}-1}{4k^{2}-1}$ for $n \geq 3$. Observe that
$$\prod^{n-1}_{k=2}\frac{16k^{2}-1}{16k^{2}-4} < \prod^{\infty}_{k=2}\frac{16k^{2}-1}{16k^{2}-4}
  < \prod^{\infty}_{k=2}\left(1+\frac{1}{5k^{2}}\right) <
  \exp\left(\sum^{\infty}_{k=2}\frac{1}{5k^{2}}\right)
  = \exp\left(\frac{\pi^{2}}{30} - \frac{1}{5}\right) < 1.14,$$
whence $f(n+1)/f(n) < (11.4) \cdot 4^{n-2} < 4^{n}$. From this it easy follows that
$f(n) < 2^{n(n-1)}/7$ for all $n \geq 5$. Now for $n \geq 5$ we have
$\dim(A) < 2^{n-1} \cdot 2^{n(n-1)}/7 < \dim(W)/2$, so $W$ is reducible on $O(V)$. The cases
$n = 3,4$ also follow easily.
\end{proof}

\medskip
{\bf Proofs of Theorems \ref{complex} and \ref{complex-m6}.}
As we have already mentioned, Theorem \ref{complex} is a direct consequence of Theorems 
\ref{closed} and \ref{main3}. In turn, Theorem \ref{complex-m6} follows from Theorem \ref{complex}
and Propositions \ref{ex-m62}, \ref{sp-m6}, \ref{su-m6}. 

\medskip
{\bf Proofs of Theorems \ref{main2} and \ref{main3}.}
Assume $G$ is a closed subgroup of $\GC$ that satisfies the assumptions of Theorem \ref{main2},
resp. \ref{main3}. First we assume that $G$ is of positive dimension. Then the results of
\S\ref{red-thms} show that either case (A) of Theorem \ref{main2} holds, or $\GC = Sp(V)$,
$d = 6$, $\ell = 2$ and $G = G_{2}(\FF)$. In the latter case, observe that the $\GC$-module
$L(\om_{3})$ is reducible on $G$, and so is
$L(\om_{1}+\om_{3}) \simeq L(\om_{1}) \otimes L(\om_{3})$, whence $\MDA$ fails. The same happens
to the subgroup $O(V)$ of $Sp(V)$ if $\ell = 2$ and $d = 6$: over $\Omega(V)$ (and also over 
$\Omega^{\pm}_{6}(q)$) the $Sp(V)$-module $L(\om_{1}+\om_{3})$ has two composition factors of 
degree $10$ and two of degree $4$. Also, if $d = 8$ then the spin module $L(\om_{4})$ of $Sp(V)$ is 
irreducible over $O(V)$ but reducible over $\Omega(V)$. Next we may assume that $G$ is finite. The 
results of \S\ref{red-thms} reduce the problems to the cases considered in \S\S \ref{extra}, 
\ref{nearly1}, \ref{nearly2}, \ref{nearly3}, whence one of the possibilities formulated in Theorem 
\ref{main2}, resp. \ref{main3}, has to occur. It remains to verify which possibilities do indeed 
happen.

It is easy to see that the case (A) of Theorem \ref{main2} occurs.
Next we show that the cases (i) -- (v) of Theorem \ref{main2}(B) do occur if we
take $q \geq 4$. Let $n = \rank(\GC)$. Then all $\GC$-composition factors of $V \otimes V^{*}$
and $V^{\otimes k}$ with $k \leq 4$ have highest weights $k\om_{1}$ or $\sum^{n}_{i=1}a_{i}\om_{i}$
with $0 \leq a_{i} \leq 2$. Hence all of them are irreducible over $G$ provided $q \geq 4$, cf.
\cite[Thm. 5.4.1]{KL} (notice that $L(4\om_{1}) \simeq \Vdd$ when $q = 4$). Similarly,
$G$ is irreducible on all $\GC$-composition factors of $V \otimes V^{*}$ for any $q$.

Now we show that $G$ is reducible on $L((q+1)\om_{1})$ if $q = 2,3$. If $S = SU_{d}(q)$,
$Sp_{d}(q)$ or $\Omega^{\pm}_{d}(q)$, then $V^{(q)} \simeq V^{*}$, whence
$L((q+1)\om_{1})\da_{S}$ contains $1_{S}$ and so $L((q+1)\om_{1})$ is reducible over $G$.
If $q = 3$ and $S = SL_{d}(q)$, then $V^{(q)} \simeq V$, whence
$L((q+1)\om_{1})\da_{S}$ is the direct sum of two irreducible $S$-modules
$\TWB(V)$ and $\TSB(V)$ of distinct dimensions, whence $L((q+1)\om_{1})$ is reducible over $G$.
If $q = 2$, $d \geq 4$, and $S = SL_{d}(q)$, then $V^{(q)} \simeq V$, whence
$L((q+1)\om_{1})\da_{S}$ contains irreducible $S$-modules
$\TW(V)$ and $V$ of distinct dimensions, whence $L((q+1)\om_{1})$ is reducible over $G$.
(Notice the case of $SL_{3}(2)$ is not included as we assume $d > 3$.)
Assume $S = \tb D_{4}(q)$. Then $L((q+1)\om_{1})$ has dimension $64$. Notice
that $|\Out(S)| = 3$, and irreducible $q$-Brauer characters of $S$ all have
degree $\neq 64$ as they come from {\it restricted} $\GC$-modules. Hence $G$ is reducible on
$L((q+1)\om_{1})$.

Finally, we have mentioned in \S\S \ref{extra}, \ref{nearly1}, \ref{nearly2}, and \ref{nearly3},
that the cases (C) -- (F) of Theorem \ref{main3} give rise to examples.

\medskip
{\bf Proof of Theorem \ref{main1}.}
Assume $G < \GC$, $\GCR$ is reductive, and $M_{8}(G,V) = M_{8}(\GC,V)$. Then 
$M_{2k}(G,V) = M_{2k}(\GC,V)$ for $k \leq 4$ by Lemma \ref{m68}. Since $G$ acts completely 
reducibly on finite dimensional $\GC$-modules, we see that $G$ satisfies $\MB$, $\MC$, and $\MD$.
Thus $G$ satisfies the assumptions of Theorem \ref{main2}. By Theorem \ref{closed} we may assume
that $G$ is finite. Notice that $\ell = 0$, so case (B) of Theorem \ref{main2} cannot occur.
Case (D) cannot occur either, see Table I and Lemma \ref{2e62}. Case (C) gives rise to an example 
(and in this case even $\ME$ holds true).

\medskip
{\bf Proof of Theorem \ref{main4}.}

1) Define $\mu(g)$ to be the largest possible dimension of $g$-eigenspaces on $V$. We need to show
that $\mu(g)/d \leq 7/8$ with equality occuring only for $G = 2 \cdot O^{+}_{8}(2) < O_{8}(\CC)$. 
Consider any $g \in G \setminus Z(\GC)$. Then $\mu(g) \leq d-1$; in particular we are done if 
$d \leq 7$. From now on we may assume $d \geq 8$; in particular, 
$G \cap Z(\GC) = Z(G)$ by Proposition \ref{irred1}. We can now apply Theorem \ref{complex} to $G$. 
Another reduction is provided by using the function $\ags(x)$ defined in \cite{GS}: for a finite 
simple group $\BS$ and $x \in \Aut(\BS)$, let $\ags(x)$ be the smallest number of 
$\Aut(\BS)$-conjugates of $x$ which generate the subgroup $\langle \BS,x \rangle$. Assuming $G$ is 
nearly simple and $\BS := S/Z(S)$ as in Theorem \ref{complex}, one can show that 
$\mu(g)/d \leq 1 - 1/\ags(g)$, cf. \cite[Lemma 3.2]{GT2}. For all the groups listed in Table I, it 
turns out that $\ags := \max\{\ags(x) \mid x \in G \setminus Z(G)\} \leq 7$, except for the three
cases $\BS \in \{ \Omega^{+}_{8}(2),\AAA_{9},F_{4}(2)\}$. In the first case we
have $\mu(g)/d = 7/8$, $G = 2 \cdot O^{+}_{8}(2) < O_{8}(\CC)$ (and 
$g$ is an involution outside of $[G,G] = 2 \cdot \Omega^{+}_{8}(2)$). In the last two cases direct 
computation using \cite{Atlas} shows that $\mu(g)/d < 7/8$. Thus we need to consider 
only the cases (B) and (C) of Theorem \ref{complex}. Let $\chi$ denote the character of the 
$G$-module $V$. Representing the action of $g$ on $V$ by a diagonal matrix and computing the trace, 
we see that $|\chi(g)| \geq \mu(g) - (d-\mu(g))$, with equality occuring only when 
$g$ has only two distinct eigenvalues $\lam$ and $-\lam$ on $V$. Thus, 
\begin{equation}\label{4mu}
  \mu(g)/d \leq (|\chi(g)|/d + 1)/2\,, 
\end{equation}
with equality occuring only when $g$ has only two distinct eigenvalues $\lam$ and $-\lam$ on $V$.
In what follows we will also assume that $\ags \geq 8$.

2) Here we consider the case (C) of Theorem \ref{complex}.  By \cite[Lemma 2.4]{GT1}, 
$|\chi(g)| \leq p^{a-1/2} = d/\sqrt{p}$. It follows by (\ref{4mu}) that 
$\mu(g)/d \leq (1/\sqrt{2} + 1)/2 < 7/8$.

3) Next we turn to the case (B) of Theorem \ref{complex} and assume $\BS = PSp_{2n}(q)$ with 
$q = 3,5$. With no loss we may assume $G = S$. Suppose that $g$ is not a $2$-element. Then by 
\cite[Lemma 4.1]{GT2}, $\mu(g) \leq (q^{n}+q^{n-1})/4 + \gamma$ and $d = (q^{n}-1)/2 + \gamma$,
with $\gamma \in \{0,1\}$. Since $\ags \geq 8$, $n \geq 4$ by \cite{GS}. Hence (\ref{4mu}) implies 
that $\mu(g)/d < 0.7$, and this last inequality holds for any non-$2$-element in $G$. Next assume 
that $g$ is a $2$-element. Then by \cite[Lemma 3.4]{GT2}, $\mu(g)/d \leq (1+0.7)/2 < 7/8$.

4) Finally, we turn to the case (C) of Theorem \ref{complex} and assume $\BS = U_{n}(2)$. 
Since $\ags \geq 8$, $n \geq 8$ by \cite{GS}. Observe that the Weil modules of $S$ extend to 
$GU_{n}(2)$. First we consider any element $h \in GU_{n}(2) \setminus Z(GU_{n}(2))$. 
Consider the natural module $W := \FF_{4}^{n}$ for $GU_{n}(2)$ and let 
$e(h) := \max \{ \dim(\Ker((h|_{W}-\delta^{i}))) \mid i = 0,1,2\}$, where $\delta$ is a primitive
element of $\FF_{4}$. Assume $e(h) \leq n-2$. Then using \cite[Lemma 4.1]{TZ2} one can show that 
$|\chi(h)| \leq (2^{n-2}+5)/3$. Since $d \geq (2^{n}-2)/3$, by (\ref{4mu}) we obtain that 
$\mu(h)/d < 0.64$. Next assume that $e(h) = n-1$. Then \cite[Lemma 4.1]{TZ2} again implies that 
$|\chi(h)| \leq (2^{n-1}+3)/3$, whence $\mu(h)/d < 0.76$. 

Assume $d = (2^{n}-(-1)^{n})/3$. Then the involutive outer automorphism $\tau$ of $S$ inverts an 
element $t \in SU_{n}(2)$ represented by the matrix $\diag(\delta,\delta,\delta,1 , \ldots ,1)$ on $V$ but 
$\chi(t) \neq \chi(t^{-1})$. Thus $V$ is not stable under $\tau$, so without loss we may assume
$G \leq Z \cdot GU_{n}(2)$, whence $\mu(g)/d < 0.76$ as shown above.

Finally, we handle the remaining possibility $d = (2^{n}+2(-1)^{n})/3$. Recall we have shown 
that $\mu(h)/d < 0.64$ if $e(h) \leq n-2$. Again consider any $h \in GU_{n}(2) \setminus Z(G)$ with 
$e(h) = n-1$. Applying \cite[Lemma 4.1]{TZ1} we get $\chi(h) = \pm ((-2)^{n-1} + \beta)/3$ with 
$\beta = 2$ or $-1$. Assume $\beta = -1$. Then $|\chi(h)| = (2^{n-1}+(-1)^{n})/3 = d/2$.
In this case,  $\mu(h)/d < 3/4$ by (\ref{4mu}) (here the equality cannot hold as otherwise
$h^{3}$ cannot act scalarly on $V$, contrary to the containment $Z(GU_{n}(2)) \ni h^{3}$). 
If $\beta = 2$, then notice that $hZ(G)$ is a $2$-element in $G/Z(G)$ (in fact $h$ is a transvection 
modulo scalars). Thus we have shown that $\mu(h)/d < 3/4$ if $hZ(G)$ is not a $2$-element in 
$G/Z(G)$. Applying \cite[Lemma 3.4]{GT2}, we get $\mu(g)/d < (1+3/4)/2 = 7/8$, which completes the 
proof of Theorem \ref{main4}.

\section{Larsen's Conjecture for Small Rank Groups}\label{smallrank}

Throughout this section, $\FF$ is an algebraically closed field of characteristic $\ell$ and
$V = \FF^{d}$.

\begin{lemma}\label{lang}
{\sl Let $\GC = SL(V)$, $\ell > 0$, and let $q$ be a power of $\ell$. Let $G$ be a subgroup of 
$\GC$ that is irreducible on $V$.

{\rm (i)} The $G$-modules $V$ and $V^{(q)}$ are isomorphic if and only if some $\GC$-conjugate of
$G$ is contained in $SL_{d}(q)$.

{\rm (ii)} The $G$-modules $V^{*}$ and $V^{(q)}$ are isomorphic if and only if some $\GC$-conjugate of
$G$ is contained in $SU_{d}(q)$.}
\end{lemma}

\begin{proof}
(i) Assume $V|_{G} \simeq V^{(q)}|_{G}$. Consider the Frobenius map $F~:~x \mapsto x^{(q)}$ that 
raises any entry of the matrix $x \in \GC$ to its $q^{\mathrm {th}}$-power. Then there is some 
$x \in \GC$ such that $xgx^{-1} = F(g)$ for any $g \in G$. By the Lang-Steinberg Theorem,
one can find $y \in \GC$ such that $y^{-1}F(y) = x$. It follows that 
$ygy^{-1}$ is $F$-stable for any $g \in G$, i.e. $yGy^{-1} \leq \GC^{F} = SL_{n}(q)$. The converse
is immediate.

(ii) Argue similarly as in (i), taking $F~:~x \mapsto (\tn x^{-1})^{(q)}$.
\end{proof}

Now we proceed to prove Theorem \ref{small}. We will use the notation for fundamental weights of 
$\GC$, as well as the dimensions of restricted irreducible $\GC$-modules as given in \cite{Lu2}; 
in particular, $V = L(\om_{1})$. 

\begin{lemma}\label{sm-pos}
{\sl Theorem \ref{small} is true if $G$ is a positive dimensional closed subgroup.}
\end{lemma}

\begin{proof}
It suffices to show that if $\dim(G^{\circ}) > 0$ and $G$ satisfies $\MC \cap \MD$ then $G = \GC$.
Assume the contrary: $G < \GC$. If $\GC$ is of type $A_{2}$ or $A_{3}$ then we get an immediate
contradiction by Theorem \ref{red-sl}, since $G$ is irreducible on $\AC(V)$. Notice that $\MC \cap \MD$
implies $G$ is irreducible on $V$. (Indeed, $\MC$ implies $G$ is irreducible on 
$L(3\om_{1})$ which is $V \otimes \Vbb$ if $\ell = 2$. If $\ell \neq 2$, then $\MD$ implies 
$\MB$, so $G$ is irreducible on $L(2\om_{1}) = \SB(V)$, whence $G$ is irreducible on $V$). But the 
unipotent radical of $G^{\circ}$ fixes some nonzero vector of $V$, hence $G^{\circ}$ is reductive. 
Since the connected component of $Z(G^{\circ})$ acts scalarly and faithfully on $V$ and 
$Z(\GC)$ is finite, we see that $G^{\circ}$ is semisimple. In particular, $\GC$ cannot be of type
$A_{1}$. Thus $\GC = Sp_{4}(\FF)$ and $G^{\circ}$ is of type $kA_{1}$ with $k = 1,2$. If $\ell \neq 2$
then $G$ is irreducible on the adjoint module for $\GC$ which has dimension $10$. On the other hand, 
$G$ fixes the adjoint module for $G^{\circ}$, of dimension $3k$. So $\ell = 2$. In this case,
the $G$-composition factors on the adjoint module for $\GC$ are of dimensions $4$ and $1$. Hence 
$k = 2$, and the $G^{\circ}$-module $V$ decomposes as $A \otimes B$ with $\dim(A) = \dim(B) = 2$.
Since $\WB(A)$ and $\WB(B)$ are of dimension $1$ and $A^{(2)}$ and $B^{(2)}$ are 
irreducible, we must have 
$\TWB(V)|_{G^{\circ}} = A^{(2)} \otimes \WB(B) + B^{(2)} \otimes \WB(A)$. Thus 
$$(\Vbb \otimes \TWB(V))|_{G^{\circ}} = 
 (A^{(4)} + 2\WB(A^{(2)})) \otimes B^{(2)} \otimes \WB(B) +  
 (B^{(4)} + 2\WB(B^{(2)})) \otimes A^{(2)} \otimes \WB(A).$$  
Thus $L(2\om_{1}+\om_{2})|_{G^{\circ}}$ has composition factors of distinct dimensions $4$ and $2$, 
so $G$ is reducible on $L(2\om_{1}+\om_{2})$, contrary to $\MD$.
\end{proof}

From now on we may assume that $G$ is a finite subgroup of $\GC$. In the case $\ell = 0$, we may choose
a prime $p$ that does not divide $|G|$ and, by reducing $V$ modulo $p$, embed $G$ in the algebraic 
group $\GC_{p}$ in characteristic $p$, of the same type as of $\GC$. It is easy to see that, the 
irreducible modules of $\GC$ and $\GC_{p}$ with the same highest weight $\sum_{i}a_{i}\om_{i}$ and 
$a_{i}$ all bounded by a constant $C$ have the same dimension, provided we choose $p$ large enough
comparatively to $C$. In particular, $G$ satisfies $\MC \cap \MD$ if and only $G$ satisfies
$M_{6}(\GC_{p}) \cap M_{8}(\GC_{p})$. Hence without loss we may assume $\ell > 0$.  

\bigskip
{\bf Case I: $d = 2$.} Here $\GC = SL(V)$. Since $G$ is finite, we can find a smallest power
$q$ of $\ell$ such that $G \leq SL_{2}(q)$.

\medskip
Assume $G$ is irreducible on $\SB(V)$ and $\SD(V)$ if $\ell \geq 5$, on $\SB(V)$ if $\ell = 3$,
and $G$ satisfies $\MC \cap \MD$ if $\ell = 2$. We prove by induction on $(SL_{2}(q):G)$ that 
either $\OL(G) = SL_{2}(r)$ for some $r|q$, or $G = SL_{2}(5)$ and $\ell \neq 2,5$. The 
induction base is clear. For the induction step, we can find a maximal subgroup $M$ of $SL_{2}(q)$ 
that contains $G$. Clearly, $G$ is irreducible on $V$. Moreover, $G$ is primitive on $V$. 
(Indeed, if $G$ preserves the decomposition $V = A \oplus B$ with $\dim(A) = \dim(B) = 1$, then 
$G$ has a proper submodule $A \otimes B$ in $\SB(V)$ if $\ell \geq 3$, 
$A \otimes A^{(2)} \oplus B \otimes B^{(2)}$ in $L(3\om_{1}) = V \otimes \Vbb$ if $\ell = 2$.) 
Inspecting the list of maximal subgroups of $SL_{2}(q)$ as listed in \cite{Kle}, we arrive at one 
of the following possibilities.

(a) $M = SL_{2}(\qn) \cdot \kappa$ with $q = \qn^{b}$, $b$ a prime and $\kappa = (2,q-1,b)$.
Observe that the condition imposed on $G$ is also inherited by $H := G \cap SL_{2}(\qn)$. Since 
$(SL_{2}(\qn):H) < (SL_{2}(q):G)$, we are done by the induction hypothesis.   

(b) $q = \ell > 3$ and $2\AAA_{4} \leq M \leq 2\SSS_{4}$. This case is impossible as 
$G$ is irreducible on $\SD(V)$ of dimension $5$.

(c) $\ell \neq 2,5$ and $M = SL_{2}(5)$. If $G = M$, we are done. Assume $G < M$. If $\ell > 5$
then the irreducibility of $G$ on $\SD(V)$ implies $|G| > 25$, a contradiction. So $\ell = 3$. 
Since $G$ is irreducible on $\SB(V)$, we conclude that $G = SL_{2}(3)$. Thus the induction step
is completed.

\medskip     
Now assume that $G$ satisfies $\MC \cap \MD$. Then we can apply the above claim to see that 
either conclusion (iii) of Theorem \ref{small} holds, or $S := \OL(G) = SL_{2}(r)$. In the latter case
$r \geq 4$. (Otherwise $L((r+1)\om_{1})|_{S} = (V \otimes V^{(r)})|_{S} \simeq (V \otimes V)|_{S}$ 
contains $1_{S}$ and so $G$ is reducible on $L((r+1)\om_{1})$.) Conversely, assume $G$ satisfies 
either conclusion (iii) or (i) of Theorem \ref{small}. In the former case, direct computation shows
that $G$ satisfies $\MC \cap \MD$ (moreover, if $\ell = 3$ then $G$ fails $M_{14}(\GC)$ as 
$V^{\otimes 7} \supset V^{\otimes 5}$ and so $L(7\om_{1})$ and $L(5\om_{1})$ are inside 
$V^{\otimes 7}$; however, $L(7\om_{1})|_{G} \simeq L(5\om_{1})|_{G}$). Assume we are in the 
latter case. Notice that $V = V^{*}$ and all $\GC$-composition factors of $V^{\otimes k}$ with 
$0 \leq k \leq 4$ have highest weights $m\om_{1}$ with $0 \leq m \leq 4$. Hence $G$ satisfies 
$\MC \cap \MD$ if $q \geq 5$. If $q = 4$, then $G = S = SL_{2}(4)$. Since the $\GC$-module 
$V^{\otimes 4}$ decomposes as $\Vdd + 4\Vbb + 6L(0)$, $G$ also satisfies $\MC \cap \MD$. 

\bigskip
{\bf Case II: $d = 4$ and $\GC = Sp(V)$}. Since $G$ is finite, we can find a smallest power
$q$ of $\ell$ such that $G \leq Sp_{4}(q)$.

\medskip
Assume $G$ is irreducible on $L(\om_{2}) = \TWB(V)$ (dimension 5) and $L(4\om_{1}) = \SD(V)$ 
(dimension $35$) if $\ell \geq 5$, on $L(\om_{2}) = \TWB(V)$ (dimension 5) and 
$L(2\om_{1}+\om_{2})$ (dimension $25$) if $\ell = 3$, on $L(\om_{1}+\om_{2}) = V \otimes \TWB(V)$ 
and $L(2\om_{1}+\om_{2}) = \Vbb \otimes \TWB(V)$ (both of dimension $16$) if $\ell = 2$. We 
prove by induction on $(Sp_{4}(q):G)$ that either $\OL(G) = Sp_{4}(r)$ for some $r|q$, or 
$\ell = 2$ and $G = \ta B_{2}(r)$ with $r \geq 8$. The induction base is clear. For the induction 
step, consider a maximal subgroup $M$ of $Sp_{4}(q)$ that contains $G$. 

Claim that $G$ is irreducible and primitive on $V$. Indeed, assume $G$ is reducible on $V$. Since 
every $G$-composition factor of $\WB(V)$ has dimension $\leq 4$, the irreducibility of $G$ on
$\TWB(V)$ implies $\ell = 2$. But then every $G$-composition factor of $V \otimes \WB(V)$ has 
dimension $\leq 12$ and so $G$ is reducible on $V \otimes \TWB(V)$ of dimension $16$. Now assume $G$ 
acts transitively on the summands of a decomposition $V = \oplus^{m}_{i=1}A_{i}$ with 
$\dim(A_{i}) \geq 1$ and $m \geq 2$. First consider the case $\ell \neq 2$ and let $K$ be the 
kernel of the action of $G$ on $\{A_{1}, \ldots ,A_{m}\}$. Since $G/K \leq \SSS_{m} \leq \SSS_{4}$ 
and $\dim(\TWB(V)) = 5$, $K$ is irreducible on $\TWB(V)$. But obviously every composition factor of 
$\WB(V)|_{K} = \sum^{m}_{i=1}\WB(A_{i}) \oplus \sum_{1 \leq i < j \leq m}A_{i} \otimes A_{j}$
has dimension $\leq 3$, a contradiction. Assume $\ell = 2$ and $m = 4$. Then the $G$-module 
$V \otimes \WB(V)$ decomposes as the sum of $\sum_{i < j < k}A_{i} \otimes A_{j} \otimes A_{k}$
(with multiplicity $3$), and $\sum_{i < j}A_{i} \otimes A_{j} \otimes (A_{i} \oplus A_{j})$. Thus 
every $G$-composition factor of $V \otimes \WB(V)$ has dimension $\leq 12$ and so $G$ is reducible 
on $V \otimes \TWB(V)$. Finally, assume $\ell = 2$ and $m = 2$. Since $G$ is irreducible on 
$V \otimes \TWB(V)$ and $\WB(V)|_{G} = (A_{1} \otimes A_{2}) \oplus (\WB(A_{1}) \oplus \WB(A_{2}))$, 
we must have $\TWB(V)|_{G} = A_{1} \otimes A_{2}$. But then $V \otimes \TWB(V)$ contains the 
$G$-submodule $\WB(A_{1}) \otimes A_{2} + \WB(A_{2}) \otimes A_{1}$ of dimension $4$, a 
contradiction.  

Now we can inspect the possibilities for $M$ that act irreducibly primitively on $V$ (as listed in 
\cite{Kle}). If $M = Sp_{4}(\qn) \cdot \kappa$ with $\kappa = (2,q-1,b)$ and $q = \qn^{b}$, then, 
since the conditions imposed on $G$ are also inherited by $G \cap Sp_{4}(\qn)$, we are done by the 
induction hypothesis. If $q$ is even and $M = GO^{\pm}_{4}(q)$, then $V|_{[M,M]} = A \otimes B$
for some modules $A,B$ of dimension $2$, and arguing as at the end of the proof of Lemma 
\ref{sm-pos} we see that $M$ is reducible on $L(2\om_{1}+\om_{2})$. Another possibility is
$q = \ell > 2$ and $M = N_{Sp_{4}(q)}(E)$ with $E = 2^{1+4}_{-}$. Since the $M$-orbits on 
nontrivial linear characters of $E$ have length $10$ or $5$, $M$ is reducible on 
$\SD(V)$ if $\ell \geq 5$ and on $L(2\om_{1}+\om_{2})$ if $\ell = 3$. By the same
reason $M$ cannot be $\AAA_{6}$ or $\SSS_{6}$. If $M = 2\AAA_{7}$, then $\ell = 7$ and $M$ is 
reducible on $\SD(V)$. There remain two possibilities for $M$.

(a) $M = \ta B_{2}(q)$ with $q \geq 8$. Notice that the only maximal subgroups of $M$ that can act
irreducibly on $L(\om_{1}+\om_{2})$ are $\ta B_{2}(\qn)$ with $q = \qn^{b}$ and $\qn \geq 8$. 
Inducting on $(M:G)$ we conclude that $G = \ta B_{2}(r)$ for some $r \geq 8$.  

(b) $M$ is a cover of $L_{2}(q)$, $\ell \geq 5$, and $M$ acts irreducibly on $V$. Since $V$ is a 
symplectic module, $V|_{M}$ cannot be tensor decomposable, whence we may assume it is isomorphic to
$\SC(U)$, where $U$ is the natural module for $SL_{2}(\FF)$. If $\om$ denotes the unique fundamental
weight of $SL_{2}(\FF)$, then any $SL_{2}(\FF)$-composition factor of $(\SC(U))^{\otimes 4}$ has 
highest weight $k\om$ with $k \leq 12$ and so it has dimension $\leq 13$. In particular, $M$ is
reducible on $\SD(V)$. Thus the induction step is completed.    

\medskip     
Now assume that $G$ satisfies $\MC \cap \MD$. Then we can apply the above claim to see that 
either $G = \ta B_{2}(r)$, or $S := \OL(G) = Sp_{4}(r)$. Notice that when $\ell = 2$,
$$V^{\otimes 3} = 4L(\om_{1}) + 2L(\om_{1}+\om_{2}) + L(3\om_{1}),$$  
$$V^{\otimes 4} = 36L(0) + 24L(\om_{2}) + 8L(2\om_{1}) + 6L(2\om_{2}) + 
 4L(2\om_{1}+\om_{2}) + L(4\om_{1}).$$ 
So in the former case $r \geq 32$, as otherwise the distinct $\GC$-composition factors $L(4\om_{1})$ 
and $L(\om_{2})$ are isomorphic over $M$, contrary to $\MD$. By the same reason, $r \neq 2$ in the 
latter case. Assume $r = 3$. Then $L(4\om_{1})|_{S} = (V \otimes \Vcc)|_{S} = (V \otimes V)|_{S}$ 
contains $1_{S}$, and so $G$ must be reducible on $L(4\om_{1})$. Conversely, assume $G$ satisfies 
either conclusion (i) or (ii) of Theorem \ref{small}. Notice that $V = V^{*}$ and all 
$\GC$-composition factors of $V^{\otimes k}$ with 
$0 \leq k \leq 4$ have highest weights $a\om_{1}+b\om_{2}$ with $0 \leq a,b \leq 4$. Hence $G$ 
satisfies $\MC \cap \MD$ if $\OL(G) = Sp_{4}(q)$ with $q \geq 5$ or $G = \ta B_{2}(q)$ with 
$q \geq 32$. The above decompositions for $V^{\otimes 3}$ and $V^{\otimes 4}$ shows that 
$Sp_{4}(4)$ also satisfies $\MC$ and $\MD$.

\bigskip
In what follows, $SL^{\eps}_{n}(q)$ denotes $SL_{n}(q)$ if $\eps = +$ and $SU_{n}(q)$ with 
$\eps = -$.

\medskip
{\bf Case III: $d = 3$ and $\GC = SL(V)$}. Since $G$ is finite, we can find a smallest power
$q$ of $\ell$ such that $G \leq SL^{\eps}_{3}(q)$ for some $\eps = \pm$. 

\medskip
Assume $G$ is irreducible on $L(\om_{1}+\om_{2}) = \AC(V)$ (dimension $8-\delta_{3,\ell}$) and,
additionally on $L(3\om_{1}) = \SC(V)$ (dimension $10$) if $\ell \geq 5$. We 
prove by induction on $(SL^{\eps}_{3}(q):G)$ that either $\OL(G) = SL^{\al}_{3}(r)$ for some $r|q$
and $\al = \pm$, or $\ell = 2$ and $G \rhd 3^{1+2}:Q_{8}$, or $\ell = 3$ and $G/Z(G) = L_{2}(7)$, 
or $\ell \neq 3$ and $G \rhd 3\AAA_{6}$,
or $\ell = 5$ and $G \rhd 3\AAA_{7}$. The induction base is clear. For the induction 
step, consider a maximal subgroup $M$ of $SL^{\eps}_{3}(q)$ that contains $G$. Applying 
Theorem \ref{red-sl} to $G$, we need to consider only the case $G$ is almost quasisimple
with $\soc(G/Z(G)) \notin Lie(\ell)$. Applying Theorem \ref{red-sl} to $M$, we see that 
either $M$ is almost quasisimple with $\soc(M/Z(M)) \notin Lie(\ell)$, or 
$M = SL^{\beta}_{3}(\qn) \cdot \kappa$ with $\beta = \pm$, $q = \qn^{b}$ and 
$\kappa |(3,q^{2}-1,b)$. In the latter case we are done by the induction hypothesis,     
since the conditions imposed on $G$ are also inherited by $G \cap SL^{\beta}_{3}(\qn)$. Thus 
$G$ and $M$ are both almost quasisimple, with a unique nonabelian composition factor not a 
Lie type group in characteristic $\ell$. Now we can inspect the list of maximal subgroups of 
$SL^{\eps}_{3}(q)$ as given in \cite{Kle} and arrive at one of the following subcases.

(a) $M/Z(M) = L_{2}(7)$ and $\ell \neq 2,7$. Here, if $\ell > 3$ then $M$ is reducible on 
$\SC(V)$. On the other hand, if $\ell = 3$ then $G = M$.
 
(b) $M = 3\AAA_{6}$ and $\ell \neq 3,5$. Here, if $G < M$, then $\AAA_{5} \lhd G$, whence 
$G$ is reducible on $\AC(V)$. Thus $G = M$.

(c) $M = 3\AAA_{7}$ or $M = 3M_{10}$ and $\ell = 5$. Now the irreducibility of
$G \cap [M,M]$ on $\SC(V)$ forces either $G = M = 3\AAA_{7}$ or $G \rhd 3\AAA_{6}$. The induction
step is completed.

\medskip
Now we assume that $G$ satisfies $\MC \cap \MD$ and apply the above claim to $G$. 
Assume $G \rhd N := 3^{1+2}: Q_{8}$ (so $\ell = 2$). Then 
$L(3\om_{1})|_{N} = (V \otimes \Vbb)|_{N} = (V \otimes V^{*})|_{N}$ contains $1_{N}$, so 
$G$ is reducible on $L(3\om_{1})$. The same argument applies to the case where $G \rhd 3\AAA_{6}$ 
and $\ell = 2$. If $\ell = 3$ and $G/Z(G) = L_{2}(7)$, then $G \rhd L := L_{2}(7)$, and 
$L(4\om_{1})|_{L} = (V \otimes \Vcc)|_{L} = (V \otimes V^{*})|_{L}$ contains $1_{L}$, so $G$ is 
reducible on $L(4\om_{1})$. If $\ell \geq 5$ and $G \rhd 3\AAA_{6}$ or $G = 3\AAA_{7}$, then $G$ 
cannot act irreducibly on $L(2\om_{1}+2\om_{2})$ (of dimension $19$ if $\ell = 5$ and $27$ if 
$\ell > 5$). We conclude that $L := \OL(G) = SL^{\al}_{3}(r)$ for some $r|q$ and $\al = \pm$. Observe 
that $r \geq 4$. Indeed, if $r = 2$, then $L(3\om_{1})|_{L} = (V \otimes \Vbb)|_{L}$ contains $1_{L}$
if $\al = -$, and contains $L$-composition factors of distinct multiplicities if $\al = +$, so $G$ 
is reducible on $L(3\om_{1})$. If $r = 3$, then $L(4\om_{1})|_{L} = (V \otimes \Vcc)|_{L}$ contains 
$L$-composition factors of distinct dimensions, so $G$ is reducible on $L(4\om_{1})$. Conversely, 
assume $\OL(G) = SL^{\al}_{3}(r)$ with $4 \leq r|q$ and $\al = \pm$. Notice that all 
$\GC$-composition factors of $V^{\otimes k} \otimes (V^{*})^{\otimes l}$ with $k+l \leq 4$ have 
highest weights $a\om_{1}+b\om_{2}$ with $a+b \leq 4$. Hence $G$ satisfies $\MC \cap \MD$ if 
$r \geq 5$. Consider the case $r = 4$ and $\ell = 2$. Here, all $\GC$-composition factors of 
$V^{\otimes 3}$, $V^{\otimes 2} \otimes V^{*}$, and 
$V^{\otimes 2} \otimes (V^{*})^{\otimes 2}$ have highest weights $a\om_{1}+b\om_{2}$ with 
$0 \leq a,b \leq 3$. Furthermore,   
$$V^{\otimes 4} = 8L(\om_{1}) + 6L(2\om_{2}) + 4L(2\om_{1}+\om_{2}) + L(4\om_{1}),$$
$$V^{\otimes 3} \otimes V^{*} = 5L(2\om_{1}) + 8L(\om_{2}) + 
  2L(\om_{1}+2\om_{2}) + L(3\om_{1} + \om_{2}).$$
It follows that $SU_{3}(4)$ satisfies $\MC \cap \MD$, but $SL_{3}(4)$ fails $\MD$ as 
$L(4\om_{1})$ and $L(\om_{1})$ are isomorphic over $SL_{3}(4)$. Notice that  
$H := N_{\GC}(K) = K \cdot 3$ for $K := SL_{3}(4)$. It remains to show that $H$ satisfies 
$\MC \cap \MD$, which is equivalent to showing that $L(4\om_{1})$ and $L(\om_{1})$ are not 
isomorphic over $H$. But the last statement holds by Lemma \ref{lang}.     

\bigskip
{\bf Case IV: $d = 4$ and $\GC = SL(V)$}. Since $G$ is finite, we can find a smallest power
$q$ of $\ell$ such that $G \leq SL^{\eps}_{4}(q)$ for some $\eps = \pm$. 

\medskip
Assume $G$ is irreducible on $L(\om_{1}+\om_{2}) = \AC(V)$ (dimension $15-\delta_{2,\ell}$) and 
additionally, on $L(4\om_{1}) = \SD(V)$ (dimension $35$) if $\ell \geq 5$, on 
$L(2(\om_{1}+\om_{3}))$ (of dimension $69$) if $\ell = 3$, on $L(2\om_{1}+\om_{2})$ (of dimension $24$) 
if $\ell = 2$. We prove by induction on $(SL^{\eps}_{4}(q):G)$ that either $\OL(G) = SL^{\al}_{4}(r)$ 
for some $r|q$ and $\al = \pm$, or $G/Z(G) = \AAA_{7}$. The induction 
base is clear. For the induction step, consider a maximal subgroup $M$ of $SL^{\eps}_{4}(q)$ that 
contains $G$. Applying Theorem \ref{red-sl} to $G$, we need to consider only the case where $G$ 
is almost quasisimple with $\soc(G/Z(G)) \notin Lie(\ell)$. Applying Theorem \ref{red-sl} to $M$, 
we see that either $M$ is almost quasisimple with $\soc(M/Z(M)) \notin Lie(\ell)$, or 
$M = SL^{\beta}_{4}(\qn) \cdot \kappa$ with $\beta = \pm$, $q = \qn^{b}$ and $\kappa$ is a $2$-power. 
In the latter case we are done by the induction hypothesis,     
since the conditions imposed on $G$ are also inherited by $G \cap SL^{\beta}_{4}(\qn)$. Thus 
$G$ and $M$ are both almost quasisimple, with a unique nonabelian composition factor not a 
Lie type group in characteristic $\ell$. Now we can inspect the list of maximal subgroups of 
$SL^{\eps}_{3}(q)$ as given in \cite{Kle} and arrive at one of the following subcases. Notice
that the irreducibility condition on $G$ implies $|G/Z(G)| \geq 24^{2}$.  

(a) $M/Z(M) = \AAA_{7}$ and $q = \ell \neq 7$. Clearly, $G = M$ by order.
 
(b) $M = Sp_{4}(3)$ and $q = \ell \neq 2,3$. But in this case $M$ is reducible on 
$\SD(V)$.

(c) $M = 4L_{3}(4)$ and $q = 3$. In this case $M$ is reducible on $L(2(\om_{1}+\om_{3}))$.  
The induction step is completed.

\medskip
Now we assume that $G$ satisfies $\MC \cap \MD$ and apply the above claim to $G$.
Assume $G/Z(G) = \AAA_{7}$. If $\ell \geq 3$, then $\dim(L(2(\om_{1}+\om_{3}))) \geq 69$ and so 
$G$ is reducible on $L(2(\om_{1}+\om_{3})))$, contrary to $\MD$. If $\ell = 2$, then 
$V^{\otimes 3} \otimes V^{*}$ contains $L(3\om_{1}+\om_{3})$ of dimension $56$, and so 
$G$ cannot satisfy $\MD$. Thus $L := \OL(G) = SL^{\al}_{4}(r)$ for some $r|q$ and $\al = \pm$.
Observe that $r \geq 4$. Indeed, if $r = 2$, then $L(3\om_{1})|_{L} = (V \otimes \Vbb)|_{L}$ contains 
$L$-composition factors of distinct dimensions, so $G$ is reducible on $L(3\om_{1})$. By 
the same reason, if $r = 3$ then $G$ is reducible on $L(4\om_{1})$. Conversely, assume 
$\OL(G) = SL^{\al}_{4}(r)$ with $4 \leq r|q$ and $\al = \pm$. Notice that all $\GC$-composition 
factors of $V^{\otimes k} \otimes (V^{*})^{\otimes l}$ with $k+l \leq 4$ have highest weights 
$a\om_{1}+b\om_{2}$ with $a+b \leq 4$. Hence $G$ satisfies $\MC \cap \MD$ if $r \geq 5$. Consider
the case $r = 4$ and $\ell = 2$. Here, all $\GC$-composition factors of $V^{\otimes 3}$, 
$V^{\otimes 2} \otimes V^{*}$, and $V^{\otimes 3} \otimes V^{*}$ have highest weights 
$a\om_{1}+b\om_{2}$ with $0 \leq a,b \leq 3$. Furthermore,   
$$V^{\otimes 4} = 8L(0) + 6L(2\om_{2}) + 4L(2\om_{1}+\om_{2}) + 8L(\om_{1}+\om_{3}) + L(4\om_{1}),$$
$$V^{\otimes 2} \otimes (V^{*})^{\otimes 2} = 10L(0) + 8L(\om_{1}+\om_{3}) + 
  2L(2\om_{1}+\om_{2}) + 2L(\om_{2}+2\om_{3}) + 4L(2\om_{2}) + L(2(\om_{1}+\om_{3})).$$
It follows that $SL^{\al}_{3}(4)$ satisfies $\MC \cap \MD$.   
  
Theorem \ref{small} has been proved.

\section{Generation Results}\label{generation}

We first note the following well known result.

\begin{lemma}\label{2var}
{\sl Let $\FF$ be an algebraically closed field and $V$ a finite dimensional vector space over
$\FF$ of dimension $d$. 

{\rm (i)} Let $\RL_{k,r}(V)$ be the set of $r$-tuples of elements in $GL(V)$ that fix a common
$k$-dimensional subspace. Then $\RL_{k,r}(V)$ is a closed subvariety of $GL(V)^{r}$.

{\rm (ii)} $\IL_{r}(V):= \{(g_{1}, \ldots, g_{r}) \in GL(V)^{r}
\mid \langle g_{1}, \ldots, g_{r}\rangle \mbox{ is irreducible on
} V\}$ is an open subvariety of $GL(V)^{r}$.}
\end{lemma}

\begin{proof}
The set in (i) is the set of $r$-tuples having a common fixed point on a projective variety
(the Grassmannian) and so is closed. The set in (ii) is the complement of a finite union of closed 
sets. Alternatively, the complement of the second condition is that the dimension of the 
$\FF$-subalgebra of $\End(V)$ generated by the $g_{i}$ is less than $d^{2}$. This is given by setting
certain determinants equal to $0$, a closed condition.
\end{proof}

In particular, if $V$ is a $G$-module, we intersect the sets above with $G^{r}$ to obtain closed
and open subvarieties of $G^{r}$. We will abuse notation and use the same notation for the
intersection of the sets above with $G^{r}$ when we are talking about $G$-modules.

Using results about generation, we easily obtain:

\begin{lemma}\label{ir}
{\sl If $\GC$ is a connected group, $V$ is irreducible and $r >1$, then $\IL_{r}(V)$
is a nonempty open subset of $\GC^{r}$ (and in particular it is dense).}
\end{lemma}

\begin{proof}
There is no loss in extending $\FF$ if necessary and assume that $\FF$ is not algebraic over a
finite field. We then use the fact that every semisimple algebraic group can be topologically
generated (in the Zariski topology) by two elements (see \cite{gurnewton}).
Alternatively, in positive characteristic, we can use the result of Steinberg \cite{steinberg}
that every simple finite Chevalley group can be generated by $2$ elements.
\end{proof}

It is sometimes convenient to use an indecomposable module to detect generators for our group 
(a particularly useful case is the orthogonal group inside the symplectic group in characteristic
$2$). We restrict attention to {\it multiplicity free} modules (i.e. those modules where each 
composition factor occurs precisely once).

\begin{lemma}\label{indecomp}
{\sl Let $\GC$ be a semisimple algebraic group and $V$ a multiplicity free rational $\GC$-module.
Let $r > 1$. The set $\DL_{r}(V)$ of $r$-tuples of elements in $\GC$ that have the same invariant 
subspaces on $V$ as $\GC$ is a non-empty open subvariety of $\GC^{r}$.}
\end{lemma}

\begin{proof}
We induct on the composition length of $V$ as a $\GC$-module.  We have already proved this for 
simple modules. Note that the multiplicity free assumption implies that there are only finitely 
many $\GC$-submodules inside $V$.

Consider the case of length $2$. We will show that $\DL_{r}(V)$ is the set of $r$-tuples that 
generate the same algebra in $\End(V)$. This is clearly an open condition (and non-empty since 
$\GC$ can be topologically generated by two elements -- if necessary passing to a field that is 
not algebraic over a finite field).

Clearly, any $r$-tuple that generates the same algebra leaves invariant only those subspaces that 
are $\GC$-invariant. If $V$ is semisimple as a $\GC$-module, the result is clear. Let $S$ be
the socle of $V$ and $T := V/S$, of dimensions $c$ and $d$. Since $S$ and $T$ are nonisomorphic, 
it follows that the algebra $A$ generated by $\GC$ has a nilpotent radical $N$ of dimension equal 
to $\dim(S) \cdot \dim (T)$ and $A/N \simeq Mat_{c}(\FF) \oplus Mat_{d}(\FF)$. If the algebra
$B$ generated by some $r$-tuple is proper, then either $(B + N)/N = A/N$ or $B$ is a complement to 
$N$ (since $N$ is a minimal $2$-sided ideal of $A$). In the former case the $r$-tuple does not 
generate an irreducible group on either $S$ or $T$ and so leaves invariant some subspace that is 
not $\GC$-invariant. In the latter case, $B$ is semisimple and so acts completely reducibly
and so also has extra invariant subspaces.

Now we assume that the composition length is at least $3$. We claim that $\DL_{r}(V)$ is the 
intersection of the varieties $\DL_{r}(V/S)$ and $\DL_{r}(W)$ where $S$ ranges over simple
submodules of $V$ and $W$ the maximal submodules of $V$. Clearly, $\DL_{r}(V)$ is contained in this 
finite intersection of open subvarieties and this intersection is a non-empty open subvariety.

So assume that there is an $r$-tuple in the intersection and let $H$ be the subgroup generated by 
it. Suppose that $HW = W$ but $W$ is not $\GC$-invariant and choose such a $W$ of minimal dimension.
Then $S + W$ is $\GC$-invariant for any simple $\GC$-module $S$, whence $V = S + W$ (otherwise
this is contained in a maximal $\GC$-submodule). Since $V/S$ is not simple (because the composition 
length of $V$ is at least $3$), it follows that $W$ is not simple for $H$. Let $U$ be a simple 
$H$-submodule of $W$. By the choice of $W$, $U$ is $\GC$-invariant and then we pass to $V/U$ and 
see that $W/U$ is $\GC$-invariant, a contradiction.
\end{proof}

Let $\GC$ be a simple algebraic group. By a {\it subfield subgroup} we mean any subgroup $G$ with
$\OL(G) = \GC^{F}$, the (finite) fixed point subgroup for some (twisted or untwisted)
Frobenius endomorphism $F$ of $\GC$ (in particular, this includes the triality groups for type
$D_{4}$), and we call it {\it good} if $\GC^{F}$ is none of the groups
$\ta B_{2}(2)$, $\ta G_{2}(3)$, or $\ta F_{4}(2)$. We write this
subgroup as $G(q)$ where $q$ is the absolute value of the eigenvalues of $F$ on the character group
of a maximal $F$-invariant torus. If the characteristic is $0$, then there are no subfield
subgroups. Let $\IL_{r}(\GC, N)$ denote the set of $r$-tuples $(g_{1}, \ldots ,g_{r}) \in \GC^{r}$
such that the closure of $\langle g_{1}, \ldots ,g_{r} \rangle$ contains a good subfield subgroup
with $q > N$. Let $\IL_{r}(\GC) := \IL_{r}(\GC, 1)$.

\begin{corol}\label{ir-clas}
{\sl Let $\FF$ be an algebraically closed field of characteristic $\ell$. Let $\GC$ be a
simple simply connected classical group with natural module $V$ of dimension $d \geq 5$.

{\rm (i)} If $\ell = 0$ or $\ell > 71$, then $\cap_{W \in \SL} \IL_{r}(W) = \IL_{r}(\GC)$,
where $\SL$ is the set of $\GC$-composition factors of $V \otimes V^{*}$, $V^{\otimes 4}$, and
$L(3\omega_{1})$.

{\rm (ii)} If $\ell > 5$, $\cap_{W \in \SL} \IL_{r}(W) = \IL_{r}(\GC)$, where $\SL$ is the set of
$\GC$-composition factors of $V^{\otimes (s-j)} \otimes (V^{*})^{\otimes j}$, $0 \leq j \leq s$
and $s = 5,6$.

{\rm (iii)} If $\ell \in \{3,5\}$, then $\cap_{W \in \SL} \IL_{r}(W) =\IL_{r}(\GC, 5)$
where $\SL$ is the set of $\GC$-composition factors of 
$V^{\otimes (s-j)} \otimes (V^{*})^{\otimes j}$, $0 \leq j \leq s$ and $s = 5,6$.

{\rm (iv)} If $\ell = 2$, $\cap_{W \in \SL} \IL_{r}(W) = \IL_{r}(\GC, 5)$,
where $\SL$ is the set of $\GC$-composition factors of 
$V^{\otimes (s-j)} \otimes (V^{*})^{\otimes j}$, $0 \leq j \leq s$ and $s = 5,6$, plus, if 
$\GC = Sp(V)$, $L(\om_{1}+\om_{d/2})$.}
\end{corol}

\begin{proof}
These results follow from Theorem \ref{main2} observing that the only groups acting irreducibly
on the modules in $\SL$ are the subfield subgroups.  In the first two cases, all modules
are restricted and so all subfield subgroups do act irreducibly.

Consider the third case. Then every module in $\SL$ is either a Frobenius twist of a restricted
module or of a tensor product $W_{1} \otimes W_{2}^{(\ell)}$ with $W_{1}$ and $W_{2}$ restricted.
Any subfield subgroup other than a subfield group over the prime field is irreducible (and those
are not) on such a module, whence the result.

Assume $\ell = 2$. By Lemma \ref{sp-o3}, any subgroup of $O(V)$ is reducible on the
$Sp(V)$-module $L(\om_{1}+\om_{d/2})$, so only subfield subgroups may act irreducibly on all
$W \in \SL$. In addition to the non-restricted modules mentioned in case (iii), we also need to
consider the additional possibility of $W_{1} \otimes W_{2}^{(4)}$.  Again all subfield subgroups
are irreducible with the exception of groups defined over the prime field or the field
of $4$ elements. One may formulate another version of (iv) using the indecomposable module
of dimension $2n+1$ for $Sp_{2n}(\FF)$, see Lemma \ref{indecomp}.
\end{proof}

In what follows we denote by $St$ the {\it basic Steinberg representation} of $\GC$, that is,
$L(\om)$ with $\om = (\ell-1)\sum^{n}_{i=1}\om_{i}$ if $\ell > 0$. Abusing the notation,
we also denote by $St$ the module $L(\om)$ with $\om = 3\sum^{n}_{i=1}\om_{i}$ if $\ell = 0$. Also,
we denote by $V_{ad}$ the set of nontrivial $\GC$-composition factors of the adjoint module for $\GC$.
Notice that $|V_{ad}| = 1$ (in which case we also let $V_{ad}$ denote this unique module), except
for the cases $(\GC,\ell) = (F_{4},2)$ , $(G_{2},3)$, in which $V_{ad}$ consists of two nontrivial
composition factors of same dimension $26$, resp. $7$. We choose natural module $V$ of $\GC$ to
have dimension $6$ for type $A_{3} \simeq D_{3}$, $5$ for type $B_{2}$ and $\ell \neq 2$.

\begin{theor}\label{restr1}
{\sl Let $\FF$ be an algebraically closed field of characteristic $\ell$. Let $\GC$ be a simple
simply connected algebraic group over $\FF$, and let $\SL := V_{ad} \cup \{St\}$. Assume
$G$ is a proper closed subgroup of $\GC$. Then $G$ is irreducible on all $X \in \SL$ if
and only if one of the following holds.

{\rm (i)} $G$ is a good subfield subgroup of $\GC$, that is, $\OL(G) = \GC^{F}$ for some
Frobenius endomorphism $F$ on $\GC$, and
$\GC^{F} \notin \{\ta B_{2}(2), \ta G_{2}(3), \ta F_{4}(2)\}$.

{\rm (ii)} $\GC=A_{1}$, $\ell=2$, $G$ is any closed subgroup not contained in a Borel subgroup.

{\rm (iii)} $\GC = A_{1}$, $\ell=3$, and $G = 2\AAA_{5}$.

{\rm (iv)} $\GC = A_{2}$, $\ell=2$, and $G = 3\AAA_{6}$.

{\rm (v)} $\GC = C_{2} = Sp(V)$, $\ell = 2$, $G = O(V)$.

{\rm (vi)} $\GC = C_{2}$, $\ell = 2$, $G = O^{\pm}_{4}(q)$ with $q= 2^{f} \geq 4$.

{\rm (vii)} $\GC = G_{2}$, $\ell = 2$, $G = J_{2}$.}
\end{theor}

\begin{proof}
1) First suppose that $\HC$ is a proper positive dimensional closed subgroup of $\GC$.
Exclude the cases $(\GC,\ell) = (A_{1}, 2)$, $(C_{n}, 2)$, $(F_{4}, 2)$ and $(G_{2}, 3)$
for a moment. Then the unique nontrivial irreducible $\GC$-constituent $V_{ad}$ of the
adjoint module has codimension at most $2$ in the adjoint module. It follows that $\HC$ does not
irreducibly on this constituent, as the adjoint module of $\HC$ is invariant and has smaller
dimension. The same is true for $N_{\GC}(\HC)$; in particular, any Lie imprimitive subgroup of
$\GC$ is reducible on $V_{ad}$.

Assume $(\GC,\ell) = (C_{n},2)$ with $n \geq 2$. Then the only $\HC$ with adjoint module of dimension
at least $\dim(V_{ad})$ is isomorphic to $D_{n} = \Omega(V)$ with $\GC = Sp(V)$. If $n \geq 3$ then
$O(V)$ is reducible on $St$ by Lemma \ref{sp-o3}, excluding the case of $D_{n}$. It remains to
consider $\HC = \Omega(V)$ (fixing some quadratic form on $V$) or $Sp_{2}(\FF) \times Sp_{2}(\FF)$
(fixing some nondegenerate $2$-space of $V$) inside $C_{2} = Sp(V)$. Observe that
$St = V \otimes \TWB(V)$. In the first case, $V|_{\HC} = A \otimes B$ with $A$, $B$ natural
$2$-dimensional modules for $A_{1}$ and $\HC \simeq A_{1} \times A_{1}$. Hence
$\TWB(V)|_{\HC} = A^{(2)} \otimes 1 \oplus 1 \otimes B^{(2)}$, and so
\setcounter{equation}{0}
\begin{equation}\label{st-o}
  St|_{\HC} = (A \otimes A^{(2)}) \otimes B \oplus A \otimes (B \otimes B^{(2)}).
\end{equation}
It follows that $St$ is reducible over $\HC$ but irreducible over $N_{\GC}(\HC) = O(V)$, leading
to the case (v). In the second case, $V|_{\HC} = A \oplus B$, so $\TWB(V)|_{\HC} = A \otimes B$.
Therefore, $St|_{\HC}$ contains irreducibles $1 \otimes B$ and $A^{(2)} \otimes B$ of dimensions
$2$ and $4$, whence $St$ is reducible over $N_{\GC}(\HC)$.

Next we assume $(\GC,\ell) = (C_{2},2)$ and $G$ is a finite subgroup of $\GC$ that fixes some
quadratic form on $V$ and is irreducible on $St$. Then $G < O(V) = \HC.2$ with
$\HC = \Omega(V)$. We will identify $\HC$ with $Sp_{2}(\FF) \times Sp_{2}(\FF)$ and write
$V|_{\HC} = A \otimes B$ as above. Since $\Omega(V)$ is reducible on $St$, $G = H.2$ with
$H := G \cap \HC$. Let $\pi_{1}$, resp. $\pi_{2}$, denote the projection of $\HC$ onto the first,
resp. second, simple component of $\HC$. Any element in $O(V) \setminus \HC$ permutes these two
components, hence $\pi_{1}(H) = \pi_{2}(H) = S$ for some finite subgroup $S < Sp_{2}(\FF)$.
Clearly, $H \leq S \times S$ is irreducible on $(A \otimes A^{(2)}) \otimes B$, cf. (\ref{st-o}),
whence $S$ is irreducible on $A \otimes A^{(2)}$. It is not difficult to see that this forces
$S = SL_{2}(r)$ for some $r = 2^{g} \geq 4$. Since $S$ is simple and since $H$ projects onto $S$
under both $\pi_{1}$, $\pi_{2}$, one sees that $H = S \times S$ or $H$ can be identified with the
diagonal subgroup $\{(s,s) \mid s \in S\}$ (after applying a suitable automorphism to the second
$S$ in $S \times S$). Suppose we are in the first case. Recall any $t \in G \setminus H$ interchanges
the two copies of $S$ in $H = S \times S$ and stabilizes the $H$-module $A \otimes B$. It follows
that $A \simeq B$ as $H$-modules. Thus $H = \Omega^{+}_{4}(q)$ with $q = r$ (see for instance
\cite[p. 45]{KL}), and one can check that $G = O^{+}_{4}(q)$. Now assume we are in the second case.
Since the {\it inner} tensor product $(A \otimes A^{(2)}) \otimes B$ is irreducible over $H$, the
$S$-module $B$ is not isomorphic to $A$ nor to $A^{(2)}$. On the other hand, any
$t \in G \setminus H$ stabilizes the $H$-module $A \otimes B$. One readily shows that
$r = q^{2} \geq 16$, and $B = A^{t} = A^{(q)}$. Thus $H = \Omega^{-}_{4}(q)$ (see for instance
\cite[p. 45]{KL}), and one can check that $G = O^{-}_{4}(q)$. Conversely, any subgroup
$O^{\pm}_{4}(q)$ with $q \geq 4$ is irreducible on all irreducible $\GC$-restricted modules. We have
arrived at the possibility (vi).

2) Here we consider the case $(\GC,\ell) = (G_{2},3)$, $(F_{4},2)$ and assume $G \leq N_{\GC}(\HC)$
is irreducible on all $X \in V_{ad}$ for some proper positive dimensional closed connected subgroup
$\HC$ of $\GC$. We already mentioned that the adjoint module $Z$ of $\GC$ consists of two composition
factors $X,Y$ of dimension $\dim(Z)/2$, and $V_{ad} = \{X,Y\}$. Now $G$ fixes the adjoint module $T$
of $\HC$, which is a proper submodule of the $\HC$-module $Z$. Since $G$ is irreducible on both $X$
and $Y$, we conclude that $\dim(T) = \dim(Z)/2$.

Observe that $\HC$ is not a parabolic subgroup of $\GC$. Indeed, $\GC$ acts faithfully on $X$, so
$G$ and $N_{\GC}(\HC)$ act faithfully and irreducibly on $X$. By Clifford's theorem, $\HC$ is
semisimple and faithful on $X$. On the other hand, the unipotent radical $\HC_{u}$ of $\HC$ acts
trivially on every irreducible $\FF\HC$-module and so on $X$, whence $\HC_{u} = 1$.

Now assume $(\GC,\ell) = (G_{2},3)$. We have shown that $\HC$ is a non-parabolic proper closed
connected subgroup of $\GC$ with adjoint module of dimension $7$. The maximal proper closed
connected subgroups of $\GC$ are known (see \cite{LiS}), and so one easily deduces that
$\HC = A_{2}$, which comes from a maximal rank subgroup. One can check that $St|_{A_{2}}$ contains a
subquotient of dimension $27 \cdot 6$. Meanwhile, $\dim(St) = 729$, so $St$ is reducible over
$N_{\GC}(\HC) = \HC \cdot 2$. Thus $G$ is reducible on $St$.

Next assume $(\GC,\ell) = (F_{4},2)$. Arguing as above we obtain $\HC = D_{4}$ (with two versions for
$D_{4}$: long-root and short-root), $B_{4}$ or $C_{4}$. First consider the case of $B_{4}$ and let
$\pi_{1}, \ldots ,\pi_{4}$ denote the fundamental weights of $B_{4}$. Since this is a maximal rank
subgroup, we can check that
$St|_{B_{4}} \simeq (L(\om_{1}+\om_{2}) \otimes L(\om_{3}+\om_{4}))|_{B_{4}}$ contains a subquotient
$$L(\pi_{1}+\pi_{2}+\pi_{3}) \otimes L(4\pi_{1}+ \pi_{4}) \simeq L(\pi_{1}+\pi_{2}+\pi_{3})
  \otimes L(\pi_{1})^{(4)} \otimes L(\pi_{4})$$
of dimension $2^{12} \cdot 2^{3} \cdot 2^{4}$. Since $\dim(St) = 2^{24}$, $St$ is reducible over
$B_{4}$. Also, since the long-root subgroup $D_{4}$ is contained in $B_{4}$, $St$ is also reducible
over the corresponding $N_{\GC}(D_{4}) = D_{4} \cdot \SSS_{3}$. A similar argument applies to the
case of $C_{4}$ and the short-root subgroup $D_{4}$. Thus $G$ is again reducible on $St$.

3) Consider the case $\GC$ is an exceptional group. By the virtue of 1) and 2), it remains to
consider finite Lie primitive subgroups $G$ of $\GC$. If $G$ is not a subfield subgroup, then we can
use \cite{LiS} to get all possibilities for $G$ and check that $|G| < (\dim(St))^{2}$, with only one
exception (vii). Notice that $J_{2}$ is irreducible on all restricted $\GC$-modules in the case of
(vii). Finally, assume that $\OL(G) = \GC^{F}$ for some Frobenius map $F$ on $\GC$. Aside from
the two exceptions $\ta G_{2}(3)$ and $\ta F_{4}(2)$ (which are reducible on $St$), all other
subfield subgroups are indeed irreducible on all irreducible restricted $\GC$-modules.

4) Now we assume that $\GC$ is classical with natural module of dimension $d \geq 5$. If $\ell > 0$
then $\dim(St) \geq \ell^{d(d-2)/4}$, and if $\ell = 0$ then $\dim(St) \geq 4^{d(d-2)/4}$ by our
choice of $St$. Assume $G$ is irreducible on all $X \in \SL$. By 1) we may assume that $G$ is finite.
Applying the results of \S\ref{red-thms} to $G$, we see that either $G$ is a subfield subgroup,
or $G$ normalizes a $p$-group of symplectic type, or $G$ normalizes an elementary abelian
$2$-group, or $E(G)$ is simple but not in $Lie(\ell)$, cf. for instance Theorem \ref{red-so}.
In the case $G$ is a subfield subgroup, it is indeed irreducible on all restricted $\GC$-modules. In
all other cases, it is straightforward to check that $\dim(St)$ is larger than the largest degree of
irreducible $\FF G$-representations, whence $G$ is reducible on $St$. Some exceptions do arise:
$A_{8} \simeq \Omega^{+}_{6}(2) < SO_{6}(\FF)$ and
$PSp_{4}(3) \simeq \Omega^{-}_{6}(2) < SO_{6}(\FF)$ for $\ell = 2$, and
$SU_{4}(2) \simeq \Omega_{5}(3) < SO_{5}(\FF)$ for $\ell = 3$, but they give rise to subfield
subgroups.

5) Finally we consider the case that $\GC$ is classical with natural module $V$ of dimension
$d \leq 4$. Suppose $G$ is proper closed subgroup of $\GC$ that is irreducible on all $X \in \SL$.

The case $(\GC,\ell) = (A_{1},2)$ leads to the possibility (ii). Assume $\GC = A_{1}$,
$\ell > 2$. Using 1) we may assume $G$ is finite, whence $G$ is contained in some $SL_{2}(q)$. Take
$q$ smallest subject to this containment, and assume $G \neq SL_{2}(q)$. Let $M$ be a maximal
subgroup of $SL_{2}(q)$ containing $G$. The possibilities for $M$ are listed in \cite{Kle}. Since
$M$ is irreducible on $St$, it is easy to see that either $M = SL_{2}(q_{0})$ with
$q = q_{0}^{b}$ for some prime $b> 2$, or $\ell = 3$ and $M = 2\AAA_{5}$ (and $q = 9$), or
$M = SL_{2}(r).2$ with $r = q^{1/2}$. The first subcase is impossible by the choice of $q$. In the
second subcase, the only proper subgroup of $M$ that is irreducible on $St$ is
$2\AAA_{4} = SL_{2}(3)$, so the choice of $q$ implies that $G = M = 2\AAA_{5}$, as listed in (iii).
In the third subcase we must have $G = H.2$ with $H \leq SL_{2}(r)$, and $H$ is also irreducible on
$St$ since $\dim(St) = \ell$ is odd. Applying the above argument to $H$, we can show that
$\OL(H) = SL_{2}(s)$ for some $s|r$. Now $\OL(G) = \OL(H)$, so $G$ is a (good)
subfield subgroup, as stated in (i).

Assume $\GC = A_{2}$. Again we may assume by 1) that $G$ is finite, and choose smallest $q$ such that
$G \leq SL_{3}(q)$. If $G = SL_{3}(q)$ then we arrive at (i). Otherwise we can find a maximal
subgroup $M$ of $SL_{3}(q)$ containing $G$. Using the irreducibility of $M$ on $\AC(V)$ and $St$ and
the list of $M$ given in \cite{Kle}, we can show that either $M = SL_{3}(q_{0})$ with
$q = q_{0}^{b}$ for some prime $b \neq 3$, or $\ell = 2$ and $M = 3\AAA_{6}$ (and $q = 4$), or
$M = SL_{3}(r).3$ with $r = q^{1/3}$, or $M = SU_{3}(r)$ with $r = q^{1/2}$. The first subcase is
ruled out by the choice of $q$. In the second subcase, no proper subgroup of $M$ can be irreducible
on $St$, so $G = M = 3\AAA_{6}$, as listed in (iv). In the third subcase we must have $G = H.3$ with
$H \leq SL_{3}(r)$, and $H$ is also irreducible on $\AC(V)$ and $St$. Applying the above argument to
$H$, we can again show that $\OL(G) = SL_{3}(s)$ for some $s|r$, as stated in (i). Consider
the fourth subcase. We may assume that $r$ is smallest subject to the containment $G \leq SU_{3}(r)$.
Inspecting the list of maximal subgroups of $SU_{3}(r)$ as given in \cite{Kle}, we can show that
$\OL(G) = SU_{3}(s)$ for some $s|r$, as listed in (i).

Finally, we assume that $\GC = C_{2}$ and $\ell = 2$. By the virtue of 1), we may assume that $G$ is
finite, irreducible and primitive on $V$, and does not fix any quadratic form on $V$. Again choose
smallest $q$ such that $G \leq Sp_{4}(q)$. Using the list of maximal subgroups of $Sp_{4}(q)$ as
reproduced in \cite{Kle}, we can show that either $G$ is as in (i), or $G \leq \ta B_{2}(q)$.
Now using the list of maximal subgroups of $\ta B_{2}(q)$ as reproduced in \cite{Kle}, we can show
that $G = \ta B_{2}(r)$ for some $r \geq 8$, as stated in (i).

Conversely, all the cases listed in (i) -- (vii) give rise to subgroups irreducible on all
irreducible restricted $\GC$-modules.
\end{proof}

Theorem \ref{restr1} immediately yields:

\begin{corol}\label{restr2}
{\sl Let $\FF$ be an algebraically closed field of characteristic $\ell$, and let $\GC$ be a simple
simply connected algebraic group over $\FF$. Assume $\rank(\GC) > 2$ if $\ell = 2$ and
$\rank(\GC) > 1$ if $\ell = 3$. Let $\SL$ be the set of restricted $\GC$-modules
defined in Theorem \ref{restr1}. Then $\cap_{W \in \SL} \IL_{r}(W) =\IL_{r}(\GC)$.
\hfill $\Box$}
\end{corol}

Taking $r = 2$ and $\LF = \SL$, we get Corollary \ref{nice}.

\medskip
We point out the following related result that requires considerably less effort than the previous
results.

\begin{theor}\label{single}
{\sl Let $\FF$ be an algebraically closed field of characteristic $\ell$. Let $\GC$ be a simple
simply connected algebraic group over $\FF$. There exists an irreducible rational $\FF\GC$-module
$W$ such that no proper closed subgroup other than subfield subgroups are irreducible on $W$.
Moreover, if $\ell > 3$, then we may take $W$ to be of dimension at most $(\dim \GC)^{2}$.}
\end{theor}

\begin{proof}
We give the proof in positive characteristic. The proof is a bit easier in characteristic zero
(but different). Alternatively, one could use the positive characteristic result to deduce the
characteristic zero result.

We first exclude $B_{n}$ and $F_{4}$ in characteristic $2$ and $G_{2}$ in characteristic $3$. Let 
$U$ be the nontrivial irreducible constituent of the adjoint module for ${\GC}$. We have already 
seen that this is not irreducible over any proper closed positive dimensional subgroup. By a result 
of Larsen and Pink \cite{larsenpink} (this also can be derived from the classification
of finite simple groups and representation theory for the classical groups; in particular, Liebeck
and Seitz proved this result for exceptional groups), it follows that every finite subgroup
of $\GC$ is either a subfield subgroup, is contained in a positive dimensional proper closed
subgroup or has order at most $N$ (here the bound $N$ depends only upon the type of $\GC$ and not
on the characteristic). Let $M = \prod_{2 \leq a \leq N, (\ell,a) = 1}\varphi(a)$, where
$\varphi(\cdot)$ is the Euler function, and let $\sigma$ be the field automorphism sending
$x \in \FF$ to $x^{\ell^{M}}$. Then $\sigma$ is trivial on all roots of unity of order at most $N$
in $\FF$.

Let $W = U^{*} \otimes U^{\sigma}$.  Then aside from subfield subgroups, $W$ is reducible. Indeed,
it is so if $H$ is reducible on $U$, in particular it is so if $H$ is a positive dimensional
proper closed subgroup. Assume $H$ is finite of order at most $N$ and $H$ is irreducible on
$U$. The construction of $\sigma$ ensures that the Brauer characters of $H$ on $U^{\sigma}$ and
$U$ are the same, whence $W|_{H} \simeq (U^{*} \otimes U)|_{H}$ and so $H$ has fixed points on $W$.

In the excluded cases, we replace $U$ by a module that is reducible over every positive
dimensional subgroup (there are only finitely many classes of maximal positive dimensional
reductive subgroups -- it is not difficult to choose one for each and take a tensor product of
Frobenius twists of these to obtain such a module). Now argue as above.
\end{proof}

\end{document}